\documentclass[11pt,reqno]{amsart}
\usepackage{amscd}
\usepackage{amsfonts}
\usepackage{amsmath}
\usepackage{amssymb}
\usepackage{amsthm}
\usepackage{caption}
\usepackage{fancyhdr}
\usepackage{latexsym}
\usepackage[colorlinks=true, pdfstartview=FitV, linkcolor=blue,
            citecolor=blue, urlcolor=blue]{hyperref}
\usepackage{enumitem}

\usepackage{mathtools}            
\usepackage{upgreek}
\usepackage{extarrows}
\usepackage{empheq}
\usepackage{cancel}

\usepackage{textcase}

\input epsf
\input texdraw
\input xy
\xyoption{all}

\usepackage{tikz,float}
\usepackage{tikz-cd}

\usetikzlibrary{patterns}
\usetikzlibrary{calc,matrix, arrows}
\usetikzlibrary{snakes}
\usepackage{subfigure}
\usetikzlibrary{decorations.pathreplacing}
\tikzstyle{level 1}=[level distance=3.5cm, sibling distance=8cm]
\tikzstyle{level 2}=[level distance=3.5cm, sibling distance=3cm]
\tikzstyle{level 3}=[level distance=5cm, sibling distance=1.4cm]
\usetikzlibrary{decorations.markings}
\tikzset{-{latex[length=3mm, width=2mm]}-/.style={decoration={
  markings,
  mark=at position #1 with {\arrow{>}}},postaction={decorate}}}
  \tikzset{middlearrow/.style={
        decoration={markings,
            mark= at position 0.55 with {\arrow{#1}} ,
        },
        postaction={decorate}
    }
}

\usepackage{color}

\definecolor{Red}{rgb}{1.00, 0.00, 0.00}
\definecolor{DarkGreen}{rgb}{0.00, 1.00, 0.00}
\definecolor{Blue}{rgb}{0.00, 0.00, 1.00}
\definecolor{Cyan}{rgb}{0.00, 1.00, 1.00}
\definecolor{Magenta}{rgb}{1.00, 0.00, 1.00}
\definecolor{DeepSkyBlue}{rgb}{0.00, 0.75, 1.00}
\definecolor{DarkGreen}{rgb}{0.00, 0.39, 0.00}
\definecolor{dgreen}{RGB}{0,200,100}
\definecolor{ddgreen}{RGB}{0,170,0}
\definecolor{SpringGreen}{rgb}{0.00, 1.00, 0.50}
\definecolor{DarkOrange}{rgb}{1.00, 0.55, 0.00}
\definecolor{OrangeRed}{rgb}{1.00, 0.27, 0.00}
\definecolor{DeepPink}{rgb}{1.00, 0.08, 0.57}
\definecolor{DarkViolet}{rgb}{0.58, 0.00, 0.82}
\definecolor{SaddleBrown}{rgb}{0.54, 0.27, 0.07}
\definecolor{Black}{rgb}{0.00, 0.00, 0.00}
\definecolor{dark-magenta}{rgb}{.5,0,.5}
\definecolor{myblack}{rgb}{0,0,0}
\definecolor{darkgray}{gray}{0.5}
\definecolor{lightgray}{gray}{0.75}

\renewcommand{\b}{\mathcolor{blue}}
\newcommand{\nn}{\nonumber}
\newcommand{\p}{\partial}

\newcommand{\supp}{\text{supp}}
\newcommand{\rr}{\mathbb{R}}
\newcommand{\ci}{\mathbb{T}}
\renewcommand{\p}{\partial}
\newcommand{\zz}{\mathbb{Z}}
\newcommand{\cc}{\mathbb{C}}
\newcommand{\NN}{\mathbb{N}}

\def\nin{\noindent}

\makeatletter
\def\mathcolor#1#{\@mathcolor{#1}}
\def\@mathcolor#1#2#3{%
  \protect\leavevmode
  \begingroup
    \color#1{#2}#3%
  \endgroup
}
\makeatother

\theoremstyle{plain}  
\newtheorem{theorem}{Theorem}[section]
\newtheorem{proposition}{Proposition}[section]
\newtheorem{lemma}{Lemma}[section]

\numberwithin{figure}{section}
\numberwithin{equation}{section}

\usepackage[breakable, theorems, skins]{tcolorbox}
\tcbset{enhanced}

\usepackage[framemethod=TikZ]{mdframed}

\makeatletter
\@namedef{subjclassname@2020}{%
  \textup{2020} Mathematics Subject Classification}
\makeatother

\begin{document}
\title{
The Majda-Biello  system on the half-line
}
\author{A. Alexandrou Himonas \& Fangchi Yan}
\date{July 6, 2022
 \mbox{}$^*$\!\textit{Corresponding author}:
himonas.1@nd.edu}



\begin{abstract}
The Majda-Biello system models the interaction of Rossby waves.
It consists of two coupled KdV equations one of which has a  parameter 
$\alpha$ as coefficient of its dispersion. This work studies this system on the half line with Robin, Neumann, and Dirichlet  boundary data. 
It shows that for $0<\alpha<1$ or $1<\alpha<4$ all these problems are well-posed for initial data in Sobolev spaces $H^s$, $s\ge 0$.
For  $\alpha=1$ or $\alpha>4$  well-posedness holds
for Dirichlet data if $s>-3/4$, while for  Neumann and Robin data it 
depends 
on the sign of the parameters 
involved in the data. 
For $\alpha=4$ well-posedness of all problems holds for $s\ge 3/4$.
The Robin and Neumann boundary data are in $H^{s/3}$ 
while the  Dirichlet  boundary data are in $H^{(s+1)/3}$.
 These are consistent with the time regularity of the Cauchy problem for the corresponding linear system.
The proof is based on linear estimates in Bourgain spaces
derived by utilizing the Fokas solution formula for the forced linear system, and appropriate bilinear estimates suggested by 
the coupled nonlinearities.  These show that the iteration map defined via the Fokas formula is a contraction 
in appropriate solution spaces.
All the well-posedness results obtained here are optimal.

 \end{abstract}

\maketitle

\markboth{The Majda-Biello  system on the half-line}{
A. Himonas and F. Yan}


%
%
%
%
%
%
%
%
\section{Introduction and Results}
\label{Introduction}
\setcounter{equation}{0}

In this work we consider the three  basic 
initial-boundary value problems (ibvp) for the Majda-Biello (MB) system on the half-line.
This system consists of  two coupled 
Korteweg-de Vries (KdV) equations 
with the dispersion of the second
equation having a coefficient $\alpha$,
which in  \cite{mb2003}  is assumed  to be
in the interval $(0, 1)$ and close to 1, although here we consider all 
$\alpha>0$.
We begin with the Robin problem:
\begin{subequations}
\label{Rob-MBsys}
\begin{align}
\label{Rob-MB-eqn1}
&
\p_tu+\p_x^3u+
\p_x(v^2)
=
0,
\quad
x\in\rr^+,
\,\,
t\in(0,T),
\\
\label{Rob-MB-eqn2}
&
\p_tv+\alpha\p_x^3v+
\p_x(uv)
=
0,
\\
\label{Rob-MBsys-ic}
&(u,v)(x,0)
=
(u_0,v_0)(x),
\quad
x\in\rr^+,
\\
\label{Rob-MBsys-bc}
&(u_x+\gamma_1 u ,v_x+\gamma_2v)(0, t) = (\varphi_1,\varphi_2)(t),
\quad
t\in(0,T),
\end{align}
\end{subequations}
where the parameters $\gamma_1$, $\gamma_2$ are
 real numbers. When $\gamma_1=\gamma_2=0$, then this 
becomes 
 the Neumann problem, which is the second ibvp
 we study here. Finally, replacing the Robin boundary
 data \eqref{Rob-MBsys-bc} with the Dirichlet boundary data
$
 (u,v)(0, t) = (g_0,h_0)(t)
 $
 gives the third problem. Here, we shall study the 
 well-posedness of these problems for initial data
 $u_0(x) ,v_0(x)$  in spatial  Sobolev  spaces $H^s(0,\infty)$ and boundary data
 in temporal  Sobolev  spaces suggested by the time 
 regularity of the  Cauchy problem of the 
 corresponding linear MB system, which for the Robin 
 problem \eqref{Rob-MBsys} we have that 
 $\varphi_1(t),\varphi_2(t)$ are in 
  $H^{s/3}(0,T)$.
 Before stating our results about  MB  we provide an
 outline of its physical significance.

System \eqref{Rob-MBsys} was  introduced by   Majda and Biello in  \cite{mb2003} to model  the nonlinear interaction of barotropic and equatorial baroclinic Rossby waves.
These waves occur in the troposphere, which is the
 lowest layer of the atmosphere, largely due to the Earth's rotation, and affect  weather and climate.
 In the system  \eqref{Rob-MBsys} $u$ is the amplitude of a (barotropic) Rossby wave packet with significant energy in the midlatitudes, and $v$ is the amplitude of an equatorially confined (baroclinic) Rossby wave packet. 
  The parameter $\alpha$ is close to 1 with 0.899, 0.960, 0.980,
  being some of its numerical estimates given in  \cite{mb2003}.
   For a complete picture about 
  the physical relevance of the MB system we refer
  the reader to   \cite{mb2003, mb2004, bm2004-1, bm2004-2, b2009}. 

  Concerning 
   conservation laws, it can be seen that the MB system conserves the quantities $\int udx$ and  $\int vdx$. Furthermore, and more importantly,  MB conserves  the total energy \cite{mb2003}
  \begin{equation}
  \label{conserved-energy}
  E=\int \Big[u^2+v^2\Big]dx,
  \end{equation}
  which is useful in the study of its global well-posedness,
  and the Hamiltonian 
  \begin{equation}
 \label{conserved-hamiltonian}
   H=\int \Big[u_x^2+\alpha v_x^2-uv^2\Big]dx,
  \end{equation}
which provides MB with a Hamiltonian structure \cite{bm2004-1}.
  There are no other conservation laws \cite{V2015}.
  Therefore, MB may not be a completely integrable system.
  This should be contrasted with KdV  \cite{b1877, kdv1895},
  which is an integrable equation having infinitely many conservation laws, a Lax pair  \cite{lax1968}, and it can be solved via the inverse scattering transform
   \cite{ggkm1967, zk1965}. We mention that although
  the Fokas unified transform method was motivated by 
  integrable equations \cite{f1997}, our work demonstrates that it can also be applied to other equations and systems, like the MB system whose  ibvp's we are studying here.

Next, we recall the spaces needed for stating our results 
precisely.   For $s\in \rr$ the Sobolev space $H^s(\rr)$ consists 
of all temperate distribution $F$ with finite norm
$
\left\| F \right\|^2_{H^s(\mathbb R)}
\doteq
\int_{\xi\in\mathbb R} \hskip-0.02in (1+ \xi^2)^s
|\widehat F(\xi)|^2 d\xi,
$
where
$\widehat F(\xi)$ is the  Fourier transform defined by 
$
  \widehat F(\xi)
  \doteq
  \int_{\rr} e^{-ix\xi} F(x)dx.
  $
Furthermore, 
for an open interval $(a, b)$ in $\rr$,  the Sobolev space $H^s(a,b)$ is defined by 
\begin{equation}
H^s(a,b) = \Big\{f: f= F\big|_{(a,b)}\
 \mathrm{where}\ F\in H^s(\mathbb R)\
  \mathrm{and}\
\left\| f \right\|_{H^s(a,b)}\doteq \inf_{F \in H^s(\mathbb R)} 
\left\| F \right\|_{H^s(\mathbb R)} <\infty
   \Big\}.
\end{equation}
\vskip-0.05in
\nin
Also, for any real numbers $s$, $b$ and $\alpha$ the Bourgain space $X^{s,b}_{\alpha}(\rr^2)$ associated with the linear part of  KdV
is defined by the norm  \cite{b1993-kdv, kpv1996}
\begin{equation}
\label{bourgain-kdv}
\|u\|_{X^{s,b}_{\alpha}}^2
=
\|u\|_{X^{s,b}_{\alpha}(\rr^2)}^2
\doteq
\int_{-\infty}^\infty\int_{-\infty}^\infty
(1+|\xi|)^{2s}
(1+|\tau-\alpha\xi^3|)^{2b}
|\widehat{u}(\xi,\tau)|^2
d\xi
d\tau.
\end{equation}
Furthermore, we shall need  the  restriction space 
$X^{s,b}_{\alpha,\rr^+\times(0, T)}$,
which is defined as follows
\begin{equation}
\label{sb-restrict}
X^{s,b}_{\alpha,\rr^+\times(0, T)}
\doteq
\{
u:
u(x,t)
=
\tilde{u}(x,t)
\;\;\mbox{on}\;\; \mathbb{R^+}\times(0,T)
\,\,
\text{with}
\,\,
\tilde{u}\in X^{s,b}_{\alpha}(\rr^2)
\},
\end{equation}
and 
which is equipped  with the norm
\begin{equation}
\label{sb-restrict-norm}
\| u\|_{X^{s,b}_{\alpha,\rr^+\times(0, T)}}
\doteq
\inf\limits_{\tilde{u}\in X^{s,b}_{\alpha}}\left\{
\| \tilde{u}\|_{X^{s,b}_{\alpha}}
:
\tilde{u}(x,t)
=
u(x,t)
\;\;\mbox{on}\;\; \mathbb{R^+}\times(0,T)
\right\}.
\end{equation}
Also, we shall need the following modification 
of the Bourgain norm  for the KdV  equation 
\cite{b1993-kdv}
\begin{align}
\label{bourgain-like-kdv}
\|u\|_{X^{s,b,\theta}_{\alpha}}^2
=&
\|u\|_{X^{s,b,\theta}_{\alpha}(\rr^2)}^2
\doteq
\|u\|_{X^{s,b}_{\alpha}(\rr^2)}^2
+
\Big[
\int_{-\infty}^\infty
\int_{-1}^1
(1+|\tau|)^{2\theta}
|\widehat{u}(\xi,\tau)|^2
d\xi
d\tau 
\Big]
\\
\simeq&
\int_{-\infty}^\infty
\int_{-\infty}^\infty
\Big[
(1+|\xi|)^{s}
(1+|\tau-\alpha\xi^3|)^{b}
+
\chi_{|\xi|< 1}(\xi)(1+|\tau|)^{\theta}
\Big]^2
|\widehat{u}(\xi,\tau)|^2
d\xi
d\tau,
\notag
\end{align}
where 
$\theta>\frac 12$ (in fact  $\theta<1$ and close to $\frac 12$),
and $\chi_{|\xi|< 1}$ is the characteristic function of the interval 
$(-1, 1)$.
Finally, we shall need the  modified restriction space 
$X_{\alpha,\rr^+\times(0, T)}^{s,b,\theta}$,
which is defined as follows
\begin{equation}
\label{sb-restrict}
X_{\alpha,\rr^+\times(0, T)}^{s,b,\theta}
\doteq
\{
u:
u(x,t)
=
\tilde{u}(x,t)
\;\;\mbox{on}\;\; \mathbb{R^+}\times(0,T)
\,\,
\text{with}
\,\,
\tilde{u}\in X^{s,b,\theta}_{\alpha}(\rr^2)
\},
\end{equation}
and which is equipped  with the norm
\begin{equation}
\label{sb-restrict-norm}
\|u\|_{X_{\alpha,\rr^+\times(0, T)}^{s,b,\theta}}
\doteq
\inf\limits_{\tilde{u}\in X^{s,b,\theta}_\alpha}\left\{
\| \tilde{u}\|_{X^{s,b,\theta}_{\alpha}}:
\;\tilde{u}(x,t)
=
u(x,t)
\;\;\mbox{on}\;\; \mathbb{R^+}\times(0,T)
\right\}.
\end{equation}
$Y^{s,b}_{\alpha}$ is a ``temporal" Bourgain space defined 
by the norm, used for the KdV ibvp by Faminskii  \cite{Fa2004,Fa2007},
\begin{equation}
\label{Y-def-recall}
\|
u
\|_{Y^{s,b}_{\alpha}}
\doteq
\left(
\int_{\rr^2}
(1+|\tau|)^{\frac{2s}{3}}
(1+|\tau-\alpha\xi^3|)^{2b}
|\widehat{u}(\xi,\tau)|^2
d\xi d\tau
\right)^{1/2}.
\end{equation}
Also, we need the  restriction space 
$Y_{\alpha,\rr^+\times(0, T)}^{s,b}$,
which is defined as follows
\begin{equation}
\label{sb-restrict-Y}
Y_{\alpha,\rr^+\times(0, T)}^{s,b}
\doteq
\{
u:
u(x,t)
=
\tilde{u}(x,t)
\;\;\mbox{on}\;\; \rr^+\times(0,T)
\,\,
\text{with}
\,\,
\tilde{u}\in Y^{s,b}_{\alpha}(\rr^2)
\},
\end{equation}
and which is equipped  with the norm
\begin{equation}
\label{sb-restrict-norm-Y}
\| u\|_{Y_{\alpha,\rr^+\times(0, T)}^{s,b}}
\doteq
\inf\limits_{\tilde{u}\in Y^{s,b}_{\alpha}}\left\{
\| \tilde{u}\|_{Y^{s,b}_{\alpha}}:
\;\tilde{u}(x,t)
=
u(x,t)
\;\;\mbox{on}\;\; \rr^+\times(0,T)
\right\}.
\end{equation}
For $\alpha=1$ we simplify the notation for the spaces defined above
as follows
$$
X^{s,b}
\doteq
X^{s,b}_{1},
\quad
X^{s,b,\theta}
\doteq
X^{s,b,\theta}_{1}
\quad
\text{and}
\quad
Y^{s,b}
\doteq
Y^{s,b}_{1}.
$$
\nin
Finally,  to state our results concisely we need 
the following definition of  critical Sobolev exponent:

\begin{minipage}{0.5\linewidth}
\begin{equation}
\label{MB-critical-sob-index}
s_{c}(\alpha)
\doteq
\begin{cases}
0,
\quad
&0<\alpha<4
\,\,
\text{and}
\,\,
 \alpha\ne 1,
\\
-\frac34^+,
\quad
&\alpha=1
\,\,
\text{or}
\,\,
\alpha>4,
\\
\frac34,
\quad
&\alpha=4.
\end{cases}
\end{equation}
\end{minipage}
\nin
\begin{minipage}{0.5\linewidth}
\begin{center}
\begin{tikzpicture}[yscale=1.3, xscale=1]
%
%
\newcommand\X{0};
\newcommand\Y{0};
\newcommand\FX{8};
\newcommand\FY{8};
\newcommand\R{0.6};
\newcommand*{\TickSize}{2pt};
%
%
\draw[black,line width=0.5pt,-{Latex[black,length=2mm,width=2mm]}]
(0,0)
--
(6,0)
node[above]
{\fontsize{\FX}{\FY} \textcolor{black}{$\alpha$}};

\draw[black,line width=0.5pt,-{Latex[black,length=2mm,width=2mm]}]
(0,-1)
--
(0,1.3)
node[right]
{\fontsize{\FX}{\FY} \textcolor{black}{$s_{c}(\alpha)$}};

\draw[dashed,black,line width=0.5pt]
(0,0)
node[xshift=-0.2cm, yshift=0cm]
{\fontsize{\FX}{\FY} $0$}

(0,-3/4)
node[xshift=-0.25cm, yshift=0cm]
{\fontsize{\FX}{\FY}$-\frac34$}
--
(0.1,-3/4)
(1,-3/4)
node[red,xshift=0cm, yshift=0cm]
{$\bullet$}
--
(1,-0.05)
node[xshift=0cm, yshift=0.3cm]
{\fontsize{\FX}{\FY} $1$}

(0,3/4)
node[xshift=-0.2cm, yshift=0cm]
{\fontsize{\FX}{\FY}\bf $\frac34$}
--
(0.1,3/4)
(4,3/4)
node[red,xshift=0cm, yshift=0cm]
{$\bullet$}
--
(4,0.05)
node[xshift=0cm, yshift=-0.3cm]
{\fontsize{\FX}{\FY}\bf $4$}

(1,0)
node[xshift=0cm, yshift=0cm]
{$\circ$}

(4,0)
node[xshift=0cm, yshift=0cm]
{$\circ$}

(4,-3/4)
node[xshift=0cm, yshift=0cm]
{$\circ$}
;

\draw[red,line width=1.5pt]
(0,0)
--
(0.9,0)

(1.1,0)
--
(3.9,0)

(4.1,-3/4)
--
(5.5,-3/4)
;

\end{tikzpicture}

%
\end{center}
\end{minipage}
\begin{theorem}
[\textcolor{blue}{Well-posedness for the Robin and Neumann ibvp's}]
\label{thm-MB-half-line-Rob}
Let $\alpha>0$. 
\vskip0.03in
\nin
$\bullet$ 
If $\gamma_1\le 0$ and  $\gamma_2\le 0$, 
then for $s_c(\alpha)\le s<\frac32$,
initial data $(u_0, v_0) \in H^s(0,\infty)\times H^s(0,\infty)$,
and  boundary data  
$(\varphi_1,\varphi_2)\in H_t^\frac{s}{3}(0,T)\times H_t^\frac{s}{3}(0,T)$, 
there  is a lifespan  $0 < T_0\leq T<\frac12$ such that  the
MB system ibvp  \eqref{Rob-MBsys} admits a unique solution $(u,v)\in X^{s,b,\theta}_{1,\rr^+\times(0,T_0)}\times X^{s,b,\theta}_{\alpha,\rr^+\times(0,T_0)}$ satisfying the size estimate
\begin{align}
\label{X-norm-def-Rob1}
\|u\|_{X^{s,b,\theta}_{1,\rr^+\times(0,T_0)}}
\hskip-0.1in
+
\|v\|_{X^{s,b,\theta}_{\alpha,\rr^+\times(0,T_0)}}
\hskip-0.05in
\lesssim
\Big[
\|
u_0
\|_{H^{s}(\rr^+)}
+
\|
v_0
\|_{H^{s}(\rr^+)}
+
\|\varphi_1\|_{H_t^\frac{s}{3}(0,T)}
+
\|\varphi_2\|_{H_t^\frac{s}{3}(0,T)}
\Big],
\end{align}
for some $b\in(0,\frac12)$ and $\theta\in (\frac12,1)$. Also, an estimate for the lifespan is given by
\begin{align}
\label{MB-lifespan-Rob1}
T_0
=c_0
\Big[
1
+
\|
u_0
\|_{H^{s}(\rr^+)}
+
\|
v_0
\|_{H^{s}(\rr^+)}
+
\|\varphi_1\|_{H_t^\frac{s}{3}(0,T)}
+
\|\varphi_2\|_{H_t^\frac{s}{3}(0,T)}
\Big]^{-4/\beta},
\end{align}
where   
$\beta$ is given by \eqref{beta-choice}.
Finally, the data to solution map  $\{u_0, v_0, \varphi_1, \varphi_2\}\mapsto (u,v)$ is locally Lipschitz continuous.
\vskip0.05in
\nin
$\bullet$  If $\gamma_1>0 $ or  $\gamma_2> 0$, 
then  the above conclusions  hold for
$\max\{0, s_c(\alpha)\}\le s<\frac32$.
\end{theorem}
Well-posedness results for the Robin and Neumann problems
of the  KdV and NLS equations  with smooth data were obtained in 
 \cite{hmy2021, hm2021}. Furthermore, these problems 
 for KdV and rough data were studied in  \cite{hy2022-Robin}.
The common ingredient in both situations was the 
employment of the Fokas method for solving 
the corresponding forced linear ibvp.

%

%
%
%
%
%
%
%
%
%
%
%
\vskip0.05in
\nin
\textbf{\large The Dirichlet ibvp.}
Concerning  the Dirichlet problem for the MB system
\begin{subequations}
\label{KdV-sys}
\begin{align}
\label{KdV-eqn1}
&
\p_tu+\p_x^3u+
\p_x(v^2)
=
0,
\quad
x\in\rr^+,
\,\,
t\in(0,T),
\\
\label{KdV-eqn2}
&
\p_tv+\alpha\p_x^3v+
\p_x(uv)
=
0,
\\
\label{KdV-sys-ic}
&(u,v)(x,0)
=
(u_0,v_0)(x),
\quad
x\in\rr^+,
\\
\label{KdV-sys-bc}
&(u,v)(0, t) = (g_0,h_0)(t),
\quad
\,
t\in(0,T),
\end{align}
\end{subequations}
 we shall need the following compatibility conditions
\begin{equation}
\label{comp-cond}
u_0(0)
=
g_0(0),
\quad
\text{and}
\quad
v_0(0)
=
h_0(0),
\qquad
\frac12<s<\frac32.
\end{equation}
Our local well-posedness result for  this problem is the following, which is stated concisely by using the critical Sobolev exponent $s_c(\alpha)$ defined in  \eqref{MB-critical-sob-index}.
\begin{theorem}
[\textcolor{blue}{Well-posedness for MB Dirichlet ibvp}]
\label{thm-MB-half-line-Dir}
Let $\alpha>0$. Then, for $s_c(\alpha)\le s<\frac32$, $s\ne \frac 12$, initial data $(u_0, v_0) \in H^s(0,\infty)\times H^s(0,\infty)$
with the  compatibility 
conditions
\eqref{comp-cond}, 
and  boundary data  $(g_0,h_0)\in H_t^\frac{s+1}{3}(0,T)\times H_t^\frac{s+1}{3}(0,T)$, 
there  is a lifespan  $0 < T_0\leq T<\frac12$ such that  
Majda-Biello system ibvp \eqref{KdV-sys}
admits a unique solution $(u,v)\in X^{s,b,\theta}_{1,\rr^+\times(0,T_0)}\times X^{s,b,\theta}_{\alpha,\rr^+\times(0,T_0)}$ satisfying the size estimate
\begin{align}
\label{X-norm-def-Dir1}
\|u\|_{X^{s,b,\theta}_{1,\rr^+\times(0,T_0)}}
\hskip-0.1in
+
\|v\|_{X^{s,b,\theta}_{\alpha,\rr^+\times(0,T_0)}}
\hskip-0.05in
\lesssim
\Big[
\|
u_0
\|_{H^{s}(\rr^+)}
+
\|
v_0
\|_{H^{s}(\rr^+)}
+
\|g_0\|_{H_t^\frac{s+1}{3}(0,T)}
+
\|h_0\|_{H_t^\frac{s+1}{3}(0,T)}
\Big],
\end{align}
for some $b\in(0,\frac12)$ and $\theta\in (\frac12,1)$. 
Also, it satisfies the lifespan estimate given by
\begin{align}
\label{MB-lifespan-Dir1}
T_0
=c_0
\Big[
1
+
\|
u_0
\|_{H^{s}(\rr^+)}
+
\|
v_0
\|_{H^{s}(\rr^+)}
+
\|g_0\|_{H_t^\frac{s+1}{3}(0,T)}
+
\|h_0\|_{H_t^\frac{s+1}{3}(0,T)}
\Big]^{-4/\beta},
\end{align}
where   
$\beta$  is given by \eqref{beta-choice}.
Moreover,  the data to solution map  $\{u_0, v_0, g_0, h_0\}\mapsto (u,v)$ is locally Lipschitz continuous.
\end{theorem}
One can show that the Dirichlet problem for the MB system 
is well-posed  on the half-line for all $s\ge s_c(\alpha)$ 
except for
$\frac{s+1}{3} \in  \NN_0 +\frac 12$, 
and with certain additional compatibility conditions.
(For the  Robin and Neumann problems the excluded values are
$\frac{s}{3}\in  \NN_0 +\frac 12$.) 
We mention that the well-posedness of the MB system
with Dirichlet data for $0<\alpha<1$, $s\ne \frac 12,  \frac 32$, has been studied in \cite{e2020} using the Laplace transform
in the time variable for solving the forced linear ibvp.
For $s>-3/4$, the analogous to Theorem \ref{thm-MB-half-line-Dir} result for KdV has been proved in \cite{hy2022-KdVm}.
For smooth data the Dirichlet problem on the half-line for the
KdV, NLS, and the ``good" Boussinesq equations with smooth
data was studied in \cite{fhm2016, fhm2017, hm2015}.
In both situations the basic ingredient in the proofs was
the use of the Fokas unified transform method 
for producing the solution formulas of the corresponding linear ibvp.
We recall that there are two other approaches for studying ibvp's.
The first method, developed by 
Bona, Sun and Zhang  \cite{bsz2002,bsz2003, bsz2018},
uses the Laplace transform
in the time variable 
for solving the  related 
to the equation forced linear ibvp.
This method in combination with a bilinear smoothing 
is also used by Erdogan and Tzirakis \cite{et2016}.
The second method, developed 
by Colliander and Kenig  \cite{ck2002}, and by Holmer  \cite{h2005, h2006}, expresses an ibvp as a superposition of initial value problems. 

%
%
%
%
%
%
%
%
%
%
%
\vskip0.05in
\nin
\textbf{\large The Linear Neumann and Robin Problems.}
In order to prove the well-posedness result for the 
MB system ibvp  \eqref{Rob-MBsys} we first 
solve the corresponding forced  linear ibvp for the $u$-equation
\begin{subequations}
\label{Rob-LKdV}
\begin{align}
\label{Rob-LKdV-eqn}
&\p_tu+\p_x^3u=f_1(x,t),
\quad
0<x<\infty,
\,\,
0<t<T,
\\
\label{Rob-LKdV-ic}
&u(x,0)
=
u_0(x)\in H_x^s(0,\infty),
\\
\label{Rob-LKdV-bc}
&u_x(0, t)+\gamma_1 u(0,t)
=
\varphi_1(t)
\in H_t^\frac{s}{3}(0,T),
\end{align}
\end{subequations}
and also for the $v$-equation with a different forcing 
\begin{subequations}
\label{Rob-LKdV-v}
\begin{align}
\label{Rob-LKdV-eqn-v}
&\p_tv+\alpha\p_x^3v=f_2(x,t),
\quad
0<x<\infty,
\,\,
0<t<T,
\\
\label{Rob-LKdV-ic-v}
&v(x,0)
=
v_0(x)\in H_x^s(0,\infty),
\\
\label{Rob-LKdV-bc-v}
&v_x(0, t)+\gamma_2 v(0,t)
=
\varphi_2(t)
\in H_t^\frac{s}{3}(0,T).
\end{align}
\end{subequations}
Using the Fokas  unified transform method (UTM)
we get the following {\bf UTM solution formulas}
for $u$ and $v$ (see \cite{hmy2021} for a detailed derivation):
\begin{align}
\label{Rob-kdv-utm}
&u(x,t)
=
S_1[u_0,\varphi_1;f_1]
\\
&\doteq
\frac{1}{2\pi}
\int_{-\infty}^\infty
 e^{i\xi x+i\xi^3 t} 
 \big[
\widehat{u}_0(\xi )
  +
 F_1(\xi , t)  
 \big]   d\xi 
 \nonumber
\\
&+
\frac{1}{2\pi}
\int_{\partial D^+}
e^{i\xi x+i\xi ^3 t}  
\frac{
(\xi +i(\sigma+1)\gamma_1)
[
F_1(\sigma \xi , T)
+
\widehat{u}_0(\sigma \xi )] 
-
(\xi (\sigma+1)+i\gamma_1)
[
F_1(\sigma^2 \xi , T)
+
\widehat{u}_0(\sigma^2 \xi )] 
}{
\sigma(\xi -i\gamma_1)
}
d\xi 
\nonumber
\\
&-
\chi_{(0,\infty)}(\gamma_1)
\frac{(2+\sigma)\gamma_1}{\sigma}
\Big[
\widehat{u}_0(i\sigma^2\gamma_1)
+
F_1(i\sigma^2\gamma_1, T)
-
\widehat{u}_0(i\sigma\gamma_1)
-
F_1(i\sigma\gamma_1, T)
\Big]
e^{-\gamma_1 x+\gamma_1^3t},
\qquad        
\nonumber
\\
&+
\frac{3i}{2\pi}
\int_{\partial D^+}
e^{i\xi x+i\xi ^3 t}  
\frac{\xi ^2}{\xi -i\gamma_1}
\widetilde{\varphi}_1(\xi ^3, T)
d\xi
\,\,-\,\,
3\chi_{(0,\infty)}(\gamma_1)
\gamma_1^2
\widetilde{\varphi}_1(-i\gamma_1^3, T)
e^{-\gamma_1 x+\gamma_1^3t},
\nonumber
\end{align}
\begin{align}
\label{Rob-kdv-utm-v}
&v(x,t)
=
S_\alpha[v_0,\varphi_2;f_2]
\\
&\doteq
\frac{1}{2\pi}
\int_{-\infty}^\infty
 e^{i\xi x+i\alpha\xi^3 t} 
 \big[
\widehat{v}_0(\xi )
  +
 F_2(\xi , t)  
 \big]   d\xi 
 \nonumber
\\
&+
\frac{1}{2\pi}
\int_{\partial D^+}
e^{i\xi x+i\alpha\xi ^3 t}  
\frac{
(\xi +i(\sigma+1)\gamma_2)
[
F_2(\sigma \xi , T)
+
\widehat{v}_0(\sigma \xi )] 
-
(\xi (\sigma+1)+i\gamma_2)
[
F_2(\sigma^2 \xi , T)
+
\widehat{v}_0(\sigma^2 \xi )] 
}{
\sigma(\xi -i\gamma_2)
}
d\xi 
\nonumber
\\
&-
\chi_{(0,\infty)}(\gamma_2)
\frac{(2+\sigma)\gamma_2}{\sigma}
\Big[
\widehat{v}_0(i\sigma^2\gamma_2)
+
F_2(i\sigma^2\gamma_2, T)
-
\widehat{v}_0(i\sigma\gamma_2)
-
F_2(i\sigma\gamma_2, T)
\Big]
e^{-\gamma_2 x+\alpha\gamma_2^3t},
\qquad        
\nonumber
\\
&+
\alpha
\frac{3i}{2\pi}
\int_{\partial D^+}
e^{i\xi x+i\alpha\xi ^3 t}  
\frac{\xi ^2}{\xi -i\gamma_2}
\widetilde{\varphi}_2(\xi ^3, T)
d\xi
\,\,-\,\,
\alpha
3\chi_{(0,\infty)}(\gamma_2)
\gamma_2^2
\widetilde{\varphi}_2(-i\gamma_2^3, T)
e^{-\gamma_2 x+\alpha\gamma_2^3t},
\nonumber
\end{align}
where
%
\vskip-0.3in
\begin{align}
\label{sigma-def}
&\hskip2in
\sigma\doteq
e^{\frac{2\pi}{3}i},
\\
\label{half-line-FT-ic}
&\hskip0.3in
\widehat {u}_0(\xi)
\doteq
\int_0^\infty
e^{-ix\xi}
u_0(x)
dx,
\qquad
\widehat {v}_0(\xi)
\doteq
\int_0^\infty
e^{-ix\xi}
v_0(x)
dx,
\,\,\,
\text{ Im} \xi\le 0,
\\
\label{time-trans-bc}
&\hskip0.3in
\widetilde \varphi_1(\xi,T)
\doteq
\int_0^Te^{-i\xi t}\varphi_1(t)dt,
\qquad
\widetilde \varphi_2(\xi,T)
\doteq
\int_0^Te^{
-i\alpha\xi t
}\varphi_2(t)dt,
\\
\label{time-trans-F1}
&F_1(\xi,t)
\doteq
\int_0^t e^{-i\xi^3\tau}\widehat f_1(\xi,\tau)d\tau
=
\int_0^t e^{-i\xi^3\tau}\int_0^\infty e^{-i\xi x}f_1(x,\tau)dxd\tau,
\quad
t\in (0,T),
\\
\label{time-trans-F2}
&F_2(\xi,t)
\doteq
\int_0^t e^{
-i\alpha\xi^3\tau
}
\widehat f_2(\xi,\tau)d\tau
=
\int_0^t e^{
-i\alpha\xi^3\tau
}\int_0^\infty e^{-i\xi x}f_2(x,\tau)dxd\tau,
\quad
t\in (0,T),
\end{align}
and $D^+$ is the domain
in the complex $\xi$-plane  displayed in 
Figure \ref{domain-D-pos-Rob} when $\gamma_j>0$
(making  $i\gamma_j$  a singularity),
and  Figure  \ref{domain-D-neg-Rob} when $\gamma_j\le 0$.

\hskip-0.3in
\begin{minipage}{0.56\linewidth}
\begin{center}
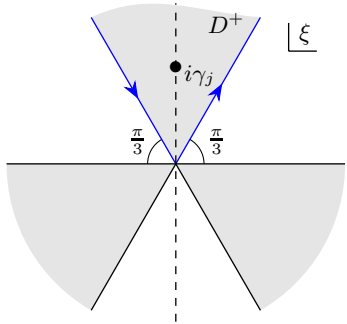

\begin{tikzpicture}[scale=0.75]
%
%
\newcommand\X{0};
\newcommand\Y{0};
\newcommand\FX{9};
\newcommand\FY{9};
\newcommand\R{0.6};
%
%
%
%
\filldraw [fill=gray, fill opacity=0.2, draw opacity=0.01] (0,0)--(-1.5,{3/2*sqrt(3)}) to[out=25,in=175] (1.5,{3/2*sqrt(3)})--cycle;
\filldraw [fill=gray, fill opacity=0.2, draw opacity=0.01] (0,0)--(-1.5,{-3/2*sqrt(3)}) arc (240:180:3cm)--cycle;
\filldraw [fill=gray, fill opacity=0.2, draw opacity=0.01] (0,0)--(1.5,{-3/2*sqrt(3)}) arc (-60:0:3cm)--cycle;
\draw[middlearrow={Stealth[scale=1.3]},blue,line width=0.5pt] (0,0) -- (60:3);
\draw[middlearrow={Stealth[scale=1.3,reversed]},blue,line width=0.5pt] (0,0) -- (120:3);

\draw[black,line width=0.5pt] (0,0) -- (180:3);
\draw[black,line width=0.5pt] (0,0) -- (0:3);

\draw[black,line width=0.5pt] (0,0) -- (-60:3);
\draw[black,line width=0.5pt] (0,0) -- (-120:3);

\draw (0.5,0)
node[xshift=0.1cm,yshift=0.3cm]
{\fontsize{\FX}{\FY} $\frac{\pi}{3}$}
arc (0:60:0.5);

\draw (-0.5,0)
node[xshift=-0.2 cm,yshift=0.3cm]
{\fontsize{\FX}{\FY} $\frac{\pi}{3}$}
arc (180:120:0.5);

\draw [] (0.8,2.5) node {\it\fontsize{\FX}{\FY} $D^+$};
\draw [] (0,1.7) node {$\bullet$};
\draw [] (0,1.7) node[yshift=-0.1cm,xshift=0.25cm] {\it\fontsize{\FX}{\FY} $i\gamma_j$};

\draw[dashed, line width=0.5pt]
(0,-2.8)
--
(0,2.85)
;

\draw[line width=0.5pt, black] 
(2,2.5)--(2,2) (2.5,2)--(2,2);
\node[] at ({2.2},{2.3}) {\fontsize{\FX}{\FY} $\xi$};

\end{tikzpicture}
\end{center}

\vskip-0.15in
\captionof{figure}{The domain $D^+$ with $\gamma_j>0$}
\label{domain-D-pos-Rob}
\end{minipage}
\hskip-0.3in
\begin{minipage}{0.56\linewidth}
\begin{center}
\begin{tikzpicture}[scale=0.75]
%
%
\newcommand\X{0};
\newcommand\Y{0};
\newcommand\FX{9};
\newcommand\FY{9};
\newcommand\R{0.6};
%
%
%
%
\filldraw [fill=gray, fill opacity=0.2, draw opacity=0.01] (0,0)--(-1.5,{3/2*sqrt(3)}) to[out=25,in=175] (1.5,{3/2*sqrt(3)})--cycle;
\filldraw [fill=gray, fill opacity=0.2, draw opacity=0.01] (0,0)--(-1.5,{-3/2*sqrt(3)}) arc (240:180:3cm)--cycle;
\filldraw [fill=gray, fill opacity=0.2, draw opacity=0.01] (0,0)--(1.5,{-3/2*sqrt(3)}) arc (-60:0:3cm)--cycle;
\draw[middlearrow={Stealth[scale=1.3]},blue,line width=0.5pt] (0,0) -- (60:3);
\draw[middlearrow={Stealth[scale=1.3,reversed]},blue,line width=0.5pt] (0,0) -- (120:3);

\draw[black,line width=0.5pt] (0,0) -- (180:3);
\draw[black,line width=0.5pt] (0,0) -- (0:3);

\draw[black,line width=0.5pt] (0,0) -- (-60:3);
\draw[black,line width=0.5pt] (0,0) -- (-120:3);

\draw (0.5,0)
node[xshift=0.1cm,yshift=0.3cm]
{\fontsize{\FX}{\FY} $\frac{\pi}{3}$}
arc (0:60:0.5);

\draw (-0.5,0)
node[xshift=-0.2 cm,yshift=0.3cm]
{\fontsize{\FX}{\FY} $\frac{\pi}{3}$}
arc (180:120:0.5);

\draw [] (0.8,2.5) node {\it\fontsize{\FX}{\FY} $D^+$};
\draw [] (0,-1.5) node {$\bullet$};
\draw [] (0,-1.5) node[yshift=-0.25cm,xshift=0.2cm] {\it\fontsize{\FX}{\FY} $i\gamma_j$};

\draw [] (0,0) node {$\bullet$};

\draw[dashed, line width=0.5pt]
(0,-2.8)
--
(0,2.85)
;

\draw[line width=0.5pt, black] 
(2,2.5)--(2,2) (2.5,2)--(2,2);
\node[] at ({2.2},{2.3}) {\fontsize{\FX}{\FY} $\xi$};

\end{tikzpicture}
\end{center}

\vskip-0.15in
\captionof{figure}{The domain $D^+$
with $\gamma_j \le0$}
\label{domain-D-neg-Rob}
\end{minipage}
%

%
%
%
%
%
%
\vskip.05in
\nin
The second step of our approach is to estimate 
the UTM solution formula \eqref{Rob-kdv-utm} and \eqref{Rob-kdv-utm-v} in Bourgain spaces. Doing so we derive 
the linear estimates contained in the following result.
\begin{theorem}
[\textcolor{blue}{Linear estimates for the Robin and Neumann ibvp}]
\label{thm-MB-forced-linear-Rob}
Let $\alpha>0$. Also, assume that $0<b<1/2$ and 
$1/2<\theta
<
\min\{
\frac{1}{2}
+\frac 13 (s+\frac32),1
\}$.
 Then, the Fokas formula  
\eqref{Rob-kdv-utm} defines a solution  $u$  to
the ibvp \eqref{Rob-LKdV}  which  is  in the space
$ X^{s,b,\theta}_{1,\rr^+\times (0,T)}$ and  
satisfies
the estimates
\begin{align}
\label{MB-forced-linear-est1-Rob}
&\|S_1\big[u_0,\varphi_1;f_1\big]\|_{X_{1,\rr^+\times(0,T)}^{s,b,\theta}}
\le
c\big[
\|u_0\|_{H_x^{s}(\rr^+)}
+
\|\varphi_1\|_{H_t^\frac{s}{3}(0,T)}
+
\|f_1\|_{X^{s,-b,\theta-1}_{1,\rr^+\times(0,T)}}
+
\|f_1\|_{Y^{s,-b}_{1,\rr^+\times(0,T)}}
\big],
\end{align}
where
$-\frac32<s<\frac32$ if   $\gamma_1\le 0$,
and
$0\le s<\frac32$  if  $\gamma_1>0$.
Also, 
the Fokas formula  
 \eqref{Rob-kdv-utm-v} define a solution  $v$  to
the ibvp  \eqref{Rob-LKdV-v}
 in the space $X^{s,b,\theta}_{\alpha,\rr^+\times (0,T)}$ satisfying the estimates
\begin{align}
\label{MB-forced-linear-est2-Rob}
&\|S_\alpha\big[v_0,\varphi_2;f_2\big]\|_{X_{\alpha,\rr^+\times(0,T)}^{s,b,\theta}}
\hskip-0.05in
\le
c\big[
\|v_0\|_{H_x^{s}(\rr^+)}
+
\|\varphi_2\|_{H_t^\frac{s}{3}(0,T)}
+
\|f_2\|_{X^{s,-b,\theta-1}_{\alpha,\rr^+\times(0,T)}}
+
\|f_2\|_{Y^{s,-b}_{\alpha,\rr^+\times(0,T)}}
\big],
\end{align}
where
$-\frac32<s<\frac32$ if   $\gamma_2\le 0$,
and
$0\le s<\frac32$  if  $\gamma_2>0$.
\end{theorem}

%
%
%
%
%
%
%
%
%
%
%
\noindent
{\bf From Robin to  Neumann.} 
 If we let $\gamma_1=\gamma_2=0$ in the Robin problem \eqref{Rob-MBsys}, then  we get the 
 Neumann problem for the MB system,
 and our 
 well-posedness
Theorem \ref{thm-MB-half-line-Rob} 
 becomes a result for the Neumann problem. 
 Also, letting $\gamma_2=0$ in the UTM solution formula \eqref{Rob-kdv-utm-v}, 
 using the facts that 
$
\frac1\sigma
=
\sigma^2
$
and
$
\sigma+1
=
-\sigma^2,
$
we get the solution to the Neumann problem 
(for the  linear $v$-equation involving $\alpha$)
\begin{align*}
v(x,t)
&=
\frac{1}{2\pi}
\int_{-\infty}^\infty
 e^{i\xi x+i\alpha\xi^3 t} 
 \big[
\widehat{v}_0(\xi)
  +
 F_2(\xi, t)  
 \big]   d\xi
\\
&+
\frac{1}{2\pi}
\int_{\partial D^+}
e^{i\xi x+i\alpha\xi^3 t} 
\left\{ 
\sigma^2
[
F_2(\sigma \xi, T)
+
\widehat{v}_0(\sigma \xi)] 
+\sigma
[
F_2(\sigma^2 \xi, T)
+
\widehat{v}_0(\sigma^2 \xi)] 
\right\}
d\xi
\\
&+
\alpha
\frac{3i}{2\pi}
\int_{\partial D^+}
e^{i\xi x+i\alpha\xi^3 t}  
\xi
\widetilde{\varphi}_2(\xi^3, T)
d\xi,
\hskip0.2in
\varphi_2(t)=v_x(0,t),
\,\,
\text{when }
\gamma_2=0.
\end{align*}

%
%
%
%
%
%
%
%
%
%
%
\vskip0.05in
\nin
\textbf{\large  The Linear Dirichlet Problem.}
As in the case of the Robin problem, the first step in our approach is to solve the corresponding forced Dirichlet problem for the linear MB system, which decouples to the $u$-equation problem
with forcing $f_1$
\begin{subequations}
\label{LKdV}
\begin{align}
\label{LKdV-eqn}
&\p_tu+\p_x^3u=f_1(x,t),
\quad
0<x<\infty,
\,\,
0<t<T,
\\
\label{LKdV-ic}
&u(x,0)
=
u_0(x)\in H_x^s(0,\infty),
\\
\label{LKdV-bc}
&u(0, t)
=
g_0(t)
\in H_t^\frac{s+1}{3}(0,T),
\end{align}
\end{subequations}
and to the $v$-equation Dirichlet problem
with  forcing $f_2$ 
\begin{subequations}
\label{LKdV-v}
\begin{align}
\label{LKdV-eqn-v}
&\p_tv+\alpha\p_x^3v=f_2(x,t),
\quad
0<x<\infty,
\,\,
0<t<T,
\\
\label{LKdV-ic-v}
&v(x,0)
=
v_0(x)\in H_x^s(0,\infty),
\\
\label{LKdV-bc-v}
&v(0, t)
=
h_0(t)
\in H_t^\frac{s+1}{3}(0,T).
\end{align}
\end{subequations}
Again, employing the Fokas method
 we obtain the following {\bf UTM solution formulas}
(see Section \ref{sec:utm-deriv} for an outline of the derivation):
\begin{align}
\label{kdv-utm}
&u(x,t)
=
S_1\big[u_0,g_0;f_1\big](x, t)
\doteq
\frac{1}{2\pi}\int_{-\infty}^\infty e^{i\xi x+i\xi^3t}
[\widehat {u}_0(\xi)+F_1(\xi,t)]d\xi
\\
+&
\frac{1}{2\pi}\int_{\p D^+} e^{i\xi x+i\xi^3t}
\big\{
\sigma[\widehat u_0(\sigma \xi)
+
F_1(\sigma \xi,t)]+\sigma^2[\widehat u_0(\sigma^2 \xi)+F_1(\sigma^2 \xi,t)]
-
3\xi^2 \tilde g_0(\xi^3,T)
\big\}d\xi,
\nonumber
\end{align}
\begin{align}
\label{kdv-utm-v}
&v(x,t)
=
S_{\alpha}\big[v_0,h_0;f_2\big](x, t)
\doteq
\frac{1}{2\pi}\int_{-\infty}^\infty e^{i\xi x+
i\alpha\xi^3t
}
[\widehat {v}_0(\xi)+F_2(\xi,t)]d\xi
\\
+&
\frac{1}{2\pi}\int_{\p D^+} 
e^{i\xi x+
i\alpha\xi^3t
}
\big\{
\sigma[\widehat v_0(\sigma \xi)
+
F_2(\sigma \xi,t)]+\sigma^2[\widehat v_0(\sigma^2 \xi)+F_2(\sigma^2 \xi,t)]
-
3
\alpha
\xi^2 \tilde h_0(\xi^3,T)
\big\}d\xi,
\nonumber
\end{align}
where $\sigma, \widehat{u}_0, \widehat{v}_0, F_1, F_2$ are given in \eqref{sigma-def}, \eqref{half-line-FT-ic}, \eqref{time-trans-F1}, \eqref{time-trans-F2}
respectively and
\begin{equation}
\label{Rob-time-trans-bc}
\tilde g_0(\xi^3,T)
\doteq
\int_0^Te^{-i\xi^3 t}
g_0(t)
dt,
\qquad
\tilde h_0(\xi^3,T)
\doteq
\int_0^Te^{
-i\alpha\xi^3 t
}
h_0
(t)dt.
\end{equation}
\begin{minipage}{1\linewidth}
\begin{center}
\begin{tikzpicture}[scale=0.75]
%
%
\newcommand\X{0};
\newcommand\Y{0};
\newcommand\FX{9};
\newcommand\FY{9};
\newcommand\R{0.6};
%
%
%
%

\filldraw [fill=gray, fill opacity=0.2, draw opacity=0.01,variable=\x,domain={pi/3}:{2*pi/3}] (0,0)--
plot({3*cos(deg(\x))},{3*sin(deg(\x))})
--cycle;

\draw[middlearrow={Stealth[scale=1.3]},blue,line width=0.5pt] (0,0) -- (60:3);
\draw[middlearrow={Stealth[scale=1.3,reversed]},blue,line width=0.5pt] (0,0) -- (120:3);

\draw[black,line width=0.5pt
] (0,0) -- (180:3);
\draw[
black,line width=0.5pt
] (0,0) -- (0:3);

\draw (0.5,0)
node[xshift=0.1cm,yshift=0.2cm]
{\fontsize{\FX}{\FY} $\frac{\pi}{3}$}
arc (0:60:0.5);

\draw (-0.5,0)
node[xshift=-0.2 cm,yshift=0.2cm]
{\fontsize{\FX}{\FY} $\frac{\pi}{3}$}
arc (180:120:0.5);

\draw [] (0,1.5) node[xshift=0.5cm,yshift=0.8cm] {\fontsize{\FX}{\FY} $D^+$};

\draw[line width=0.5pt, black] 
(2,3)--(2,2.5) (2.5,2.5)--(2,2.5);
\node[] at ({2.2},{2.8}) {\fontsize{\FX}{\FY} $\xi$};

\draw[dashed, line width=0.5pt,
-{latex[length=3mm, width=2mm]}
]
(0,0)
--
(0,3.5)

;

\end{tikzpicture}
\end{center}
\vskip-0.15in
\captionof{figure}{Domain $D^+$}
\label{kdv-domain}
\end{minipage}

\vskip0.05in
Next, we state the linear estimates for the Dirichlet Problem.
\begin{theorem}
[\textcolor{blue}{Linear estimates for  Dirichlet  ibvp}]
\label{thm-MB-forced-linear-Dir}
Let $\alpha>0$.
Suppose that
$-\frac32< s<\frac32$, $s\neq\frac12$,
$0<b<1/2$ and 
$1/2<\theta
<
\min\{
\frac{1}{2}
+\frac 13 (s+\frac32),1
\}$.
 Then, under the compatibility conditions \eqref{comp-cond}, the Fokas formula  
\eqref{kdv-utm} and \eqref{kdv-utm-v} define a solution  $u,v$  to
the ibvp \eqref{LKdV} and \eqref{LKdV-v} respectively, 
which  are in the space
$ X^{s,b,\theta}_{1,\rr^+\times (0,T)}$ and $X^{s,b,\theta}_{\alpha,\rr^+\times (0,T)}$, satisfying the estimates
\begin{align}
\label{MB-forced-linear-est1-Dir}
&\|S_{1}\big[u_0,g_0;f_1\big]\|_{X_{1,\rr^+\times(0,T)}^{s,b,\theta}}
\le
c\big[
\|u_0\|_{H_x^{s}(\rr^+)}
+
\|g_0\|_{H_t^\frac{s+1}{3}(0,T)}
+
\|f_1\|_{X^{s,-b,\theta-1}_{1,\rr^+\times(0,T)}}
+
\|f_1\|_{Y^{s,-b}_{1,\rr^+\times(0,T)}}
\big],
\\
\label{MB-forced-linear-est2-Dir}
&\|S_{\alpha}\big[v_0,h_0;f_2\big]\|_{X_{\alpha,\rr^+\times(0,T)}^{s,b,\theta}}
\le
c\big[
\|v_0\|_{H_x^{s}(\rr^+)}
+
\|h_0\|_{H_t^\frac{s+1}{3}(0,T)}
+
\|f_2\|_{X^{s,-b,\theta-1}_{\alpha,\rr^+\times(0,T)}}
+
\|f_2\|_{Y^{s,-b}_{\alpha,\rr^+\times(0,T)}}
\big].
\end{align}
\end{theorem}

The third and final step of our method for proving 
Theorems \ref{thm-MB-half-line-Rob} and \ref{thm-MB-half-line-Dir}
consists of deriving appropriate bilinear estimates for the 
coupled nonlinearities, 
considering that the iteration map defined 
by the UTM formulas
 \eqref{Rob-kdv-utm}, \eqref{Rob-kdv-utm-v}  and
 \eqref{kdv-utm},  \eqref{kdv-utm-v}
of the corresponding
forced linear MB system with  forcing terms $f_1$ and $f_2$  
are replaced by the nonlinearities $\p_x(v^2)$
and $\p_x(uv)$. Looking at the linear estimates 
\eqref{MB-forced-linear-est1-Rob},
\eqref{MB-forced-linear-est2-Rob}
for the Robin problem, 
and the linear estimates
\eqref{MB-forced-linear-est1-Dir},
 \eqref{MB-forced-linear-est2-Dir}
  for  Dirichlet problem, 
we see that in order to obtain a contraction map we must 
estimate the spatial modified Bourgain norms 
$\|\p_x(v^2)\|_{X^{s,-b,\theta-1}_{1,\rr^+\times(0,T)}}$,
$\|\p_x(uv)\|_{X^{s,-b,\theta-1}_{\alpha,\rr^+\times(0,T)}}$ 
and the temporal Bourgain norms 
$\|\p_x(v^2)\|_{Y^{s,-b}_{1}}$,  $\|\p_x(uv)\|_{Y^{s,-b}_{\alpha}}$
appropriately.  All these estimates are included in our next result,
which is stated concisely by using 
critical exponent $s_c(\alpha)$ defined in  \eqref{MB-critical-sob-index}.
\begin{theorem}
[\textcolor{blue}{Bilinear estimates}]
\label{mb-bilinear-est}
Let $\alpha>0$. 
If $s\ge s_c(\alpha)$, then for any  
 $\frac12-\beta_1\le b'\le b<\frac12<\theta'\le \theta\le \frac12+\beta_1$, where $\beta_1=\beta_1(s)$ is given below,  
\begin{equation}
\label{beta1-choice}
\beta_1
\doteq
\begin{cases}
\frac1{36},
\quad
&
s\ge 0,
\\
\frac1{96}[s+\frac34],
\quad
&
-\frac34<s<0,
\end{cases}
\end{equation}
we have the bilinear estimates in the modified Bourgain spaces
\begin{equation}
\label{bi-est-X-1}
\|
\p_x(fg)
\|_{X^{s,-b,\theta-1}_{1}}
\le
c_2
\|
f
\|_{X^{s,b',\theta'}_{\alpha}}
\|
g
\|_{X^{s,b',\theta'}_{\alpha}},
\end{equation}
\begin{equation}
\label{bi-est-X-2}
\|
\p_x(fg)
\|_{X^{s,-b,\theta-1}_{\alpha}}
\le
c_2
\|
f
\|_{X^{s,b',\theta'}_{\alpha}}
\|
g
\|_{X^{s,b',\theta'}_{1}}.
\end{equation}
In addition, if $ s_c(\alpha)\le s<3$, then for any 
 $\frac12-\beta\le b'\le b<\frac12<\theta'\le \theta\le \frac12+\beta$, where 
 \begin{equation}
\label{beta-choice}
 \beta
\doteq
 \min\{\beta_1,\frac{3-s}{36}\},
 \end{equation}
we have the bilinear estimates in the ``temporal" Bourgain spaces
\begin{equation}
\label{bi-est-Y-1}
\|
\p_x(fg)
\|_{Y^{s,-b}_{1}}
\le
c_2
\|
\p_x(fg)
\|_{X^{s,-b}_{1}}
+
c_2
\|
f
\|_{X^{s,b'}_{\alpha}}
\|
g
\|_{X^{s,b'}_{\alpha}},
\end{equation}
\begin{equation}
\label{bi-est-Y-2}
\|
\p_x(fg)
\|_{Y^{s,-b}_{\alpha}}
\le
c_2
\|
\p_x(fg)
\|_{X^{s,-b}_{\alpha}}
+
c_2
\|
f
\|_{X^{s,b'}_{\alpha}}
\|
g
\|_{X^{s,b'}_{1}}.
\end{equation}
\end{theorem}
Such bilinear estimates in Bourgain spaces  for  KdV 
have been used by many authors 
(see for example
 \cite{b1993-kdv, kpv1996,  ck2002, ckstt2003, Fa2004,  h2006})
 and play a central role in the analytic theory 
 of dispersive equations.
The key ingredients in the proof of these bilinear estimates 
are the following two Bourgain quantities
\begin{align}
\label{u-Bourgain-q}
&d_{\alpha}(\xi,\xi_1)
\doteq
(\tau-\xi^3)-(\tau_1-\alpha\xi_1^3)-[\tau-\tau_1-\alpha(\xi-\xi_1)^3]
=
-\xi^3+\alpha\xi_1^3+\alpha(\xi-\xi_1)^3,
\\
&\widetilde{d}_{\alpha}(\xi,\xi_1)
\doteq
(\tau-\alpha\xi^3)-(\tau_1-\xi_1^3)-[\tau-\tau_1-\alpha(\xi-\xi_1)^3]
=
-\alpha\xi^3+\xi_1^3+\alpha(\xi-\xi_1)^3,
\end{align}
that appear in the denominator of the multipliers 
when these estimates are expressed in their $L^2$ form 
(see, for example,  \eqref{bilinear-est-L2-form}).
Note that  $d_{\alpha}$ is associated with  the $u$-equation
and $\widetilde{d}_{\alpha}$ with the $v$-equation.
Their zeros, which are lines (at most three) in the $(\xi_1,\xi)$-plane 
passing through the origin, are causing resonances influencing
the value of the critical exponent $s$ above which 
the related bilinear estimates hold. Notice that for $\alpha =1$,
which is the KdV-KdV system case, the two quantities are the same and this results to less resonances and to the lower critical exponent
$-3/4$.

\vskip0.05in
Next, we state our last result which follows from counterexamples to  bilinear estimates.
%
%
%
%
%
%
%
%
%
\begin{theorem}
\label{optimal-thm}
The bilinear estimates in Theorem \ref{mb-bilinear-est} 
 are optimal, i.e. they hold only for $s\ge s_c(\alpha)$.
\end{theorem}
\nin
A direct implication of Theorem \ref{optimal-thm} is that
 our well-posedness Theorems  \ref{thm-MB-half-line-Rob} and \ref{thm-MB-half-line-Dir} are optimal. 

\vskip0.05in
The unified transform method (UTM), also known as the Fokas method  \cite{f1997},  provides a novel approach for solving initial-boundary value problems for linear and integrable nonlinear partial differential equations. In particular, it gives solution formulas for forced linear ibvp's  like the ones used in our work here.
For further results on ibvp's for KdV,  NLS, Boussinesq, heat, and related equations  we refer the reader to \cite{bfo2020,  dtv2014, f1997, f2002, fpbook2015, fis2005, fs2012, hmy2019, hm2015, hm2020, hm2021, fl2012, lf2012a, lenells2013, oy2019, f2002, y2020}
and the references therein.
Finally, for a thorough introduction to the Fokas method 
we refer to the monograph \cite{f2008}. 

The initial value problem for the Majda-Biello system 
on both the line $\rr$ and the circle $\ci$
was investigated by Oh \cite{o2009-1, o2009-2}
using Bourgain spaces and bilinear estimates.
For data in Sobolev space in $H^s(\rr)\times H^s(\rr)$ local and global
well-posedness was proved for $s\ge 0$ when $0<\alpha<1$,
and  for $s> -3/4$ when $\alpha=1$.
On the circle, local and global well-posedness was proved in 
 $H^s(\ci)\times H^s(\ci)$, $s\ge -1/2$, when $\alpha=1$.
 When $0<\alpha<1$, then the Majda-Biello system 
 is locally well-posed on $H^s(\ci)\times H^s(\ci)$
 for $s\ge s^*(\alpha)$ where $s^*(\alpha)$
is described by using Diophantine conditions
that are generated by the resonance relations introduced via the bilinear estimates in the two Bourgain norms corresponding to the two symbols. In particular, this implies 
well-posedness for $s>1/2$ for almost every $\alpha\in (0, 1)$.
Again on the circle,  global well-posedness of the MB system in homogeneous Sobolev spaces  for $s\ge 0$ was proved in 
 \cite{gst2015} using a different approach based on a successive time-averaging method. 

Next we recall the well-posedness results for the KdV equation
with rough data. Its well-posedness in $H^s(\ci)$
(and also in $H^s(\rr)$) when $s\ge 0$ was proved  by Bourgain \cite{b1993-kdv} by introducing the well-known
$X^{s,b}$ spaces. Utilizing these spaces
  Kenig, Ponce and Vega \cite{kpv1996} prove the 
  local well-posedness of KdV on the line for
$s> -3/4$, and  on the circle for $s\ge -1/2$.
Furthermore, Colliander, Keel, Staffilani, Takaoka, Tao
  \cite{ckstt2003} proved  its global well-posedness 
  for the same range of Sobolev exponents.
  It is  worth noticing that for $\alpha=1$ the MB system and the KdV equation  have the same critical well-posedness exponent on $\rr$ or  $\ci$. 
  For additional results we refer the reader to
  \cite{cks1992, bs1975, b1993-nls, {flpv2018}, kato1983, kpv1989, kpv1991-reg, kpv1991, kpv1993, kpv2001, kv2019, lpbook,ckt2003, st1976}
  and the references therein.

%
%
%
%
%
%
\nin
{\it Organization.} The rest of the paper is organized as follows.
In Section 2 we provide the proof of the linear estimates for the MB  Robin and Neumann problems utilizing its UTM solution formulas,
starting with the reduced pure ibvp.
In Section 3, we derive the linear estimates for the 
Dirichlet problem of the forced linear MB system.
Using the linear estimates and the needed bilinear estimates,
in Section 4 we provide  the proof of  well-posedness for the Robin and Neumann problems. The proof for the Dirichlet problem is similar. 
In Section 5 we derive the required bilinear estimates in 
spatial modified Bourgain spaces,
and in Section 6 we derive the corresponding 
bilinear estimates in temporal Bourgain spaces.
Then, in Section 7 we demonstrate the optimality 
of the bilinear estimates proved in Sections 5 and 6.
This implies the optimality of our well-posedness results
(Theorems \ref{thm-MB-half-line-Rob} and \ref{thm-MB-half-line-Dir}).
Finally, in Section 8, we provide an outline for the derivation  of the UTM formula  in the case of the Dirichlet problem and for 
the equation containing the parameter $\alpha$.

%
%
%
%
%
%
%
%
%
%
%
\section{Linear estimates for  Robin and Neumann problems -- proof of Theorem \ref{thm-MB-forced-linear-Rob}
}
\label{Rob-Neu-linear-sec}
\setcounter{equation}{0}
Throughout this work, we shall use  the 
familiar {\bf time localizer} 
 $\psi(t)$, which  is defined as follows:
\begin{equation}
\label{time-localizer}
\psi \in C^{\infty}_0(-1, 1), 
\,\,\,
0\le \psi \le 1
\,\,
\text{ and }
\,\,
\psi(t)=1
\,\,
\text{ for }
\,\,
|t|\le 1/2.
\end{equation}
We prove the linear estimates 
\eqref{MB-forced-linear-est1-Rob}, 
\eqref{MB-forced-linear-est2-Rob} 
for the Robin ibvp,
by decomposing  its forced  linear ibvp to simpler
 problems and estimate each one of them separately.
 We begin with the most important one.
%
%
%
%
%
%
%
%
%
%
%
%
%
\subsection{The reduced pure linear ibvp.}
This ibvp has forcing and initial data equal to zero, and 
boundary data compactly supported in the time interval 
$(0, 2)$, which for the $v$-equation is
\begin{subequations}
\label{Rob-kdv-reduced}
\begin{align}
\label{Rob-kdv-reduced-eqn}
&\p_tv+\alpha\p_x^{3}v=0,
\quad
0<x<\infty,
\,\,
0<t<2,
\\
\label{Rob-kdv-reduced-ic}
&v(x,0)
=
0,
\\
\label{Rob-kdv-reduced-bc}
&v_x(0,t) + \gamma_2 v(0, t) 
= 
h(t),
\quad
\text{supp}(h)\subset (0, 2).
\end{align}
\end{subequations}
 We call this problem the {\bf reduced pure ibvp}.
Since $\text{supp}(h)\subset (0, 2)$ we can write
its time transform
as a Fourier transform over $\rr$, that is
\begin{equation}
\label{Rob-h-time-tran}
\widetilde h(\xi,2)
=
\int_0^2e^{-i\alpha\xi\tau}h(\tau)d\tau
=
\int_\rr e^{-i\alpha\xi\tau}h(\tau)d\tau
=
\hat h(\alpha\xi),
\end{equation}
and using it we write the Fokas solution formula \eqref{Rob-kdv-utm-v} as follows
\begin{align}
\label{Rob-pure-kdv-sln}
v(x, t)
=&
S_\alpha[0,h;0](x,t)
\nn
\\
=&
\frac{3\alpha i}{2\pi}
\int_{\partial D^+}
e^{i\xi x+i\alpha\xi^3 t}  
\frac{\xi^2}{\xi-i\gamma_2}
\widehat{h}(\alpha\xi^3)
d\xi
-
3\alpha\chi_{(0,\infty)}(\gamma_2)
\gamma_2^2
\widehat{h}(-i\alpha\gamma_2^3)
e^{-\gamma_2 x+\gamma_2^3\alpha t}.
\end{align}
Utilizing formula \eqref{Rob-pure-kdv-sln} 
we derive the following estimates in 
Bourgain spaces. 
\begin{theorem}
[\b{Estimates for pure ibvp of linear $\alpha$-KdV (Robin)}]
\label{Rob-reduced-thm-kdv}
For boundary data  a test function $h(t)$  with 
$\text{supp}(h)\subset (0, 2)$, $b\in[0,\frac12)$ 
and
$
0
\le
\theta
\le
1+s/3$,
the solution of the  linear $\alpha$-KdV reduced pure ibvp \eqref{Rob-kdv-reduced} satisfies the  following estimate
in modified Bourgain spaces
\begin{align}
\label{Rob-reduced-pure-ibvp-B-es}
\|
S_{\alpha}[0,h;0]
\|_{X^{s, b,\theta}_{\alpha,\rr^+\times (0,2)}} 
\leq
c_{s,b,\theta,\alpha}
\|h\|_{H_t^{\frac{s}{3}}(\rr)},
\,\,\,
s>-\frac32, \,\, 
\text{if} \,\, \gamma_2\le 0,
\,\,\,
\text{and}
\,\,\,
s\ge 0, \,\, \text{if} \,\, \gamma_2>0.
\end{align}
\end{theorem}
\nin
{\bf Proof of Theorem \ref{Rob-reduced-thm-kdv}.}
We follow the lines of the Robin ibvp for KdV  ($\alpha=1$) in \cite{hy2022-Robin}.
 Parametrizing  the right (left) contour of $\p D^+$ via
$[0,\infty)\ni \xi\to a \xi$ ($a^2\xi$) with 
$a
=
e^{i\frac{\pi}{3}}
$,
we write the solution as 
$
v(x,t)
\simeq
v_{r}(x,t)+v_{\ell}(x,t)+v_i(x,t),
$
where
\begin{align}
\label{Rob-def-vr}
v_r(x,t)
=
\int_0^\infty e^{ia \xi x-i\alpha\xi^3t}
\frac{a^2\xi^2}{a\xi-i\gamma_2}
\hat h(-\alpha\xi^3)d\xi,
\end{align}
\vskip-0.15in
\begin{align*}
v_{\ell}(x,t)
=
\int_0^\infty e^{ia^2\xi x+i\alpha\xi^3t}
\frac{a^4\xi^2}{a^2\xi-i\gamma_2}
\hat h(\alpha\xi^3)d\xi,
\end{align*}
\vskip-0.15in
\begin{align*}
v_i(x,t)
=
\chi_{(0,\infty)}(\gamma_2)
\gamma_2^2
\widehat{h}(-i\alpha\gamma_2^3)
e^{-\gamma_2 x+\gamma_2^3\alpha t},
\end{align*}
with $v_i$ corresponding to the residue
at the pole $i\gamma_2$ inside  $D^+$,
occurring only when $\gamma_2>0$.
Below we estimate only $v_{r}$.
 The estimate $v_\ell$ is similar and we omit it. 
Furthermore,  we split the $\xi$-integral 
defining  \eqref{Rob-def-vr}
 for $\xi$ near 0 and  for $\xi$ away from 0,
 that is, we write 
$
v_r
=
v_0
+
v_1,
$
where
\begin{equation*}
v_0(x,t)
\doteq
\int_0^1 e^{-i\alpha\xi ^3t}e^{ia_R\xi x}e^{-a_I \xi x} 
\frac{a^2\xi^2}{a\xi-i\gamma_2}
 \, \widehat h(-\alpha\xi ^3)d\xi,
\quad
x\in\rr^+,
\quad
t\in(0,2),
\end{equation*}
\begin{equation}
\label{Rob-def-v1-kdv}
v_1(x,t)
\doteq
\int_1^\infty e^{-i\alpha\xi ^3t}e^{ia_R\xi x}e^{-a_I \xi x} 
\frac{a^2\xi^2}{a\xi-i\gamma_2}
\,
\widehat h(-\alpha\xi ^3)d\xi,
\quad
x\in\rr^+
\quad
t\in(0,2),
\end{equation}
with  $a=a_R+ia_I=\frac12+i\frac{\sqrt3}{2}$, and prove
estimate \eqref{Rob-reduced-pure-ibvp-B-es} 
for each one separately. We begin with $v_0$.

%
%
%
%
%
%
%
%
%
%
%
%
\vskip.05in
\noindent
{\bf Proof of estimate \eqref{Rob-reduced-pure-ibvp-B-es} for $v_0$.}
First we extend $v_0$  as an ``almost" even function $V_0$
of $x$ from $\rr^+\times (0,2)$ to $\rr\times\rr$  as follows
\begin{equation*}
V_0(x,t)
\doteq
\int_0^1 e^{ia\xi \varphi_1(x)}e^{-i\alpha\xi ^3t}
\frac{a^2\xi^2}{a\xi-i\gamma_2}
\,\widehat {h}(-\alpha\xi ^3)d\xi,
\quad
x\in\rr,
\,\,
t\in\rr,
\end{equation*}
where $\varphi_1(x)$ is a smooth version of 
$|x|$ described below and drawn in Figure \ref{fig:phi1}

\vskip-0.05in
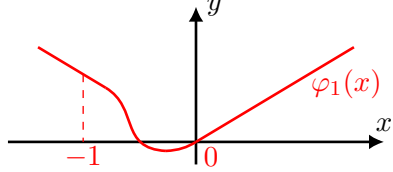
\begin{minipage}{0.49\linewidth}
\begin{center}
\begin{tikzpicture}[yscale=0.6, xscale=1]
%
%
\newcommand\X{0};
\newcommand\Y{0};
\newcommand\FX{11};
\newcommand\FY{11};
\newcommand\R{0.6};
\newcommand*{\TickSize}{2pt};
%
%
\draw[black,line width=1pt,-{Latex[black,length=2mm,width=2mm]}]
(-2.5,0)
--
(2.5,0)
node[above]
{\fontsize{\FX}{\FY}\bf \textcolor{black}{$x$}};

\draw[black,line width=1pt,-{Latex[black,length=2mm,width=2mm]}]
(0,-0.5)
--
(0,3)
node[right]
{\fontsize{\FX}{\FY}\bf \textcolor{black}{$y$}};

\draw[line width=1pt, yscale=1,domain=-2.1:-1.2,smooth,variable=\x,red]  plot ({\x},{-\x});

\draw[line width=1pt, yscale=1,domain=0:2.1,smooth,variable=\x,red]  plot ({\x},{\x});

\draw[smooth,line width=1pt, red]
(-0.4,-0.2)
to[out=5,in=-135]
(0,0)
;

\draw[smooth,line width=1pt, red]
(-1.2,1.2)
to[out=-45,in=175]
(-0.4,-0.2)
;

\draw[red,dashed, line width=0.5pt]
(2,1.3)
node[]
{\fontsize{\FX}{\FY}$\varphi_1(x)$}

(0,0)
node[yshift=-0.2cm,xshift=0.2cm]
{\fontsize{\FX}{\FY}$0$}

(-1.5,1.5)
--
(-1.5,0)
node[yshift=-0.2cm]
{\fontsize{\FX}{\FY}$-1$};

\end{tikzpicture}

\vskip-0.15in
\captionof{figure}{$\varphi_1(x)$}
\label{fig:phi1}
\end{center}
\end{minipage}
\begin{minipage}{0.48\linewidth}
\begin{equation}
\label{phi1-x-even-ext}
\varphi_1(x)
=
\begin{cases}
x,
\quad
\,\,\,\,
x\geq
0,
\\
-x,
\quad
x\leq
-1,
\\
\text{ smooth on } \rr.
\quad
\end{cases}
\end{equation}
\end{minipage}

\noindent
From the  inequalities
$
(1+|\xi|)^{2s}
\lesssim
(|\xi|^0+|\xi|^{2|s|}),
$
$
(1+|\tau|)^{2\theta}
\lesssim
(|\tau|^0+|\tau|^{2\theta}),
$
and
$
(1+|\tau-\alpha\xi^3|)^{2b}
\leq
(1+|\tau|+\alpha|\xi|^3)^{2b}
\leq
(1+|\tau|)^{2b}
(1+\alpha|\xi|^3)^{2b}
\lesssim
(|\tau|^0+|\tau|^{2b})
(|\xi|^0+|\xi|^{6b}),
$
we get 
\begin{equation}
\label{Rob-est-sb-norm-fact-no3}
(1+|\xi|)^{2s}
(1+|\tau-\alpha\xi^{3}|)^{2b}
+
\chi_{|\xi|\le 1}
(1+|\tau|)^{2\theta}
\lesssim
(|\xi|^0+|\xi|^{6b+2|s|})
(|\tau|^0+|\tau|^{2b+2\theta}).
\end{equation}
Using \eqref{Rob-est-sb-norm-fact-no3} we obtain  the following  $L^2$ estimate for
$\psi_4(t)V_0$, with  $\psi_4(t)=\psi(t/4)$,
\begin{align}
\label{Rob-v0-L2-est}
\|\psi_4 V_0\|^2_{X^{s,b,\theta}_\alpha}
\lesssim
\|\psi_4 V_0\|^2_{L^2_{x,t}}
+
\|\p_x^{n_1}[\psi_4 V_0]\|^2_{L^2_{x,t}}
+
\|\p_t^{n_2}[\psi_4 V_0]\|^2_{L^2_{x,t}}
+
\|\p_x^{n_1}\p_t^{n_2}[\psi_4 V_0]\|^2_{L^2_{x,t}},
\end{align}
where 
$n_1
=
n_1(s,b)
\doteq
6\lfloor b\rfloor+2\lfloor |s|\rfloor+8
$
and
$
n_2
=
n_2(b,\theta)
\doteq
2\lfloor b\rfloor+2\lfloor\theta\rfloor+4.
$
Now, from \eqref{Rob-v0-L2-est}  and an appropriate
application of the $L^2$ boundness of Laplace transform 
\cite{fhm2017},  \cite{hardy1933}
(for $\alpha=1$ a detailed  proof is provided in \cite{hy2022-Robin})
we get
\begin{equation}
\label{Rob-L2-V0-est-kdv}
\|\p_x^{n_1}\p_t^{n_2}[\psi_4\,V_0]\|_{L^2_{x,t}}^2
\le
C_{n_1, n_2}
\int_{-1}^{0}
|
\widehat{h}(\alpha\tau)
|^2
d\tau,
\end{equation}
\vskip-0.15in
\nin
which gives 
\begin{align*}
\|\psi_4 V_0\|^2_{X^{s,b,\theta}_{\alpha}}
\le
C_{n_1, n_2}
\int_{-1}^{0}
(1+|\alpha\tau|)^{-2\cdot\frac{s}{3}}
(1+|\alpha\tau|)^{2\cdot\frac{s}{3}}
|
\widehat{h}(\alpha\tau)
|^2
d\tau
\lesssim
\|h\|_{H_t^{\frac{s}{3}}}^2,
\quad
s\in\rr,
b\ge 0,
\theta\ge 0,
\end{align*}
thus arriving at the desired estimate \eqref{Rob-reduced-pure-ibvp-B-es} for $v_0$, that is
$
\|
v_0
\|_{X^{s, b,\theta}_{\alpha,\rr^+\times (0,2)}}\le \|\psi_4 V_0\|_{X^{s,b,\theta}_{\alpha}}\lesssim \|h\|_{H_t^{\frac{s}{3}}}.
$
\,\,
$\square$
%

%
%
%
%
%
%
%
%
%
%
%
%
\noindent
{\bf Proof of estimate \eqref{Rob-reduced-pure-ibvp-B-es} for $v_1$.}
Using the identity
$
\xi
[
e^{ia_R\xi x}
e^{-a_I\xi x}
]
=
\frac{1}{ia}
\p_x
[
e^{ia_R\xi x}
e^{-a_I\xi x}
]
$
and the fact that  $e^{-a_I \xi x}$ decays exponentially 
 in $\xi$, since  $x>0$, we
move the $\p_x$-derivative outside the integral sign 
in \eqref{Rob-def-v1-kdv} and 
rewrite $v_1$ as follows
\begin{equation*}
v_1(x,t)
=
\frac{1}{ia}
\p_x
\int_1^\infty e^{-i\alpha\xi ^3t}e^{ia_R\xi x}e^{-a_I \xi x} 
\frac{a^2\xi}{a\xi-i\gamma_2}
\,
\widehat h(-\alpha\xi ^3)d\xi,
\quad
x\in\rr^+,
\quad
t\in(0,2).
\end{equation*}
Next, we need an one-sided  $C^\infty$ cutoff function $\rho(x)$,  $0\le\rho(x)\le1$, $x\in\rr$, defined as follows outside
the interval $(-1, 0)$, and is increasing in $(-1, 0)$
\vskip.01in
\noindent
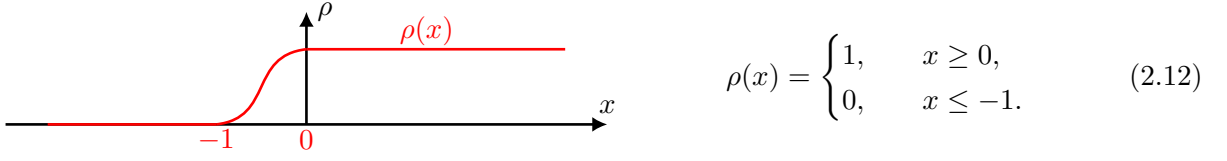
\begin{minipage}{0.6\linewidth}
\begin{center}
\begin{tikzpicture}[yscale=0.5, xscale=0.8]
%
%
\newcommand\X{0};
\newcommand\Y{0};
\newcommand\FX{11};
\newcommand\FY{11};
\newcommand\R{0.6};
\newcommand*{\TickSize}{2pt};
%
%
\draw[black,line width=1pt,-{Latex[black,length=2mm,width=2mm]}]
(-5,0)
--
(5,0)
node[above]
{\fontsize{\FX}{\FY}\bf \textcolor{black}{$x$}};

\draw[black,line width=1pt,-{Latex[black,length=2mm,width=2mm]}]
(0,0)
--
(0,3)
node[right]
{\fontsize{\FX}{\FY}\bf \textcolor{black}{$\rho$}};

\draw[line width=1pt, yscale=2,domain=-1.5:-4.3,smooth,variable=\x,red]  plot ({\x},{0});

\draw[line width=1pt, yscale=2,domain=0:4.3,smooth,variable=\x,red]  plot ({\x},{1});

\draw[smooth,line width=1pt, red]
(0,2)
to[out=-170,in=10]
(-1.5,0)
;

\draw[red]
(2,2.5)
node[]
{\fontsize{\FX}{\FY}$\rho(x)$}

(0,0)
node[yshift=-0.2cm]
{\fontsize{\FX}{\FY}$0$}

(-1.5,0)
node[yshift=-0.2cm]
{\fontsize{\FX}{\FY}$-1$};

\end{tikzpicture}
\end{center}
\end{minipage}
\hskip-0.88in
\begin{minipage}{0.5\linewidth}
\begin{equation}
\label{Rob-rho-def}
\rho(x)
=
\begin{cases}
1,
&\quad
x\ge0,
\\
0,
&\quad
x\le -1.
\end{cases}
\end{equation}
\end{minipage}

\nin
Using $\rho(x)$ we extend $v_1$  from $\rr^+\times(0,2)$ to 
 $V_1$ on $\rr\times\rr$ by the formula
%
\begin{align*}
V_1(x,t)
\doteq&
\frac{1}{ia}
\p_x
\int_1^\infty 
e^{-i\alpha\xi ^3t}e^{ia_R\xi x}e^{-a_I \xi x} \rho(a_I\xi x)
\frac{a^2\xi}{a\xi-i\gamma_2}
\,
\widehat h(-\alpha\xi ^3)d\xi,
\quad
x\in\rr,
\quad
t\in\rr.
\end{align*}
Since 
$
\|v_1\|_{X^{s,b,\theta}_{\alpha,\rr^+\times(0,2)}}
\le 
\|V_1\|_{X^{s,b,\theta}_{\alpha}},
$
 to prove  estimate \eqref{Rob-reduced-pure-ibvp-B-es} 
for $v_1$ it suffices to show that
\begin{align}
\label{Rob-reduced-ibvp-Bourgain-est-1}
\|V_1\|_{X^{s,b,\theta}_{\alpha}} 
\le
c_{s,b,\alpha}
\|h\|_{H_t^{\frac{s}{3}}(\rr)}.
\end{align}
To prove this,  we rewrite the modified Bourgain norm  as follows
\begin{align*}
\|V_1\|_{X^{s,b,\theta}_{\alpha}}^2
=
\|V_1\|_{X^{s,b}_{\alpha}}^2
+
\|V_1\|_{D^\theta}^2,
\,\,
\text{with}
\,\,
\|V_1\|_{D^\theta}^2
\doteq
\int_{-\infty}^\infty
\int_{-1}^1
(1+|\tau|)^{2\theta}
|\widehat{V}_1(\xi,\tau)|^2
d\xi
d\tau.
\end{align*}
\underline{Estimate of the Bourgain $X^{s,b}_{\alpha}$ norm for $V_1$}.
Making the change of variables $\tau=-\xi^3$, 
and defining  the Schwartz function
\begin{equation}
\label{Rob-def-eta}
\eta(x)
\doteq
e^{i\frac{a_R}{a_I}x}
e^{-x} \rho(x),
\end{equation}
we write $V_1$ in the form
\begin{equation}
\label{Rob-eqn-V1-1}
V_1(x,t)
\simeq
\p_x
\int_{-\infty}^{-1} e^{i\alpha\tau t}
\eta(-a_I\tau^{1/3}x)
\frac{a^2\tau^{1/3}}{a\tau^{1/3}-i\gamma_2}
\,
\tau^{-\frac23}
\,
\widehat h(\alpha\tau)
d\tau,
\quad
x\in\rr,
\quad
t\in\rr.
\end{equation}
Also, letting $V_{1,\alpha}(x,t)\doteq V_1(x,\frac{t}{\alpha})$ and 
using the  identity
$
\widehat{V}_{1,\alpha}(\xi,\tau)
=
\alpha
\widehat{V}_{1}(\xi,\alpha\tau),
$
we get
$
\|V_1\|_{X_\alpha^{s,b}}
\lesssim
\|V_{1,\alpha}\|_{X^{s,b}_1}
$
($b\ge 0$),
which reduces estimate \eqref{Rob-reduced-ibvp-Bourgain-est-1} to the following one
\begin{align}
\label{Rob-reduced-ibvp-Bourgain-est-2}
\|V_{1,\alpha}\|_{X^{s,b}_{1}} 
\le
c_{s,b,\alpha}
\|h\|_{H_t^{\frac{s}{3}}(\rr)},
\end{align}
\vskip-0.1in
\nin
where, by  \eqref{Rob-eqn-V1-1}, 
\begin{equation}
\label{Rob-eqn-V1-alpha}
V_{1,\alpha}(x,t)
\simeq
\p_x
\int_{-\infty}^{-1} e^{i\tau t}
\eta(-a_I\tau^{1/3}x)
\frac{a^2\tau^{1/3}}{a\tau^{1/3}-i\gamma_2}
\,
\tau^{-\frac23}
\,
\widehat h(\alpha\tau)
d\tau,
\quad
x\in\rr,
\quad
t\in\rr.
\end{equation}
Now we claim that the  Bourgain space estimate 
\eqref{Rob-reduced-ibvp-Bourgain-est-2}  follows from next result. 
\begin{lemma}
\label{Rob-reduced-pibvp-thm1}
Let $\varepsilon>0$. If $s\geq -\frac32-\varepsilon$, 
and $b\ge 0$,
then   $V_{1,\alpha}$  in  \eqref{Rob-eqn-V1-alpha}
satisfies the estimate
\begin{align}
\label{Rob-reduced-ibvp-Bourgain-est}
\|V_{1,\alpha}\|_{X^{s,b}_1} 
\leq
c_{s,b,\varepsilon,\rho,a,\alpha} \|h\|_{H_t^{\frac{s+3b-\frac32+\varepsilon}{3}}(\rr)}.
\end{align}
\end{lemma}
\vskip-0.05in
\noindent
In fact, choosing  $\varepsilon$ such that 
$
\varepsilon
=
3(\frac12-b),
$
which is possible if $b<1/2$, we get $(s+3b-3/2+\varepsilon)/3= s/3$
and estimate  \eqref{Rob-reduced-ibvp-Bourgain-est}
gives \eqref{Rob-reduced-ibvp-Bourgain-est-2}. And, since
$\|V_1\|_{X_\alpha^{s,b}}\lesssim \|V_{1,\alpha}\|_{X^{s,b}_1}$ 
we get the desired estimate \eqref{Rob-reduced-ibvp-Bourgain-est-1}.

\vskip0.1in
\nin
{\bf Proof of Lemma \ref{Rob-reduced-pibvp-thm1}.} 
To estimate $\|V_{1,\alpha}\|_{X^{s,b}_1}$ we need 
 the Fourier transform of $V_{1,\alpha}$. 
Using formula \eqref{Rob-eqn-V1-alpha} we see that its time
Fourier transform is
\begin{align*}
\widehat{V}_{1,\alpha}^t(x,\tau)
\simeq
\p_x
\big[
\chi_{(-\infty,-1)}(\tau)
\,
\eta(-a_I\tau^{1/3}x)
\,
\frac{a^2\tau^{1/3}}{a\tau^{1/3}-i\gamma_2}
\,
\tau^{-\frac23}
\,
\widehat h(\alpha\tau)
\big],
\end{align*}
and therefore its  full Fourier transform is
\begin{align}
\label{Rob-v1-xt-FT}
\widehat{V}_{1,\alpha}(\xi,\tau)
\simeq
\xi
\,
\chi_{(-\infty,-1)}(\tau)
\,
F(\xi,\tau)
\,
\frac{a^2\tau^{1/3}}{a\tau^{1/3}-i\gamma_2}
\,
\tau^{-\frac23}
\,
\widehat h(\alpha\tau),
\,\,
\text{with}
\,\,
F(\xi,\tau)
\doteq
\int_{\rr}
e^{-i\xi x}\eta(-a_I\tau^\frac13 x)
dx.
\end{align}
Also, we need the following estimate for $F$,
which follows from the fact
that $\eta$ is a Schwartz function.
\begin{lemma}
\label{Rob-F-bound-lem}
For any 
$n\ge 0$,  $\tau<-1$, and $\xi\in\rr$, we have
\begin{align*}
|F(\xi,\tau)|
\leq
c_{\rho,a,n}
\,
\frac{1}{|\tau|^{1/3}}
\,
\big(
\frac{|\tau|^{1/3}}{|\xi|+|\tau|^{1/3}}
\big)^n,
\end{align*}
where $c_{\rho,a,n}$ is a constant depending on 
$a$, $n$ and $\rho$, which is described in \eqref{Rob-rho-def}.
\end{lemma}
\nin
Now,  using \eqref{Rob-v1-xt-FT} we write $\|V_{1,\alpha}\|_{X^{s,b}_1}^2$ as follows
\begin{align}
\label{Rob-X-norm-v1}
\|V_{1,\alpha}\|_{X^{s,b}_1}^2
\lesssim&
\int_{-\infty}^{-1}
\Big[
\int_{-\infty}^{\infty}
(1+|\xi|)^{2s}
(1+|\tau-\xi^3|)^{2b}
\xi^2
|
F(\xi,\tau)
|^2
d\xi
\Big]
\cdot
\frac{\tau^{2/3}}{|a\tau^{1/3}-i\gamma_2|^2}
\,
\tau^{-\frac43}
\big|
\widehat h(\alpha\tau)
\big|^2
d\tau,
\end{align}
and estimate  the $d\xi$-integral in the following result.
\begin{lemma}
\label{Rob-v1-mult-lema}
For any $\varepsilon>0$ and $\tau<-1$, if $s\geq -\frac32-\varepsilon$, $b\ge 0$ then we have
\begin{equation}
\label{Rob-v1-mult-ine-1}
\int_{-\infty}^{\infty}
(1+|\xi|)^{2s}
(1+|\tau-\xi^3|)^{2b}
\xi^2
\big|
F(\xi,\tau)
\big|^2
d\xi
\leq
c_{s,b,\varepsilon,\rho,a}
|\tau|^{\frac{2s+6b+1+2\varepsilon}{3}},
\end{equation}
where $c_{s,b,\varepsilon,\rho}$ is a constant depending on $s$, $b$, $\varepsilon$, $\rho$ and $a$.
\end{lemma}
\noindent
Estimate  \eqref{Rob-v1-mult-ine-1} is proved in \cite{hy2022-KdVm} (see Lemma 2.3). Now,
combining  it with \eqref{Rob-X-norm-v1},
 we get
\begin{align*}
\|V_{1,\alpha}\|_{X^{s,b}_{1}}^2
\le
c_{s,b,\varepsilon,\rho,a}
\int_{-\infty}^{-1}
\big|
\tau^{\frac{s+3b-\frac32+\varepsilon}{3}}
\widehat h(\alpha\tau)
\big|^2
d\tau
\le
c_{s,b,\varepsilon,\rho,a} \|h\|_{H_t^{\frac{s+3b-\frac32+\varepsilon}{3}}}^2,
\end{align*}
which is  the desired estimate \eqref{Rob-reduced-ibvp-Bourgain-est}
in Lemma \ref{Rob-reduced-pibvp-thm1}.
\,\,
$\square$

\vskip0.1in
\nin
\underline{Estimation of $D^\theta$ norm for $V_1$.}  
As before, letting $V_{1,\alpha}(x,t)\doteq V_1(x,\frac{t}{\alpha})$ and 
using the  identity
$
\widehat{V}_{1,\alpha}(\xi,\tau)
=
\alpha
\widehat{V}_{1}(\xi,\alpha\tau),
$
we get
$
\|V_1\|_{D^\theta}
\lesssim
\|V_{1,\alpha}\|_{D^\theta}
$, if $\theta\ge0$.
Using the Fourier transform of $V_{1,\alpha}$, i.e. \eqref{Rob-v1-xt-FT} we get
\begin{align*}
\|V_{1,\alpha}\|_{D^\theta}^2
=
\int_{-\infty}^{-1}
\int_{-1}^1
(1+|\tau|)^{2\theta}
|
\xi
F(\xi,\tau)
\frac{a^2\tau^{1/3}}{a\tau^{1/3}-i\gamma_2}
\tau^{-2/3}
\widehat h(\alpha\tau)|^2
d\xi
d\tau,
\end{align*}
where 
$
F(\xi,\tau)
$
is defined in \eqref{Rob-v1-xt-FT}.
Next, applying Lemma \ref{Rob-F-bound-lem} with $n=0$, we get
$
|F(\xi,\tau)|
\lesssim
|\tau|^{-1/3},
$
 which implies that the integrand is bounded in $\xi$.
Thus, integrating $\xi$ we get
\begin{align*}
\|V_{1,\alpha}\|_{D^\theta}^2
\lesssim
\int_{-\infty}^{-1}
(1+|\tau|)^{2\theta}
|
\tau^{-1}
\widehat h(\alpha\tau)|^2
d\tau
\lesssim
\int_{-\infty}^{-1}
(1+|\tau|)^{2\theta-2}
|
\widehat h(\alpha\tau)|^2
d\tau
\lesssim
\|h\|_{H_t^{\theta-1}}^2,
\end{align*}
which gives the desired estimate 
$\|V_{1,\alpha}\|_{D^\theta}^2 \lesssim \|h\|_{H_t^{s/3}}^2$,
if we choose  $\theta-1\le s/3$ or 
$
0
\le
\theta
\le
1+s/3.
$
\,\,
$\square$

%
%
%
%
%
%
%
%
%
%
%
%
\vskip.05in
\noindent
{\bf Proof of estimate \eqref{Rob-reduced-pure-ibvp-B-es} for $v_i$
(where the signs of $\gamma_1$, $\gamma_2$ appear).}
When $\gamma_2\le 0$, then $v_i=0$ and  estimate
 \eqref{Rob-reduced-pure-ibvp-B-es} holds trivially.
So, we focus on  $\gamma_2> 0$ in which case 
$v_i(x,t)
=
\gamma_2^2
\widehat{h}(-i\alpha\gamma_2^3)
e^{-\gamma_2 x+\gamma_2^3\alpha t}.
$
Replacing $x$ with $ \varphi_1(x)$ defined in \eqref{phi1-x-even-ext}
we extend  $v_i$  as an ``almost" even function
of $x$ from $\rr^+\times(0,2)$ to $\rr\times\rr$ 
 as follows 
\begin{equation}
\label{Rob-def-vi-extension}
V_i(x,t)
\doteq
\gamma_2^2
\widehat{h}(-i\alpha\gamma_2^3)
e^{-\gamma_2 \varphi_1(x)+\gamma_2^3\alpha t},
\quad
x\in\rr,
\,\,
t\in\rr.
\end{equation}
Thus, for its modified Bourgain norm 
we have 
(by splitting integration for $|\tau|\le |\xi|$ and $|\tau|>|\xi|$)
\begin{align}
\label{vi-est}
\nn
\|v_i\|_{X^{s, b,\theta}_{\alpha,\rr^+\times (0,2)}}^2
&\le
\|\psi_4 V_i\|_{X^{s, b,\theta}_{\alpha}}^2
\le c_{\gamma_2} 3 \gamma_2^4
\Big(\|e^{-\gamma_2 \varphi_1}\|_{H^{|s|+3b}}^2
\|\psi_4 e^{\gamma_2^3\alpha t}\|_{H^{|s|+3b+\theta}}^2
\Big)
\cdot
|\widehat{h}(-i\alpha\gamma_2^3)|^2
\\
&\le C_{\psi, \phi_1,\gamma_2}
|\widehat{h}(-i\alpha\gamma_2^3)|^2.
\end{align}
Also, using the definition \eqref{Rob-h-time-tran} and the Cauchy--Schwarz inequality  we have
\begin{equation}
\label{h-l2-estt}
|\widehat{h}(-i\alpha\gamma_2^3)|^2
=
\Big|
\int_0^2e^{-i (-i\alpha\gamma_2^3)\tau}h(\tau)d\tau
\Big|^2
=
\Big|
\int_0^2e^{-\alpha \gamma_2^3\tau}h(\tau)d\tau
\Big|^2
\le
\Big|
\int_0^2 h(\tau)d\tau
\Big|^2
\le
2
\|h\|_{L^2}^2.
\end{equation}
Therefore,  combining estimate \eqref{vi-est} with \eqref{h-l2-estt}, 
for $s\ge 0$ we get 
$
\|v_i\|_{X^{s, b,\theta}_{\alpha,\rr^+\times (0,2)}}^2
\le
\|\psi_4 V_i\|_{X^{s, b,\theta}_{\alpha}}^2
\lesssim
\|h\|_{H^{s/3}}^2,
$
which completes the proof of estimate  \eqref{Rob-reduced-pure-ibvp-B-es} for $v_i$.
\,\,
$\square$

%
%
%
%
%
%
%
%
%
%
\subsection{Proof of Linear estimates for Robin ibvp (Theorem \ref{thm-MB-forced-linear-Rob})}
It suffices to prove estimate \eqref{MB-forced-linear-est2-Rob}, since estimate \eqref{MB-forced-linear-est1-Rob} follows 
from this one by letting $\alpha=1$.
We start by decomposing the forced  problem  \eqref{Rob-LKdV-v} for the linear $\alpha$-KdV (i.e. the KdV with dispersion coefficient $\alpha$) 
into the homogeneous ibvp 
\begin{subequations}
\label{homo-ibvp-kdv}
\begin{align}
\label{homo-ibvp-kdv:eqn}
&v_t+\alpha v_{xxx}
=
0,
\\
\label{homo-ibvp:ic}
&v(x,0) 
=
v_0(x)\in H^{s}(\rr^+),
\\
\label{homo-ibvp-kdv:bc}
&\p_xv(0,t)+\gamma_2v(0,t)
=
\varphi_2(t)\in H^{\frac{s}{3}}(0,T),
\end{align}
\end{subequations}
and the remaining forced ibvp with zero  data
\begin{subequations}
\label{inhomo-ibvp-kdv}
\begin{align}
\label{inhomo-ibvp-kdv:eqn}
&v_t+\alpha v_{xxx}
=
f_2(x,t)
\in 
X_{\alpha,\rr^+\times(0,T)}^{s,-b,\theta-1}
\cap
Y_{\alpha,\rr^+\times(0,T)}^{s,-b},
\\
\label{inhomo-ibvp-kdv:ic}
&v(x,0) 
= 
0, 
\\
\label{inhomo-ibvp-kdv:bc}
&\p_xv(0,t)+\gamma_2v(0,t)
= 
0,
\end{align}
\end{subequations}
thus having the  relation 
$
S_\alpha\big[v_0, \varphi_2; f_2\big]
=
S_\alpha\big[v_0, \varphi_2; 0\big]
+
S_\alpha\big[0, 0; f_2\big]
$
between their solutions. 
Furthermore, 
we  decompose the homogeneous ibvp \eqref{homo-ibvp-kdv}
into an ivp  and the remaining  pure ibvp having zero initial data. The ivp reads as follows
\begin{subequations}
\label{homo-ivp-kdv}
\begin{align}
\label{homo-ivp-kdv:eqn}
&V_t+  \alpha V_{xxx}
=
0,  
\\
\label{homo-ivp-kdv:ic}
&V (x,0)
= 
V_0(x)\in H^{s}(\rr), 
\end{align}
\end{subequations}
where $V_0\in H^{s}(\rr)$ is an extension of the initial data $v_0\in H^{s}(\rr^+)$ with
$
\|V_0\|_{H^{s}(\rr)} \le 2\|v_0\|_{H^{s}(\rr^+)}
$.
Its solution   is given by the familiar formula 
\begin{equation}
\label{homo-ivp-sln-kdv}
V(x,t) 
= 
S_\alpha\big[V_0; 0\big] (x, t) 
= 
\frac{1}{2\pi} \int_{\xi\in \rr} e^{i\xi x+i\alpha\xi^3t}\, \widehat{V}_0(\xi) d\xi,\quad (x,t)\in \rr^2,
\end{equation}
where
$
\widehat{V}_0(\xi) 
=
\int_{x\in \rr} e^{-i\xi x}\, V_0(x) dx
$,
and satisfies the following result.
\begin{proposition}
\label{pro-ivp-kdv}
For $s, b,\theta\in \rr$, the solution \eqref{homo-ivp-sln-kdv} to ivp  \eqref{homo-ivp-kdv}
satisfies the space estimates
\begin{align}
\label{homo-ivp-est-B-kdv}
\|\psi\, S_\alpha\big[V_0; 0\big]\|_{X^{{s},b,\theta}_\alpha}
\le
c_{s}
\|V_0\|_{H^{s}(\rr)}, 
\end{align}
and the following  time regularity estimates
\begin{align}
\label{homo-ivp-est-time-kdv}
&\sup_{x\in\rr}\|
\psi
\,
S_\alpha\big[V_0; 0\big](x)\|_{H_t^{\frac{s+1}{3}}(\rr)} 
\le
C_{s}\,  \|V_0\|_{H^{s}(\rr)},
\\
\label{homo-ivp-est-time-kdv-der}
&\sup_{x\in\rr}\|
\psi
\,
\p_x
S_\alpha\big[V_0; 0\big](x)\|_{H_t^{\frac{s}{3}}(\rr)} 
\le
C_{s}\,  \|V_0\|_{H^{s}(\rr)}.
\end{align}
\end{proposition}

\noindent
{\bf Proof of Proposition \ref{pro-ivp-kdv}.} 
For $\alpha=1$,
estimate \eqref{homo-ivp-est-time-kdv} (with $s\ge -1$) is proved in \cite{fhm2016} and estimate \eqref{homo-ivp-est-time-kdv-der} (with $s\ge 0$) is proved in \cite{hmy2019}. 
For  $\alpha\neq 1$ and $s\in\rr$, the proof is similar.
Estimate \eqref{homo-ivp-est-B-kdv} follows from the
definition of  Bourgain norm (also, see Lemma 5.5 in \cite{h2006}). 
\,\,
$\square$

\vskip.05in
\noindent
The remaining pure ibvp mentioned above  is given by
\begin{subequations}
\label{pure-ibvp-kdv}
\begin{align}
\label{pure-ibvp-kdv:eqn}
&v_t+\alpha v_{xxx} 
= 
0, 
\\
\label{pure-ibvp-kdv:ic} 
&v(x,0)
= 
0, 
\\
\label{pure-ibvp-kdv:bc}
&\p_xv(0,t)+\gamma_2v(0,t)
= 
\varphi_2(t)
-
[\p_xV(0,t)+\gamma_2V(0,t)]
\doteq 
\varphi(t)
\in H^{\frac{s}{3}}(0,T).
\end{align}
\end{subequations}
For its solution
$
v
\doteq
S_\alpha \big[0, \varphi; 0\big]
$, defined by UTM formula \eqref{Rob-kdv-utm-v}, we have the following result.
\begin{proposition}
\label{pro-p-ibvp-kdv}
For $b\in [0,\frac12)$,
 $\frac12<\theta
 \le
1+s/3$, we have the estimate in Bourgain spaces
\begin{align}
\label{pure-ibvp-est-B-kdv}
\|
S_\alpha [0,\varphi;0]\|_{X^{s,b,\theta}_{\alpha ,\rr^+\times(0,T)}}
\le 
C_{s}\, \|\varphi\|_{H^{\frac{s}{3}}(0,T)},
\,\,\,
-\frac32 < s< \frac32,
 \,\, 
\text{if} \,\, \gamma_2\le 0,
\,\,\,
\text{and}
\,\,\,
0\le s<\frac32, \,\, \text{if} \,\, \gamma_2>0.
\end{align}
\end{proposition}

\nin
{\bf Proof of Proposition \ref{pro-p-ibvp-kdv}.} Since $-\frac32<s<\frac32$, we have   $-\frac12<\frac{s}{3}<\frac12$. 
Now, we  extend $\varphi(t)$ from $(0,T)$ to $\rr$ such that the extension $\widetilde{\varphi}(t)$ satisfies
$
\|\widetilde{\varphi}\|_{H^{\frac{s}{3}}(\rr)}
\le
2
\|\varphi\|_{H^{\frac{s}{3}}(0,T)}.
$
Furthermore, multiplying $\varphi$ by the characteristic functions $\chi_{(-\infty,T)}$ and $\chi_{(0,\infty)}$ we get
$
\widetilde{\varphi}_1
\doteq 
\chi_{(-\infty,T)}\chi_{(0,\infty)}\widetilde{\varphi}.
$
Now, we need the following result  \cite{jk1995, h2006}.

\begin{lemma}
\label{ext-lem-neg-new}
For $-\frac12<s<\frac12$ and $g\in H^s(\rr)$, we have 
$
\|\chi_{(0,\infty)}g\|_{H^s(\rr)}
\le
c_s\|g\|_{H^s(\rr)}.
$
\end{lemma}
\nin
Applying this result twice, we obtain
$
\|\widetilde{\varphi}_1\|_{H^{\frac{s}{3}}(\rr)}
\lesssim
\|\varphi\|_{H^{\frac{s}{3}}(0,T)}.
$
Therefore, since $\widetilde{\varphi}_1$
is an extension of $\varphi$ from $(0,T)$ to $\rr$ and is  compactly support in $(0,T)$, we can apply
the reduced pure ibvp estimate \eqref{Rob-reduced-pure-ibvp-B-es} 
to get the string of  inequalities 
\begin{align*}
\|
S_\alpha [0,\varphi;0]\|_{X^{s,b,\theta}_{\alpha ,\rr^+\times(0,T)}}
=
\|
S_\alpha [0,\widetilde{\varphi}_1;0]\|_{X^{s,b,\theta}_{\alpha ,\rr^+\times(0,T)}}
\le
\|
S_\alpha [0,\widetilde{\varphi}_1;0]\|_{X^{s,b,\theta}_{\alpha ,\rr^+\times(0,2)}}
\lesssim
\,\| \widetilde{\varphi}_1\|_{H^{\frac{s}{3}}(\rr)}
\lesssim
\|\varphi\|_{H^{\frac{s}{3}}(0,T)},
\end{align*}
for 
$-\frac32 < s< \frac32$, if  $\gamma_2\le 0$,
and 
for $0\le s<\frac32$, if $\gamma_2>0$,
which gives the desired estimates \eqref{pure-ibvp-est-B-kdv}. 
\,\,
$\square$

\vskip.05in
\nin
Next, we decompose the forced ibvp with zero initial and boundary data 
\eqref{inhomo-ibvp-kdv}
 to a forced ivp with  zero initial data and a
 companion pure ibvp. 
 The forced ivp  with zero initial data is given by
\begin{subequations}
\label{inhomo-ivp-kdv}
\begin{align}
\label{inhomo-ivp-kdv:eqn}
&W_t+  \alpha W_{xxx}
= 
w(x, t) 
\in 
X^{s,-b,\theta-1}_\alpha 
\cap
Y^{s,-b}_\alpha ,
\\
\label{inhomo-ivp-kdv:ic}
& W(x,0)
=
0, 
\end{align}
\end{subequations}
where $w \in X^{s,-b,\theta-1}_\alpha \cap Y^{s,-b}_\alpha $ is an extension of the forcing $f_2 \in X_{\alpha ,\rr^+\times(0,T)}^{s,-b,\theta-1}\cap Y_{\alpha ,\rr^+\times(0,T)}^{s,-b}$ such that
$
\|w\|_{X^{s,-b,\theta-1}_\alpha }
+
\|w\|_{Y^{s,-b}_\alpha }
\leqslant
2 (\|f_2\|_{X_{\alpha ,\rr^+\times(0,T)}^{s,-b,\theta-1}}
+
\|f_2\|_{Y_{\alpha ,\rr^+\times(0,T)}^{s,-b}}
).
$
Its solution,  by the Duhamel formula, is
\begin{align}
\label{inhomo-ivp:sln-kdv}
W(x,t) 
\doteq
S_\alpha \big[0; w\big](x, t)
=&
\frac{1}{2\pi} \int_{\xi \in \rr} \int_{t'=0}^t  e^{i\xi x+i\alpha \xi^3(t-t')}
\widehat w(\xi, t') dt' d\xi
\\
\label{inhomo-ivp:sln-d-kdv}
=&
\int_{t'=0}^t  S_\alpha \big[w(\cdot, t'); 0\big](x, t-t') dt',
\end{align}
where $\widehat{w}$ is the $x$-Fourier transform of $w$, 
and  $S_\alpha \big[w(\cdot, t'); 0\big]$ in the solution formula \eqref{inhomo-ivp:sln-d-kdv} is the solution of  the ivp \eqref{homo-ivp-kdv} with initial data $V_0(\cdot)=w(\cdot, t')$. For this problem we
have the following result.
\begin{proposition}
\label{pro-f-ivp-kdv}
The solution $S_\alpha \big[0; w\big]$ 
of ivp \eqref{inhomo-ivp-kdv}  satisfies the following
estimate
\begin{equation}
\label{forced-ivp-bourgain-est-kdv}
\|\psi\, S_\alpha \big[0; w\big]\|_{X^{s,b,\theta}_\alpha }
\lesssim
\|w\|_{X^{s,-b,\theta-1}_\alpha },
\quad
s\in\rr,
\,\,
0< b<\frac12<\theta <1.
\end{equation}
Also,  satisfies the following time regularity estimates
\begin{align}
\label{inhomo-ivp-est-time-kdv}
\sup_{x\in\rr}\|\psi S_\alpha [0;w](x)\|_{H^{\frac{s+1}{3}}(\rr)}
\lesssim
\begin{cases}
\|w\|_{X^{s,-b}_\alpha },
\quad
-1
\le
 s  \le \frac12,
\,\,
0< b<\frac12,
\\
\|w\|_{X^{s,-b}_\alpha }
+
\|w\|_{Y^{s,-b}_\alpha },
\quad
s\in\rr,
 \,\,
 0< b<\frac12,
\end{cases}
\\
\label{inhomo-ivp-est-time-kdv-der}
\sup_{x\in\rr}\|\psi \p_xS_\alpha [0;w](x)\|_{H^{\frac{s}{3}}(\rr)}
\lesssim
\begin{cases}
\|w\|_{X^{s,-b}_\alpha },
\quad
0
\le
 s  \le \frac32,
\,\,
0< b<\frac12,
\\
\|w\|_{X^{s,-b}_\alpha }
+
\|w\|_{Y^{s,-b}_\alpha },
\quad
s\in\rr,
 \,\,
 0< b<\frac12.
\end{cases}
\end{align}
\end{proposition}

\nin
{\bf Proof of Proposition \ref{pro-f-ivp-kdv}.} 
Here we prove estimates \eqref{forced-ivp-bourgain-est-kdv} and \eqref{inhomo-ivp-est-time-kdv}.
The proof of estimate \eqref{inhomo-ivp-est-time-kdv-der} is similar to that of estimate \eqref{inhomo-ivp-est-time-kdv}.
\vskip.05in
\noindent
\underline{Proof of estimate \eqref{forced-ivp-bourgain-est-kdv}.}
For  the $X^{s,b}_\alpha$ part,  we use the inequality
\begin{align}
\label{elementary-est-kdv}
\|
\psi\, S_\alpha\big[0; w\big]
\|_{X^{s,b}_\alpha }^2
\lesssim
\|w\|_{X^{s,b-1}_\alpha }^2
+
\int_{\rr}(1+|\xi|)^{2s}\Big(\int_{\tau\in\rr}
\frac{|\widehat{w}(\xi,\tau)|}{1+|\tau-\alpha \xi^3|}d\tau\Big)^2d\xi,
\,\,
0< b<1,
\end{align}
whose proof for  $\alpha=1$ is given in \cite{fhy2020} (see estimate (2.17)). The proof for $\alpha\neq 1$ is similar.
Then, using the fact $0< b<\frac12$  and applying the Cauchy--Schwarz inequality in $\tau$-integral we get it bounded by $\|w\|_{X^{s,-b}_\alpha }$
(also, for $\alpha=1$ see Lemma 5.6 in  \cite{h2006}).
 For the $D^\theta$ part we use the following inequality
\begin{align*}
\|
\psi\, S_\alpha\big[0; w\big]
\|_{D^\theta}^2
\doteq&
\int_{\rr}
\int_{-1}^1
(1+|\tau|)^{2\theta}
|
\mathcal{F}\big(
\psi S_\alpha\big[0; w\big]
\big)(\xi,\tau)|^2
d\xi
d\tau
\\
\lesssim&
\int_{\rr}
\int_{-1}^1
(1+|\tau|)^{2\theta-2}
|
\widehat{
w}
(\xi,\tau)|^2
d\xi
d\tau
=
\|w\|_{D^{\theta-1}}^2,
\end{align*}
whose proof is again  similar to the proof of \eqref{elementary-est-kdv}. 
In the formula above, $\mathcal{F}$ denotes full Fourier transform.

\vskip.05in
\nin
\underline
{Proof of estimate \eqref{inhomo-ivp-est-time-kdv}.} For this,
we use 
 the decomposition
\begin{subequations}
\label{f-kdv-decomp}
\begin{align}
\label{Tfg2-term-kdv}
\psi(t)S_\alpha \big[0; w\big](x,t)
=&
 \frac{i}{4\pi^2}
 \psi(t)
\int_{\rr^2}
e^{i(\xi x+\tau t)}
\frac{1-\psi(\tau-\alpha \xi^3)}{\tau-\alpha \xi^3}
\widehat{w}(\xi,\tau)  d\tau  d\xi \\
\label{Tfg3-term-kdv}
-&
 \frac{i}{4\pi^2}
 \psi(t)
\int_{\rr^2}
e^{i(\xi x+\alpha \xi^3t)}
\frac{1-\psi(\tau-\alpha \xi^3)}{\tau-\alpha \xi^3}
\widehat{w}(\xi,\tau)   d\tau  d\xi \\
\label{Tfg4-term-kdv}
+&
 \frac{i}{4\pi^2}
  \psi(t)
\int_{\rr^2}
  e^{i(\xi x+\alpha \xi^3 t)} 
 \frac{\psi(\tau-\alpha \xi^3)[e^{i(\tau-\alpha \xi^3)t} - 1]}{\tau-\alpha \xi^3}
\widehat{w}(\xi,\tau)   d\tau  d\xi,
\end{align}
\end{subequations}
and estimate each term separately like in \eqref{elementary-est-kdv}
(also,  for $\alpha=1$ see Lemma 5.6 in  \cite{h2006}).
\,\,
$\square$

\vskip.05in
\nin
The companion pure ibvp from the decomposition of
ibvp \eqref{inhomo-ibvp-kdv}  is
\begin{subequations}
\label{forced-pure-ibvp-kdv}
\begin{align}
\label{forced-pure-ibvp:eqn-kdv}
&v_t+\alpha v_{xxx} 
= 
0, 
\\
\label{forced-pure-ibvp:ic-kdv}
&v(x,0)
= 
0,  
\\
\label{forced-pure-ibvp:bc-kdv}
&\p_xv(0,t)+\gamma_2v(0,t)
=
-
[\p_xW(0,t)+\gamma_2W(0,t)]
\doteq 
W_0(t)
\in
H^{\frac{s}{3}}(0,T),
\end{align}
\end{subequations}
where by the KdV time estimate \eqref{inhomo-ivp-est-time-kdv} and \eqref{inhomo-ivp-est-time-kdv-der}, $W_0\in H^{\frac{s}{3}}(0,T)$.
Its solution
$
v(x, t) 
\doteq
S_\alpha \big[0, W_0; 0\big](x, t)
$
is defined by UTM formula \eqref{Rob-kdv-utm-v}
and satisfies 
 an estimate like the one in Proposition \ref{pro-p-ibvp-kdv} with $\varphi$ being replaced
with $W_0$.

\vskip0.05in
\nin
{\bf Proof of estimate \eqref{MB-forced-linear-est2-Rob}.} 
 For this estimate
using the superposition principle  we express the Fokas solution formula  $S_\alpha\big[v_0, \varphi_2; f_2\big]$
for the  forced ibvp \eqref{Rob-LKdV-v}  as follows
\begin{equation*}
S_\alpha\big[v_0, \varphi_2; f_2\big] 
=
\psi(t)
S_\alpha\big[
V_0; 0\big]
+
S_\alpha\big[0, \varphi; 0\big]
+
\psi(t)
S_\alpha\big[0; w\big]
+
 S_\alpha\big[0,  W_0; 0\big],
 \,\,
 x>0,
 \,\,
 t\in(0,T).
\end{equation*}
Applying the triangular inequality 
for the norm  $\|\cdot\|_{X^{s,b,\theta}_{\alpha,\rr^+\times(0,T)}}$, and 
 using Propositions \ref{pro-ivp-kdv}--\ref{pro-f-ivp-kdv}
we get
\begin{align}
\label{Rob-wp-fin}
\big\|
S_\alpha\big[v_0, \varphi_2; f_2\big] 
\big\|_{X^{s,b,\theta}_{\alpha,\rr^+\times(0,T)}}
\lesssim
\big\|V_0\big\|_{H^{s}(\rr)}
+
\|
\varphi
\|_{H_t^{\frac{s}{3}}(0,T)} 
+
\|
w
\|_{X_\alpha^{s,-b,\theta-1}}
+
\|
W_0
\|_{H_t^{\frac{s}{3}}(0,T)}.
\end{align}
Finally,  considering that $\varphi(t)=\varphi_2(t)
-
[\p_xV(0,t)+\gamma_2V(0,t)]$, 
$
W_0(t)
=-
[\p_xW(0,t)+\gamma_2W(0,t)],
$ 
using the time regularity estimates  
\eqref{homo-ivp-est-time-kdv}, \eqref{homo-ivp-est-time-kdv-der},
\eqref{inhomo-ivp-est-time-kdv}, 
\eqref{inhomo-ivp-est-time-kdv-der}
and using the inequalities
$\|V_0\|_{H^{s}(\rr))}\le  2\|v_0\|_{H^{s}(\rr^+)}$, 
$
\|w\|_{X_\alpha^{s,-b,\theta-1}}
+
\|w\|_{Y_\alpha^{s,-b}}
\le
2 \|f_2\|_{X_{\alpha,\rr^+\times(0,T)}^{s,-b,\theta-1}}
+
2 \|f_2\|_{Y_{\alpha,\rr^+\times(0,T)}^{s,-b}}
$,
from \eqref{Rob-wp-fin}
we complete  the proof of  estimate \eqref{MB-forced-linear-est2-Rob}.
\,\,
$\square$

%
%
%
%
%
%
%
%
%
%
%
\section{Linear estimates for Dirichlet problem
--
Proof of Theorem \ref{thm-MB-forced-linear-Dir}
}
\label{sec:linear-Dir}
\setcounter{equation}{0}
\subsection{Reduced pure ibvp for Dirichlet problem}
We begin with the reduced pure  ibvp  
corresponding to the $v$-equation ($\alpha$-KdV),
which is our main tool in proving the linear 
estimates for the  Dirichlet problem.
This problem reads as follows
\begin{subequations}
\label{kdv-reduced}
\begin{align}
\label{kdv-reduced-eqn}
&\p_tv+\alpha\p_x^{3}v=0,
\quad
0<x<\infty,
\,\,
0<t<2,
\\
\label{kdv-reduced-ic}
&v(x,0)
=
0,
\\
\label{kdv-reduced-bc}
&v(0,t)
=
h(t),
\,\,\,  \text{supp}(h)\subset (0, 2).
\end{align}
\end{subequations}
Since the boundary data $h$ are  compactly supported in time interval 
$(0, 2)$, we see the time transform  $\widetilde h$ over the interval $(0, 2)$
is the following Fourier transform over $\rr$
\begin{equation*}
\widetilde h(\xi,2)
=
\int_0^2e^{-i\alpha\xi\tau}h(\tau)d\tau
=
\int_\rr e^{-i\alpha\xi\tau}h(\tau)d\tau
=
\widehat h(\alpha\xi).
\end{equation*}
 Thus, the Fokas solution formula  \eqref{kdv-utm-v} for this
 ibvp  on $\rr^+\times (0, 2)$ reads as follows
\begin{align}
\label{pure-kdv-sln}
v(x,t)
\doteq
S_\alpha[0,h;0]
=
-
\frac{3\alpha}{2\pi}\int_{\p D^+} e^{i\xi x+i\alpha\xi^3t}\xi^2\,
 \widehat h(\alpha\xi^3)d\xi.
\end{align}
Next, using formula \eqref{pure-kdv-sln} 
we prove the following key result.
\begin{theorem}
[\b{Estimates for pure ibvp of linear $\alpha$-KdV (Dirichlet)}]
\label{reduced-thm-kdv}
For boundary data $h$ a test function with
$\text{supp}(h)\subset (0, 2)$
 the solution to  reduced pure linear ibvp \eqref{kdv-reduced} satisfies the  Bourgain space estimate
\begin{align}
\label{reduced-pure-ibvp-B-es}
\|
S_\alpha[0,h;0]
\|_{X^{s,b,\theta}_{\alpha,\rr^+\times(0,2)}} 
\leq
c_{s,b,\theta,\alpha}
\|h\|_{H_t^{\frac{s+1}{3}}(\rr)},
\,\,
s>-\frac32,
\,\,
0\le b<\frac12,
\,\,
0\le \theta< 1+\frac13s.
\end{align}
\end{theorem}

\nin
{\bf Proof of Theorem \ref{reduced-thm-kdv}.}
Like in the Robin problem, 
parametrizing  the right (left) contour of $\p D^+$ via
$[0,\infty)\ni \xi\to a \xi$ ($a^2\xi$) with 
$a
=
e^{i\frac{\pi}{3}}
$,
we write the solution as 
$
v
\simeq
v_{r}+v_{\ell},
$
where
\begin{align}
\label{def-vr-vl}
v_r(x,t)
=
\int_0^\infty e^{ia\xi x} e^{-i\alpha\xi ^3t} \xi ^{2}
\,
\widehat h(-\alpha\xi ^3)d\xi
\quad
\text{and}
\quad
v_\ell(x,t)
=
\int_0^\infty e^{ia^2\xi x}  e^{i\alpha\xi^3t}
\,
\xi^{2}
\widehat h(\alpha\xi^3)d\xi.
\end{align}
Here we only estimate $v_{r}$. 
The estimation  of $v_{\ell}$ is similar.
We begin by we splitting  the integration 
in formula \eqref{def-vr-vl}
 for $\xi$ near 0 and  for $\xi$ near $\infty$.
 Thus, we have
$
v_r
=
v_0
+
v_1,
$
where
\begin{align*}
v_0(x,t)
\doteq
\int_0^1 
\hskip-0.05in
e^{ia\xi x} e^{-i\alpha\xi ^3t} \xi ^{2}
\widehat h(-\alpha\xi ^3)d\xi
\,\,\,
\text{and}
\,\,\,
v_1(x,t)
\doteq
\int_1^\infty 
\hskip-0.05in
e^{ia\xi x} e^{-i\alpha\xi ^3t} \xi ^{2}
\widehat h(-\alpha\xi ^3)d\xi,
\,\,
x>0,
\,\,
t\in\rr.
\end{align*}
{\bf Proof of estimate \eqref{reduced-pure-ibvp-B-es} for $v_0$.}
Using $\varphi_1(x)$, which is a smooth version of 
$|x|$ defined by \eqref{phi1-x-even-ext}, 
we extend $v_0$
as an ``almost" even function
of $x$ from $\rr^+\times (0,2)$ to $\rr\times\rr$ 
as follows 
\begin{equation*}
V_0(x,t)
\doteq
\int_0^1 e^{ia\xi \varphi_1(x)}e^{-i\alpha\xi ^3t}
\xi^2
\,\widehat {h}(-\alpha\xi ^3)d\xi,
\quad
x\in\rr,
\,\,
t\in\rr.
\end{equation*}
Then,  using  \eqref{Rob-est-sb-norm-fact-no3} we get
the following  $L^2$ estimate for the extension
$\psi_4(t)V_0$
\begin{align}
\label{v0-L2-est}
\|\psi_4 V_0\|^2_{X^{s,b,\theta}_{\alpha}}
\lesssim
\|\psi_4 V_0\|^2_{L^2_{x,t}}
+
\|\p_x^{n_1}[\psi_4 V_0]\|^2_{L^2_{x,t}}
+
\|\p_t^{n_2}[\psi_4 V_0]\|^2_{L^2_{x,t}}
+
\|\p_x^{n_1}\p_t^{n_2}[\psi_4 V_0]\|^2_{L^2_{x,t}},
\end{align}
where 
$n_1
=
n_1(s,b)
\doteq
6\lfloor b\rfloor+2\lfloor |s|\rfloor+8
$
and
$
n_2
=
n_2(b,\theta)
\doteq
2\lfloor b\rfloor+2\lfloor\theta\rfloor+4.
$
Also, working like in the case of estimate \eqref{Rob-L2-V0-est-kdv}, we get
\begin{equation}
\label{L2-V0-est-kdv}
\|\p_x^{n_1}\p_t^{n_2}[\psi_4\,V_0]\|_{L^2_{x,t}}^2
\le
C_{n_1, n_2}
\int_{-1}^{0}
|
\widehat{h}(\alpha\tau)
|^2
d\tau,
\quad
\text{for any nonnegative integers}
\,\,
n_1,n_2.
\end{equation}
Finally, combining \eqref{v0-L2-est} and \eqref{L2-V0-est-kdv} 
 we get 
\begin{align*}
\|\psi_4 V_0\|^2_{X^{s,b,\theta}_{\alpha}}
\lesssim
\int_{-1}^{0}
(1+|\alpha\tau|)^{\frac{-2(s+1)}{3}}
(1+|\alpha\tau|)^{\frac{2(s+1)}{3}}
|
\widehat{h}(\alpha\tau)
|^2
d\tau
\lesssim
\|h\|_{H_t^{\frac{s+1}{3}}}^2,
\,\,
s\in\rr,b\ge0,\theta \ge 0,
\end{align*}
which  implies
$\|
v_0
\|_{X^{s, b,\theta}_{\alpha,\rr^+\times (0,2)}}\le \|\psi_4 V_0\|_{X^{s,b,\theta}_{\alpha}}\lesssim \|h\|_{H_t^{\frac{s+1}{3}}}$ and
gives  the 
desired estimate \eqref{reduced-pure-ibvp-B-es} for $v_0$.
\,\,
$\square$

\vskip0.05in
\nin
{\bf Proof of estimate \eqref{reduced-pure-ibvp-B-es} for $v_1$.} Since $a=a_R+ia_I=\frac12+i\frac{\sqrt3}{2}$, we rewrite $v_1$ as 
\begin{equation*}
\label{v1-formula}
v_1(x,t)
=
\int_1^\infty 
 e^{-i\alpha\xi ^3t}
 e^{ia_R\xi x}e^{-a_I\xi x}
  \xi ^{2}
\widehat h(-\alpha\xi ^3)d\xi,
\quad
x>0,
\,\,
t\in\rr.
\end{equation*}
Then, using  the identity
$
\xi
[
 e^{ia_R\xi x}e^{-a_I\xi x}
]
=
\frac{\p_x}{ia}
[
 e^{ia_R\xi x}e^{-a_I\xi x}
]
$
and the fact that  $e^{-a_I \xi x}$ is exponentially 
decaying in $\xi$ for $x>0$, we move
the $\p_x$-derivative outside the integral sign 
to rewrite $v_1$ as follows
\begin{equation*}
v_1(x,t)
=
\frac{1}{ia}
\p_x
\int_1^\infty e^{-i\alpha\xi^3t}e^{ia_R\xi x}e^{-a_I \xi x} \xi 
\,
\widehat h(-\alpha\xi ^3)d\xi,
\quad
x>0,
\quad
t\in\rr.
\end{equation*}
Furthermore, using   
the one-sided  $C^\infty$ cutoff function $\rho(x)$, which is
defined by 
\eqref{Rob-rho-def},
we extend $v_1$  from $\rr^+\times \rr$ to $\rr^2$ 
via the formula below 
\begin{align*}
V_1(x,t)
\doteq&
\frac{1}{ia}
\p_x
\int_1^\infty 
e^{-i\alpha\xi^3t}e^{ia_R\xi x}e^{-a_I \xi x} \rho(a_I\xi x)
\xi
\,
\widehat h(-\alpha\xi ^3)d\xi,
\quad
x\in\rr,
\quad
t\in\rr.
\end{align*}
Since 
$
\|v_1\|_{X^{s,b,\theta}_{\alpha,\rr^+\times(0,2)}}
\le 
\|V_1\|_{X^{s,b,\theta}_{\alpha}},
$
 to prove  estimate \eqref{reduced-pure-ibvp-B-es} for $v_1$ it suffices to show that
\begin{align}
\label{reduced-ibvp-Bourgain-est-1}
\|V_1\|_{X^{s,b,\theta}_{\alpha}} 
\le
c_{s,b,\theta,\alpha}
\|h\|_{H_t^{\frac{s+1}{3}}(\rr)}.
\end{align}
We do this for both parts of modified Bourgain norm
 shown below
\begin{align*}
\|V_1\|_{X^{s,b,\theta}_{\alpha}}^2
=
\|V_1\|_{X^{s,b}_{\alpha}}^2
+
\|V_1\|_{D^\theta}^2,
\,\,
\text{with}
\,\,
\|V_1\|_{D^\theta}^2
\doteq
\int_{-\infty}^\infty
\int_{-1}^1
(1+|\tau|)^{2\theta}
|\widehat{V}_1(\xi,\tau)|^2
d\xi
d\tau.
\end{align*}
\underline{Estimate of the Bourgain $X^{s,b}_{\alpha}$ norm for $V_1$}.
Using the Schwartz function  $\eta(x)$ defined by \eqref{Rob-def-eta} and 
making the change of variables $\tau=-\xi^3$ we see that $V_1$ takes the form
\begin{equation}
\label{eqn-V1-1}
V_1(x,t)
\simeq
\p_x
\int_{-\infty}^{-1} e^{i\alpha\tau t}
\eta(-a_I\tau^{1/3}x)
\,
\tau^{-\frac13}
\,
\widehat h(\alpha\tau)
d\tau,
\quad
x\in\rr,
\quad
t\in\rr.
\end{equation}
Now, letting $V_{1,\alpha}(x,t)=V_1(x,\frac{t}{\alpha})$ and using 
$
\widehat{V}_{1,\alpha}(\xi,\tau)
=
\alpha
\widehat{V}_{1}(\xi,\alpha\tau),
$
we get 
$
\|V_1\|_{X_\alpha^{s,b}}
\lesssim
\|V_{1,\alpha}\|_{X^{s,b}_1} 
$
($b\ge 0$),
which reduces estimate \eqref{reduced-ibvp-Bourgain-est-1} 
to the inequality
\begin{align}
\label{reduced-ibvp-Bourgain-est-2}
\|V_{1,\alpha}\|_{X^{s,b}_{1}} 
\le
c_{s,b,\alpha}
\|h\|_{H_t^{\frac{s+1}{3}}(\rr)},
\end{align}
where, by \eqref{eqn-V1-1},
\begin{equation}
\label{eqn-V1-alpha}
V_{1,\alpha}(x,t)
\simeq
\p_x
\int_{-\infty}^{-1} e^{i\tau t}
\eta(-a_I\tau^{1/3}x)
\,
\tau^{-\frac13}
\,
\widehat h(\alpha\tau)
d\tau,
\quad
x\in\rr,
\quad
t\in\rr.
\end{equation}
Now, we see that Bourgain space estimate \eqref{reduced-ibvp-Bourgain-est-2} follows from the next result.
\begin{lemma}
\label{reduced-pibvp-thm1}
For any $\varepsilon>0$, if $s\geq -\frac32-\varepsilon$, 
and $b\ge 0$,
then  the function $V_{1,\alpha}(x,t)$, which is  defined by \eqref{eqn-V1-alpha}
satisfies the Bourgain  spaces estimate
\begin{align}
\label{reduced-ibvp-Bourgain-est}
\|V_{1,\alpha}\|_{X^{s,b}_{1}} 
\leq
c_{s,b,\varepsilon,\rho,a,\alpha} \|h\|_{H_t^{\frac{s+3b-\frac12+\varepsilon}{3}}(\rr)}.
\end{align}
\end{lemma}
\noindent
In fact, 
choosing 
$
\varepsilon
=
3(\frac12-b),
$
which is possible if $b<1/2$, we get $(s+3b-\frac12+\varepsilon)/3=\frac{s+1}{3}$ and estimates \eqref{reduced-ibvp-Bourgain-est} gives  \eqref{reduced-ibvp-Bourgain-est-2}. Combining this with 
$
\|V_1\|_{X_\alpha^{s,b}}
\lesssim
\|V_{1,\alpha}\|_{X^{s,b}_1}, 
$
gives 
the  desired  estimate \eqref{reduced-ibvp-Bourgain-est-1}.

\vskip0.05in
\nin
{\bf Proof of Lemma \ref{reduced-pibvp-thm1}.}
This proof is similar to the proof of Lemma \ref{Rob-reduced-pibvp-thm1}. By  \eqref{eqn-V1-alpha} we have the following formula
for the time Fourier transform 
$
\widehat{V}_{1,\alpha}^t(x,\tau)
\simeq
\p_x
\big[
\chi_{(-\infty,-1)}(\tau)
\eta(-a_I\tau^{1/3}x)
\,
\tau^{-1/3}
\,
\widehat h(\alpha\tau)
\big],
$
which gives us the full Fourier transform 
\begin{align}
\label{v1-xt-FT}
\widehat{V}_{1,\alpha}(\xi,\tau)
\simeq
\xi
\chi_{(-\infty,-1)}(\tau)
\,
F(\xi,\tau)
\,
\tau^{-1/3}
\,
\widehat h(\alpha\tau),
\end{align}
where 
$
F(\xi,\tau)
\doteq
\int_{\rr}
e^{-i\xi x}\eta(-a_I\tau^\frac13 x)
dx
$ 
is defined by \eqref{Rob-v1-xt-FT}. Now,
using estimate \eqref{Rob-v1-mult-ine-1},
we get
\begin{align*}
\|V_{1,\alpha}\|_{X^{s,b}_1}^2
=&
\int_{-\infty}^{-1}
\Big(
\int_{\rr}
(1+|\xi|)^{2s}
(1+|\tau-\alpha\xi^3|)^{2b}
\xi^2
|F(\xi,\tau)|^2
d\xi
\Big)
|\tau^{-1/3}\widehat{h}(\alpha\tau)|^2
d\tau
\\
\lesssim&
\int_{-\infty}^{-1}
|\tau|^{\frac{2s+6b+1+2\varepsilon}{3}}
|\tau^{-1/3}
\widehat h(\alpha\tau)
|^2
d\tau
=
\int_{-\infty}^{-1}
\Big|
\tau^{\frac{s+3b-\frac12+\varepsilon}{3}}
\widehat h(\alpha\tau)
\Big|^2
d\tau
\lesssim
\|h\|_{H_t^{\frac{s+3b-\frac12+\varepsilon}{3}}}^2,
\end{align*}
which is  the desired estimate \eqref{reduced-ibvp-Bourgain-est} 
in Lemma \ref{reduced-pibvp-thm1}.
\,\,
$\square$

\vskip0.05in
\nin
\underline{Estimate of $D^\theta$ norm for $V_1$}. 
As before, letting $V_{1,\alpha}(x,t)\doteq V_1(x,\frac{t}{\alpha})$ and 
using the  identity
$
\widehat{V}_{1,\alpha}(\xi,\tau)
=
\alpha
\widehat{V}_{1}(\xi,\alpha\tau),
$
we get
$
\|V_1\|_{D^\theta}
\lesssim
\|V_{1,\alpha}\|_{D^\theta}
$, if $\theta\ge0$.
Using Fourier transform formula \eqref{v1-xt-FT} we get
\begin{align*}
\|V_{1,\alpha}\|_{D^\theta}^2
=
\int_{-\infty}^{-1}
\int_{-1}^1
(1+|\tau|)^{2\theta}
|
\xi
F(\xi,\tau)
\tau^{-1/3}
\widehat h(\alpha\tau)|^2
d\xi
d\tau,
\end{align*}
where 
$
F(\xi,\tau)
$
is defined in \eqref{Rob-v1-xt-FT}.
Also, applying Lemma \ref{Rob-F-bound-lem} with $n=0$, we get
estimate
$
|F(\xi,\tau)|
\lesssim
|\tau|^{-1/3},
$
which we use in the last relation  to  integrate $\xi$ and
to obtain the inequality 
\begin{align}
\label{V1-D}
\|V_{1,\alpha}\|_{D^\theta}^2
\lesssim
\int_{-\infty}^{-1}
(1+|\tau|)^{2\theta}
|
\tau^{-2/3}
\widehat h(\alpha\tau)|^2
d\tau
\lesssim
\int_{-\infty}^{-1}
(1+|\tau|)^{2\theta-4/3}
|
\widehat h(\alpha\tau)|^2
d\tau
\lesssim
\|h\|_{H_t^{\theta-2/3}}^2.
\end{align}
Choosing  $\theta-2/3\le (s+1)/3$ or $\theta\le 1+s/3$
gives
$ \|V_1\|_{D^\theta}
\lesssim
\|V_{1,\alpha}\|_{D^\theta} \lesssim \|h\|_{H_t^{(s+1)/3}}$,
i.e. \eqref{reduced-ibvp-Bourgain-est-1}.\, $\square$

%
%
%
%
%
%
%
%
%
%
%
%
\subsection{Proof of Linear estimates for Dirichlet  ibvp
(Theorem \ref{thm-MB-forced-linear-Dir})}
It suffices to prove linear estimate \eqref{MB-forced-linear-est2-Dir}
for  $\alpha-KdV$, the second equation of our MB system.
Linear estimate \eqref{MB-forced-linear-est1-Dir} for the first equation is obtained by letting
$\alpha=1$, replacing data $(v_0,h_0)$ by $(u_0,g_0)$, and replacing forcing $f_2$ by $f_1$.
We begin by decomposing the forced linear
 $\alpha$-KdV ibvp \eqref{LKdV-v}
into the homogeneous ibvp 
\begin{subequations}
\label{Dir-homo-ibvp-kdv}
\begin{align}
\label{Dir-homo-ibvp-kdv:eqn}
&v_t+\alpha v_{xxx}
=
0,
\\
\label{Dir-homo-ibvp:ic}
&v(x,0) 
=
v_0(x)\in H^{s}(\rr^+),
\\
\label{Dir-homo-ibvp-kdv:bc}
&v(0,t)
=
h_0(t)\in H^{\frac{s+1}{3}}(0,T),
\end{align}
\end{subequations}
and the companion forced ibvp with zero initial and boundary  data
\begin{subequations}
\label{Dir-inhomo-ibvp-kdv}
\begin{align}
\label{Dir-inhomo-ibvp-kdv:eqn}
&v_t+\alpha v_{xxx}
=
f_2(x,t)
\in 
X_{\alpha,\rr^+\times(0,T)}^{s,-b,\theta-1}
\cap
Y_{\alpha,\rr^+\times(0,T)}^{s,-b},
\\
\label{Dir-inhomo-ibvp-kdv:ic}
&v(x,0) 
= 
0, 
\\
\label{Dir-inhomo-ibvp-kdv:bc}
&v(0,t)
= 
0.
\end{align}
\end{subequations}
Their solutions, by linearity, satisfy the relation 
$
S_\alpha\big[v_0, h_0; f_2\big]
=
S_\alpha\big[v_0, h_0; 0\big]
+
S_\alpha\big[0, 0; f_2\big].
$
Then, 
we  decompose the homogeneous ibvp \eqref{Dir-homo-ibvp-kdv}
into the  ivp \eqref{homo-ivp-kdv},  and the pure ibvp 
\begin{subequations}
\label{Dir-pure-ibvp-kdv}
\begin{align}
\label{Dir-pure-ibvp-kdv:eqn}
&v_t+\alpha v_{xxx} 
= 
0, 
\\
\label{Dir-pure-ibvp-kdv:ic} 
&v(x,0)
= 
0, 
\\
\label{Dir-pure-ibvp-kdv:bc}
&v(0,t)
= 
h_0(t)
-
V(0,t)
\doteq 
h_1(t)
\in H^{\frac{s+1}{3}}(0,T),
\end{align}
\end{subequations}
whose  solution
$
v
\doteq
S_\alpha \big[0, h_1; 0\big]
$, defined by UTM formula \eqref{kdv-utm-v}, satisfies the following result.
\begin{proposition}
\label{Dir-pro-p-ibvp-kdv}
For $b\in [0,\frac12)$,
 $\frac12<\theta
 \le
1+\frac13s$, the solution to ibvp 
  \eqref{Dir-pure-ibvp-kdv} satisfies  the estimate
\begin{align}
\label{Dir-pure-ibvp-est-B-kdv}
\|
S_\alpha [0,h_1;0]\|_{X^{s,b,\theta}_{\alpha ,\rr^+\times(0,T)}}
\le 
C_{s}\, \|h_1\|_{H^{\frac{s+1}{3}}(0,T)},
\quad
-\frac32 < s< \frac12.
\end{align}
In addition, if $h_1(0)=0$, then estimates \eqref{Dir-pure-ibvp-est-B-kdv} holds  for $\frac12<s<2$.
\end{proposition}

\nin
{\bf Proof of Proposition \ref{Dir-pro-p-ibvp-kdv}.} Like in the Robin problem, we transform the pure ibvp \eqref{Dir-pure-ibvp-kdv} to the reduced pure ibvp \eqref{kdv-reduced}. For this, we extend the boundary data $h_1(t)$ to  a function $\widetilde{h}_1$ on $\rr$ supported in $(0,2)$ and such that $\|\widetilde{h}_1\|_{H^{(s+1)/3}(\rr)}\lesssim \|h_1\|_{H^{(s+1)/3}(0,T)}$,
via using compatibility condition \eqref{comp-cond},  Lemma \ref{ext-lem-neg-new}
and
the following result \cite{LM-book1972, T-book1996, fhm2017}. 
\begin{lemma}
\label{Extension-lemma}
For a  general function $h^*(t)\in H_t^s(0,2)$, $s\ge 0$,  let 
\begin{align*}
\widetilde h^*(t)
\doteq
\begin{cases}
h^*(t),
\quad
t\in(0,2),
\\
0,
\quad
\text{elsewhere}.
\end{cases}
\end{align*} 
If $0\leq s<1/2$,  then for some $c>0$, which is depending on $s$, we have
\begin{align}
\label{ext-est-Rob}
\|\widetilde h^*\|_{H_t^s(\rr)}
\leq
c\|h^*\|_{H_t^s(0, 2)}.
\end{align}
If  $1/2<s<1$, then for estimate  \eqref{ext-est-Rob} to hold
we must have the condition  $h^*(0)=h^*(2)=0$.
\end{lemma}

\nin
Since $\widetilde{h}_1$ extends $h_1$ from $(0,T)$ to $\rr$, we have $\|
S_\alpha [0,h_1;0]\|_{X^{s,b,\theta}_{\alpha ,\rr^+\times(0,T)}}=\|
S_\alpha [0,\widetilde{h}_1;0]\|_{X^{s,b,\theta}_{\alpha ,\rr^+\times(0,T)}}$, which combined with estimate \eqref{reduced-pure-ibvp-B-es} implies that 
$
\|
S_\alpha [0,h_1;0]\|_{X^{s,b,\theta}_{\alpha ,\rr^+\times(0,T)}}
\lesssim
\,\| \widetilde{h}_1\|_{H^{\frac{s+1}{3}}(\rr)}
\lesssim
\|h_1\|_{H^{\frac{s+1}{3}}(0,T)}.
$
This completes the proof of Proposition \ref{Dir-pro-p-ibvp-kdv}.
\,\,
$\square$

\vskip.05in
\nin
Next, we decompose the forced ibvp  
\eqref{Dir-inhomo-ibvp-kdv}
 to the forced ivp \eqref{inhomo-ivp-kdv} and the
 companion pure ibvp
\begin{subequations}
\label{Dir-forced-pure-ibvp-kdv}
\begin{align}
\label{Dir-forced-pure-ibvp:eqn-kdv}
&v_t+\alpha v_{xxx} 
= 
0, 
\\
\label{Dir-forced-pure-ibvp:ic-kdv}
&v(x,0)
= 
0,  
\\
\label{Dir-forced-pure-ibvp:bc-kdv}
&v(0,t)
=
- W(0,t)
\doteq 
W_0(t)
\in
H^{\frac{s+1}{3}}(0,T),
\end{align}
\end{subequations}
where, by  the time regularity  estimate \eqref{inhomo-ivp-est-time-kdv}, $W_0\in H^{\frac{s+1}{3}}(0,T)$.
Its solution
$
v
\doteq
S_\alpha \big[0, W_0; 0\big]
$,
defined by the Fokas UTM formula \eqref{kdv-utm-v},
satisfies 
 an estimate like the one in Proposition \ref{Dir-pro-p-ibvp-kdv} with $h_1$ being replaced
with $W_0$.

\vskip0.05in
\nin
{\bf Proof of estimate \eqref{MB-forced-linear-est2-Dir}.} 
First, using the superposition principle,  we express the Fokas solution formula  $S_\alpha\big[v_0, h_0; f_2\big]$
for the  forced ibvp \eqref{kdv-utm}  as follows
\begin{equation*}
S_\alpha\big[v_0, h_0; f_2\big] 
=
\psi(t)
S_\alpha\big[
V_0; 0\big]
+
S_\alpha\big[0, h_1; 0\big]
+
\psi(t)
S_\alpha\big[0; w\big]
+
 S_\alpha\big[0,  W_0; 0\big],
 \,\,
 x>0,
 \,\,
 t\in(0,T).
\end{equation*}
Then, applying the triangular inequality 
and utilizing  Propositions \ref{pro-ivp-kdv},  \ref{Dir-pro-p-ibvp-kdv} and \ref{pro-f-ivp-kdv}
we get
\begin{align}
\label{Dir-wp-fin}
\big\|
S_\alpha\big[v_0, h_0; f_2\big] 
\big\|_{X^{s,b,\theta}_{\alpha,\rr^+\times(0,T)}}
\lesssim
\big\|V_0\big\|_{H^{s}(\rr)}
+
\|
h_1
\|_{H_t^{\frac{s+1}{3}}(0,T)} 
+
\|
w
\|_{X_\alpha^{s,-b,\theta-1}}
+
\|
W_0
\|_{H_t^{\frac{s+1}{3}}(0,T)}.
\end{align}
Finally,  considering 
$
h_1(t)
=
h_0(t)-V(0,t),
$
$W_0(t)=-W(0,t)$,
using the time regularity estimate \eqref{inhomo-ivp-est-time-kdv}
and
using
the inequalities
$\|V_0\|_{H^{s}(\rr))}\le  2\|v_0\|_{H^{s}(\rr^+)}$
and 
$\|w\|_{X_\alpha^{s,-b,\theta-1}}
+
\|
w
\|_{Y_\alpha^{s,-b}}
\le
2 \|f_2\|_{X_{\alpha,\rr^+\times(0,T)}^{s,-b,\theta-1}}
+
2 \|f_2\|_{Y_{\alpha,\rr^+\times(0,T)}^{s,-b}},
$
from \eqref{Dir-wp-fin}
we complete the proof of  estimate \eqref{MB-forced-linear-est2-Dir}.
\,\,
$\square$

%
%
%
%
%
%
%
%
%
%
\section{
Proof of Well-posedness  of  the Majda-Biello System
}
\label{sec:wp}
\setcounter{equation}{0}
Here we  present only the proof of well-posedness for the Robin and Neumann ibvp's (Theorem \ref{thm-MB-half-line-Rob}).
The proof  of well-posedness for the Dirichlet ibvp (Theorem \ref{thm-MB-half-line-Dir}) is similar.

 We begin by recalling the  critical Sobolev exponent
 $s_{c}(\alpha)$, $\alpha>0$,  defined 
in the introduction by \eqref{MB-critical-sob-index}. 
Also, we assume that 
$s_c(\alpha)\le s<\frac32$ when $\gamma_1\le 0$ and  $\gamma_2\le 0$,
and 
$\max\{0, s_c(\alpha)\}\le s<\frac32$ when  
$\gamma_1>0 $ or  $\gamma_2> 0$.
Then, for initial data $u_0\in H^{s}(\rr^+)$, $v_0\in H^{s}(\rr^+)$,  boundary data 
$\varphi_1\in H^{\frac{s}3}(0,T)$, $\varphi_2\in H^{\frac{s}3}(0,T)$,  
where  $0<T<1/2$,
and   $T^*$, (to be chosen later),  satisfying the condition
\begin{equation*}
0<T^*<T<1/2,
\end{equation*}
  we define  the iteration map
  $
  (u,v)
\mapsto 
  \Phi_{T^*}\times\Psi_{T^*}(u,v)
  \doteq (\Phi_{T^*}(u,v), \Psi_{T^*}(u,v))
  $
  which is obtained from the Fokas solution formula
\eqref{Rob-kdv-utm} and \eqref{Rob-kdv-utm-v}
for the forced linear MB  ibvp with the forcing terms
$f_1$, $f_2$ replaced by the nonlinearities, and localized
appropriately. More precisely, we have
\begin{align}
\label{Rob2-iteration-map-T-loc}
\Phi_{T^*} (u,v)
\doteq
S_1\big[u_0, \varphi_1; -\frac12\psi_{2T^*}[\p_x(v^2)]\big]
\quad
\text{and}
\quad
\Psi_{T^*}(u,v)
\doteq
S_{\alpha}\big[v_0, \varphi_2;-\frac12\psi_{2T^*}[\p_x(uv)]\big],
\end{align}
where   $\psi_{T^*}(t)=\psi(t/T^*)$ with $\psi(t)$  being 
localizer  \eqref{time-localizer}. Since,  
$\psi_{2T^*}(t)=1$ 
for $|t|\le T^*$, 
the fixed point of the iteration map \eqref{Rob2-iteration-map-T-loc} is the solution to the MB system ibvp \eqref{Rob-MBsys}.  
Next, we will show that  our iteration map \eqref{Rob2-iteration-map-T-loc} is a contraction in the complete metric space 
$\mathcal{X}_{MN}\subset X^{s,b,\theta}_{1,\rr^+\times(0,T)}
\times
X^{s,b,\theta}_{\alpha,\rr^+\times(0,T)}$ defined by
\begin{align*}
\mathcal{X}_{MN}
=
\{
(u,v):
u\in X^{s,b,\theta}_{1,\rr^+\times(0,T)},
\,\,
v\in X^{s,b,\theta}_{\alpha,\rr^+\times(0,T)}
\,\,
\text{such that}
\,\,
\|u\|_{X^{s,b,\theta}_{1,\rr^+\times(0,T)}}
\le
M,
\,\,
\|v\|_{X^{s,b,\theta}_{\alpha,\rr^+\times(0,T)}}
\le
N
\},
\end{align*}
and having norm
$
\|(u,v)\|_{\mathcal{X}_{MN}}
\doteq
\|u\|_{X^{s,b,\theta}_{1,\rr^+\times(0,T)}}
+
\|v\|_{X^{s,b,\theta}_{\alpha,\rr^+\times(0,T)}},
$
if we choose  $ {T^*}$, $M$ and $N$ appropriately. 
\vskip.05in
\noindent
{\bf $\Phi_{T^*}\times\Psi_{T^*}$ is onto.}
We  estimate the components of the iteration map
separately. Concerning  $\Phi_{T^*}$,
using  the linear estimate \eqref{MB-forced-linear-est1-Rob} with forcing replaced by $-\frac12\psi_{2{T^*}}[\p_x(v^2)]$,   for $0<b\le b^*<\frac12$ we get
\begin{align}
\label{Rob2-onto-map-fin-est-large}
\mathcal{N}_{\Phi}
\doteq
&\|\Phi_ {T^*}(u,v)\|_{X^{s,b,\theta}_{1,\rr^+\times(0,T)}}
\le
\|\Phi_ {T^*}(u,v)\|_{X^{s,b^*,\theta}_{1,\rr^+\times(0,T)}}
\nn
\\
\le&
c
\big(
\|
u_0
\|_{H^{s}(\rr^+)}
+
\|\varphi_1\|_{H_t^\frac{s}{3}(0,T)}
+
\frac12
\big\|
\psi_{2{T^*}}[\p_x(v^2)]
\big\|_{X^{s,-b^*,\theta-1}_{1,\rr^+\times(0,T)}}
+
\frac12
\big\|
\psi_{2{T^*}}[\p_x(v^2)]
\big\|_{Y^{s,-b^*}_{1,\rr^+\times(0,T)}}
\big)
\nn
\\
\le&
c
\big(
\|
u_0
\|_{H^{s}(\rr^+)}
+
\|\varphi_1\|_{H_t^\frac{s}{3}(0,T)}
+
\frac12
\big\|
\psi_{2{T^*}}[\p_x(\widetilde{v}^2)]
\big\|_{X^{s,-b^*,\theta-1}_1}
+
\frac12
\big\|
\psi_{2{T^*}}[\p_x(\widetilde{v}^2)]
\big\|_{Y^{s,-b^*}_1}
\big),
\end{align}
where $\widetilde{v}$ is the extension of  $v$  from $\rr^+\times(0,T)$ to $\rr^2$ such that
$
\|\widetilde v\|_{X^{s,b,\theta}_\alpha}
\le
2
\|v\|_{X^{s,b,\theta}_{\alpha,\rr^+\times(0,T)}}.
$
For $\big\|
\psi_{2{T^*}}[\p_x(\widetilde{v}^2)]
\big\|_{Y^{s,-b^*}_1}$, using bilinear estimate \eqref{bi-est-Y-1} we get 
$
\big\|
\psi_{2{T^*}}[\p_x(\widetilde{v}^2)]
\big\|_{Y^{s,-b^*}_1}
\le
$
$
c_2
\big\|
\psi_{2{T^*}}\p_x(\widetilde{v}^2)
\big\|_{X^{s,-b^*}_1}
$
$
+
c_2
\|
\psi_{2{T^*}}\widetilde{v}
\|_{X^{s,b'}_\alpha}
\|
\widetilde{v}
\|_{X^{s,b'}_\alpha},
$
which combined with estimate \eqref{Rob2-onto-map-fin-est-large} implies that 
\begin{align}
\label{Rob2-onto-map-fin-est-no1}
\mathcal{N}_{\Phi}
\le
c
\big(
\|
u_0
\|_{H^{s}(\rr^+)}
+
\|\varphi_1\|_{H_t^\frac{s}{3}(0,T)}
+
\frac{1+c_2}2
\big\|
\psi_{2{T^*}}[\p_x(\widetilde{v}^2)]
\big\|_{X^{s,-b^*,\theta-1}_1}
+
\frac{c_2}2
\|
\psi_{2{T^*}}\widetilde{v}
\|_{X^{s,b'}_\alpha}
\|
\widetilde{v}
\|_{X^{s,b'}_\alpha}
\big).
\end{align}
Now we need the following result about multipliers.
\begin{lemma}
\label{tao-lemma-modified}
Let $\eta(t)$ be a function in the Schwartz space $\mathcal{S}(\rr)$ and $\alpha>0$. If 
$
-\frac12
<
b'
\le 
b
<
\frac12
$
and
$-\frac12<\theta'-1\le \theta-1<\frac12$ ($\frac12<\theta'\le \theta<1$ is sufficient condition)
then for any $0<{T^*}\le 1$ we have
\begin{align}
\label{tao-est}
&\|\eta(t/{T^*})u\|_{X^{s,b'}_{\alpha}}
\le
c_1(\eta,b',b,\alpha)
\,\,
{T^*}^{b-b'}\|u\|_{X^{s,b}_{\alpha}},
\\
\label{tao-est-modified}
&\|\eta(t/{T^*})u\|_{X^{s,b',\theta'-1}_\alpha}
\le
c_1(\eta,b,b',\theta,\theta',\alpha)
\,\,
\max\{
{T^*}^{b-b'},
{T^*}^{\theta-\theta'}
\}
\|u\|_{X^{s,b,\theta-1}_\alpha}.
\end{align}
\end{lemma}
\nin
For $\alpha=1$  the proof of estimate \eqref{tao-est} can be found in  \cite{tao-book}
 (see Lemma 2.11).
For  $\alpha\neq 1$, the estimate is similar. Also, the proof of estimate \eqref{tao-est-modified} is similar to the proof of estimate \eqref{tao-est}. 
Applying   estimate \eqref{tao-est-modified} with 
the  choices:
\begin{align}
\label{Rob2-b-b'-choice-wp-2}
&b
=
\frac12-\frac12\beta
\text{  ($-b$ is in place of $b$) }
\quad
\text{and}
\quad
b^*
=
\frac12-\frac14\beta
\text{  ($-b^*$ is in place of $b'$)},
\\
\label{Rob2-a-a1-choose}
&\theta
=
\frac12+\frac14\beta
\text{  (in place of $\theta'$) }
\quad
\text{and}
\quad
\theta^*
=
\frac12+\frac12\beta
\text{  (in place of $\theta$)},
\end{align}
where $\beta$ is defined in \eqref{beta-choice}, we get 
$
\big\|
\psi_{2{T^*}}[\p_x(\widetilde{v}^2)]
\big\|_{X^{s,-b^*,\theta-1}_1}
\le
c_1
{T^*}^{\beta/4}
\big\|
\p_x(\widetilde{v}^2)
\big\|_{X^{s,-b,\theta^*-1}_1}.
$
Also, using estimate \eqref{tao-est} with $b'=\frac12-\frac34\beta$ we get 
$\|
\psi_{2{T^*}}\widetilde{v}
\|_{X^{s,b'}_\alpha}
\le
c_1
{T^*}^{\beta/4}
\|
\widetilde{v}
\|_{X^{s,b}_\alpha}
$.
Combining these   with  \eqref{Rob2-onto-map-fin-est-no1}
gives
\begin{align}
\label{Rob2-onto-map-fin-est-no2}
\mathcal{N}_{\Phi}
\le
c
\big(
\|
u_0
\|_{H^{s}(\rr^+)}
+
\|\varphi_1\|_{H_t^\frac{s}{3}(0,T)}
+
\frac{c_1(1+c_2)}2
{T^*}^{\beta/4}
\big\|
\p_x(\widetilde{v}^2)
\big\|_{X^{s,-b,\theta^*-1}_1}
+
\frac{c_1c_2}2
{T^*}^{\beta/4}
\|
\widetilde{v}
\|_{X^{s,b}_\alpha}^2
\big).
\end{align}
Furthermore, applying bilinear estimate \eqref{bi-est-X-1} with  $b=b'=\frac12-\frac12\beta$ and $f=g=\widetilde{v}$
 for $\big\|
\p_x(\widetilde{v}^2)
\big\|_{X^{s,-b,\theta^*-1}_1}$, from \eqref{Rob2-onto-map-fin-est-no2} we get
$
\mathcal{N}_{\Phi}
\le
c_3
\big(
\|
u_0
\|_{H^{s}(\rr^+)}
+
\|\varphi_1\|_{H_t^\frac{s}{3}(0,T)}
+
{T^*}^{\frac14\beta}
\|
\widetilde v
\|_{X^{s,b,\theta}_{\alpha}}^2
\big)
$ 
or,
\begin{align}
\label{Rob2-onto-map-fin-est-1}
\|\Phi_ {T^*}(u,v)\|_{X^{s,b,\theta}_{1,\rr^+\times(0,T)}}
\le
c_3
\big(
\|
u_0
\|_{H^{s}(\rr^+)}
+
\|\varphi_1\|_{H_t^\frac{s}{3}(0,T)}
+
4
{T^*}^{\frac14\beta}
\|
v
\|_{X^{s,b,\theta}_{\alpha,\rr^+\times(0,T)}}^2 \big),
\end{align}
where $c_3\doteq c+cc_1(1+c_2)c_2$.
Concerning  $\Psi_{T^*}$, working similarly we obtain the estimate
\begin{align}
\label{Rob2-onto-map-fin-est-1-psi}
\hskip-0.1in
\|\Psi_ {T^*}(u,v)\|_{X^{s,b,\theta}_{\alpha,\rr^+\times(0,T)}}
\le
c_3
\big(
\|
v_0
\|_{H^{s}(\rr^+)}
+
\|\varphi_2\|_{H_t^\frac{s}{3}(0,T)}
\hskip-0.05in +
4
{T^*}^{\frac14\beta}
\|
u
\|_{X^{s,b,\theta}_{1,\rr^+\times(0,T)}}
\|
v
\|_{X^{s,b,\theta}_{\alpha,\rr^+\times(0,T)}}
\big).
\end{align}
Now, for $(u,v)\in \mathcal{X}_{MN}$, by estimates \eqref{Rob2-onto-map-fin-est-1} and  \eqref{Rob2-onto-map-fin-est-1-psi},
we see that
$\Phi\times \Psi$  is onto  $\mathcal{X}_{MN}$
if  the following two conditions are satisfied
\begin{align}
\label{Rob2-onto-map-condition-nls}
&c_3
\big(
\|
u_0
\|_{H^{s}(\rr^+)}
+
\|\varphi_1\|_{H_t^\frac{s}{3}(0,T)}
\big)
+
4
c_3
{T^*}^{\frac14\beta}
N^2
\le 
M,
\\
\label{Rob2-onto-map-condition-kdv}
&c_3
\big(
\|
v_0
\|_{H^{s}(\rr^+)}
+
\|\varphi_2\|_{H_t^\frac{s}{3}(0,T)}
\big)
+
4
c_3
{T^*}^{\frac14\beta}
MN
\le 
N.
\end{align}
{\bf $\Phi_{T^*}\times\Psi_{T^*}$  is contraction.}
Again,
using  the linear 
estimate \eqref{MB-forced-linear-est1-Rob} with forcing replaced by $\psi_{2{T^*}}[-\frac12\p_x(v_1^2
-
v_2^2)]$,  
for $0<b\le b^*<\frac12$  we have
\begin{align}
\label{Rob2-contraction-1}
\widetilde{\mathcal{N}}_{\Phi}
\doteq
&\|\Phi_{T^*}(u_1,v_1)-\Phi_{T^*}(u_2,v_2)\|_{X^{s,b,\theta}_{1,\rr^+\times(0,T)}}
\le
\|\Phi_{T^*}(u_1,v_1)-\Phi_{T^*}(u_2,v_2)\|_{X^{s,b^*,\theta}_{1,\rr^+\times(0,T)}}
\nn
\\
\le&
\frac12
c
\big\|
\psi_{2{T^*}}[\p_x(v_1^2
-
v_2^2)]
\big\|_{X^{s,-b^*,\theta-1}_{1,\rr^+\times(0,T)}}
+
\frac12
c
\big\|
\psi_{2{T^*}}[\p_x(v_1^2
-
v_2^2)]
\big\|_{Y^{s,-b^*}_{1,\rr^+\times(0,T)}}
\nn
\\
\le&
\frac12
c
\big\|
\psi_{2{T^*}}
[\p_x(\widetilde{v}_1^2
-
\widetilde{v}_2^2)]
\big\|_{X^{s,-b^*,\theta-1}_1}
+
\frac12
c
\big\|
\psi_{2{T^*}}[\p_x(\widetilde{v}_1^2
-
\widetilde{v}_2^2)]
\big\|_{Y^{s,-b^*}_{1}},
\end{align}
where   $\widetilde{v}_1$,    $\widetilde{v}_2$ are  extensions of  $v_1$,  $v_2$ from $\rr^+\times(0,T)$ to $\rr^2$ respectively, 
with
$
\|\widetilde{v}_2\|_{X^{s,b,\theta}_\alpha}
\le
2
\|v_2\|_{X^{s,b,\theta}_{\alpha,\rr^+\times(0,T)}}.
$
The extension $\widetilde{v}_1$ is obtained as follows. 
First, we extend
$v=v_1-v_2$ from $\rr^+\times(0,T)$ to $\rr^2$ such that 
$
\|\widetilde v_1-\widetilde v_2\|_{X^{s,b,\theta}_\alpha}
=
\|\widetilde v\|_{X^{s,b,\theta}_\alpha}
\le
2
\|v_1-v_2\|_{X^{s,b,\theta}_{\alpha,\rr^+\times(0,T)}}.
$
Then defining $\widetilde{v}_1\doteq \widetilde v+\widetilde{v}_2$, we see that 
$\widetilde{v}_1$ extends $v_1$ from $\rr^+\times(0,T)$ to $\rr^2$
and it satisfies
$
\|\widetilde{v}_1\|_{X^{s,b,\theta}_\alpha}
\le
\|\widetilde{v}\|_{X^{s,b,\theta}_\alpha}
+
\|\widetilde{v}_2\|_{X^{s,b,\theta}_\alpha}
\le
2
\|v_1-v_2\|_{X^{s,b,\theta}_{\alpha,\rr^+\times(0,T)}}
+
2
\|v_2\|_{X^{s,b,\theta}_{\alpha,\rr^+\times(0,T)}}
\le
6N.
$
Using basic identities 
about the nonlinearities of our MB system, like
\begin{align}
\label{product-algebral}
v_1^2
-
v_2^2
=
(v_1-v_2)(v_1+v_2)
\quad
\text{and}
\quad
u_1v_1
-
u_2v_2
=
(u_1-u_2)v_1
+
(v_1-v_2)u_2,
\end{align}
from the  bilinear estimates \eqref{bi-est-X-1}-- \eqref{bi-est-Y-2}
we get the following useful estimates.
\begin{lemma}
\label{mb-bi-lem-diff}
If $s_c(\alpha)\le s<3$, $\frac12-\beta\le b'\le b<\frac12<\theta'\le \theta\le \frac12+\beta$, then we have 
\begin{align}
&\|
\p_x(v_1^2-v_2^2)
\|_{X^{s,-b,\theta-1}_1}
\le
c_2
(
\|
v_1
\|_{X^{s,b',\theta'}_\alpha}
+
\|
v_2
\|_{X^{s,b',\theta'}_\alpha}
)
\|
v_1-v_2
\|_{X^{s,b',\theta'}_\alpha},
\nn
\\
&\|
\p_x(u_1v_1
-
u_2v_2)
\|_{X^{s,-b,\theta-1}_\alpha}
\le
c_2
\|
v_1
\|_{X^{s,b',\theta'}_\alpha}
\|
u_1-u_2
\|_{X^{s,b',\theta'}_1}
+
c_2
\|
u_2
\|_{X^{s,b',\theta'}_1}
\|
v_1-v_2
\|_{X^{s,b',\theta'}_\alpha},
\nn
\\
\label{mb-bilinear-est-diff}
&\|
\p_x(v_1^2-v_2^2)
\|_{Y^{s,-b}_1}
\le
c_2
\|
\p_x(v_1^2-v_2^2)
\|_{X^{s,-b}_1}
+
c_2
(
\|
v_1
\|_{X^{s,b'}_\alpha}
+
\|
v_2
\|_{X^{s,b'}_\alpha}
)
\|
v_1-v_2
\|_{X^{s,b'}_\alpha},
\\
%
%
&\|
\p_x(u_1v_1
\hskip-0.03in - \hskip-0.03in
u_2v_2)
\|_{Y^{s,-b}_\alpha}
\le
c_2
\big[
\|
\p_x(u_1v_1
\hskip-0.03in - \hskip-0.03in
u_2v_2)
\|_{X^{s,-b}_\alpha}
\hskip-0.03in + \hskip-0.03in
\|
v_1
\|_{X^{s,b'}_\alpha}
\|
u_1
\hskip-0.03in - \hskip-0.03in
u_2
\|_{X^{s,b'}_1}
\hskip-0.03in + \hskip-0.03in
\|
u_2
\|_{X^{s,b'}_1}
\|
v_1
\hskip-0.03in - \hskip-0.03in
v_2
\|_{X^{s,b'}_\alpha}
\big].
\nn
\end{align}
\end{lemma}

Next,   applying bilinear estimates \eqref{mb-bilinear-est-diff} we get 
$$
\big\|
\psi_{2{T^*}}[\p_x(\widetilde{v}_1^2
-
\widetilde{v}_2^2)]
\big\|_{Y^{s,-b^*}_{1}}
\le
c_2
\|
\psi_{2{T^*}}\p_x(\widetilde{v}_1^2-\widetilde{v}_2^2)
\|_{X^{s,-b^*}_1}
+
c_2
(
\|
\widetilde{v}_1
\|_{X^{s,b'}_\alpha}
+
\|
\widetilde{v}_2
\|_{X^{s,b'}_\alpha}
)
\|
\psi_{2{T^*}}(
\widetilde{v}_1-\widetilde{v}_2
)
\|_{X^{s,b'}_\alpha},
$$
which combined with estimate \eqref{Rob2-contraction-1} gives us that 
\begin{align*}
\widetilde{\mathcal{N}}_{\Phi}
\le
\frac{c(1+c_2)}2
\big\|
\psi_{2{T^*}}
[\p_x(\widetilde{v}_1^2
-
\widetilde{v}_2^2)]
\big\|_{X^{s,-b^*,\theta-1}_1}
+
\frac{cc_2}2
(
\|
\widetilde{v}_1
\|_{X^{s,b'}_\alpha}
+
\|
\widetilde{v}_2
\|_{X^{s,b'}_\alpha}
)
\|
\psi_{2{T^*}}(
\widetilde{v}_1-\widetilde{v}_2
)
\|_{X^{s,b'}_\alpha}.
\end{align*}
Applying 
the multiplier estimate \eqref{tao-est-modified} with $b$, $b^*$ chosen in \eqref{Rob2-b-b'-choice-wp-2} and $\theta,\theta^*$ chosen in \eqref{Rob2-a-a1-choose} for $\big\|
\psi_{2{T^*}}
[\p_x(\widetilde{v}_1^2
-
\widetilde{v}_2^2)]
\big\|_{X^{s,-b^*,\theta-1}_1}$ and applying estimate \eqref{tao-est} with $b'=\frac12-\frac34\beta$ for $\|
\psi_{2{T^*}}(
\widetilde{v}_1-\widetilde{v}_2
)
\|_{X^{s,b'}_\alpha}$ we get
\begin{align}
\label{Rob2-contraction-map-fin-est}
\widetilde{\mathcal{N}}_{\Phi}
\le
\frac{cc_1(1+c_2)}2
{T^*}^{\frac14\beta}
\big\|
\p_x(\widetilde{v}_1^2
-
\widetilde{v}_2^2)
\big\|_{X^{s,-b,\theta^*-1}_1}
+
\frac{cc_1c_2}2
{T^*}^{\frac14\beta}
(
\|
\widetilde{v}_1
\|_{X^{s,b'}_\alpha}
+
\|
\widetilde{v}_2
\|_{X^{s,b'}_\alpha}
)
\|
\widetilde{v}_1-\widetilde{v}_2
\|_{X^{s,b}_\alpha}
\end{align}
Also, using  estimates \eqref{mb-bilinear-est-diff},  from \eqref{Rob2-contraction-map-fin-est} we get
$\widetilde{\mathcal{N}}_{\Phi}
\le
c_3
{T^*}^{\frac14\beta}
(
\|
\widetilde{v}_1
\|_{X^{s,b,\theta}_\alpha}
+
\|
\widetilde v_2
\|_{X^{s,b,\theta}_\alpha}
)
\|
\widetilde v_1-\widetilde v_2
\|_{X^{s,b,\theta}_\alpha}
$,
or
\begin{align}
\label{Rob2-contraction-map-fin-est-1}
\|\Phi_{T^*}(u_1,v_1)-\Phi_{T^*}(u_2,v_2)\|_{X^{s,b,\theta}_{1,\rr^+\times(0,T)}}
\le
2c_3
{T^*}^{\frac14\beta}
8N
\|
v_1-v_2
\|_{X^{s,b,\theta}_{\alpha,\rr^+\times(0,T)}}.
\end{align}
For   $\Psi_{T^*}$, working similarly we obtain the companion estimate
\begin{align}
\label{Rob2-contraction-map-fin-est-1-psi}
\hskip-0.1in
\|\Psi_{T^*}(u_1,v_1)
\hskip-0.02in
-
\hskip-0.02in
\Phi_{T^*}(u_2,v_2)\|_{X^{s,b,\theta}_{\alpha,\rr^+\times(0,T)}}
\hskip-0.1in
\le
2c_3
{T^*}^{\frac14\beta}
(6N
\|
u_1
\hskip-0.02in
-
\hskip-0.02in
u_2
\|_{X^{s,b,\theta}_{1,\rr^+\times(0,T)}}
\hskip-0.1in
+
2M
\|
v_1
\hskip-0.02in
-
\hskip-0.02in
v_2
\|_{X^{s,b,\theta}_{\alpha,\rr^+\times(0,T)}}
\hskip-0.03in
).
\end{align}
Combining estimates
 \eqref{Rob2-onto-map-condition-nls}, 
 \eqref{Rob2-onto-map-condition-kdv},
 \eqref{Rob2-contraction-map-fin-est-1},
  \eqref{Rob2-contraction-map-fin-est-1-psi}, 
  and choosing $M$ and $N$ as follows 
$$
M
=
N
=
2c_3\big(
\|
u_0
\|_{H^{s}(\rr^+)}
+
\|
v_0
\|_{H^{s}(\rr^+)}
+
\|\varphi_1\|_{H_t^\frac{s}{3}(0,T)}
+
\|\varphi_2\|_{H_t^\frac{s}{3}(0,T)}
\big),
$$ 
we see that the iteration map 
$\Phi_{T^*}\times\Psi_{T^*}$ is a contraction,
if  $T^*$ satisfies the conditions:
$
4c_3
{T^*}^{\frac14\beta}
N^2
\le
\frac{M}{2},
$
$
4c_3
{T^*}^{\frac14\beta}
MN
\le
\frac{N}{2},
$
$
16c_3N
{T^*}^{\frac14\beta}
\le
\frac14
$
and
$
2c_3(
6N
+
2M
)
{T^*}^{\frac14\beta}
\le
\frac14.
$
Solving these inequalities, we get
$$
{T^*}^{\frac14\beta}
\le
\frac{1}{64c_3 M}
\Longrightarrow
T^*
\lesssim
\big[
2c_3\big(
\|
u_0
\|_{H^{s}(\rr^+)}
+
\|
v_0
\|_{H^{s}(\rr^+)}
+
\|\varphi_1\|_{H_t^\frac{s}{3}(0,T)}
+
\|\varphi_2\|_{H_t^\frac{s}{3}(0,T)}
\big)
+1
\big]^{-4/\beta},
$$
which holds if we choose  
 a lifespan $T_0<\frac12$ satisfying 
estimate \eqref{MB-lifespan-Rob1}
stated in  our main Theorem  \ref{thm-MB-half-line-Rob}.
The proof for the Lip-continuity of the data to solution map and 
for the uniqueness of solution is similar to the one presented in \cite{bop1998} for  well-posedness on the line.
\,\, $\square$

%
%
%
%
%
%
%
%
%
%
%
%
\section{Proof of Bilinear Estimates in $X^{s,b,\theta}$}
\label{sec:bilinear-est}
\setcounter{equation}{0}
In this section we prove the  spatial bilinear estimates 
in modified Bourgain spaces. We focus mostly 
on estimate  \eqref{bi-est-X-1} for 
the coupled nonlinearity of the 
 $u$-equation.
Since the bilinear estimates for the $v$-equation are similar,
we only provide a brief outline of their main 
differences at the end of this section.
 Also, the bilinear estimates for $\alpha=1$ (i.e. for the KdV equation) 
 can be found  in \cite{h2006}.
 
We begin by recalling the following 
calculus  estimates \cite{kpv1996, h2006}.

\begin{lemma}
\label{lem:calc_ineq}  
If $1>\ell>1/2, 1>l'>1/2$ then
\begin{subequations}
\begin{align}
\label{eq:calc_1a}
\int_{\mathbb{R}}\frac{dx}{(1+|x-a|)^{2\ell}(1+|x-c|)^{2\ell'}}&\lesssim\frac{1}{(1+|a-c|)^{2\min\{\ell',\ell\}}},
\\
\label{eq:calc_4}
\int_{|x|\leq c}\frac{dx}{(1+|x|)^{2(1-l)}\sqrt{|a-x|}}&\lesssim\frac{(1+c)^{2(\ell-1/2)}}{(1+|a|)^{1/2}}.
\end{align}
\end{subequations}
In addition, if $\frac14<\ell'\leq \ell<\frac12$, then
\begin{equation}
\label{eq:calc_5}
\int_{\mathbb{R}}\frac{dx}{(1+|x-a|)^{2\ell}(1+|x-c|)^{2\ell'}}\lesssim\frac{1}{(1+|a-c|)^{2\ell+2\ell'-1}}.
\end{equation}

\end{lemma}

\noindent
Next, we express our bilinear estimate \eqref{bi-est-X-1},
i.e.
$
\|
\p_x(fg)
\|_{X^{s,-b,\theta-1}_{1}}
\le
c_2
\|
f
\|_{X^{s,b',\theta'}_{\alpha}}
\|
g
\|_{X^{s,b',\theta'}_{\alpha}}
$
in its $L^2$ form. Using the fact that $a^2+b^2\leq(|a|+|b|)^2\leq 2(a^2+b^2)$, we get
$
\|w\|_{X^{s,b,\theta}_{\alpha}}^2
\simeq
\int_{\rr^2}
\big[
(1+|\xi|)^{s}
(1+|\tau-\alpha\xi^3|)^{b}
+
$
$
\chi_{|\xi|< 1}(1+|\tau|)^{\theta}
\big]^2
|\widehat{w}(\xi,\tau)|^2
d\xi
d\tau.
$
Following Colliander and Kenig's notation \cite{ck2002}, 
for a function $h$ 
we denote the integrand under the integral of the modified Bourgain norm by
\begin{equation}
\label{eq:c_u}
c_{h,\alpha}(\xi,\tau)
\doteq
\Big[
(1+|\xi|)^{s}
(1+|\tau-\alpha\xi^3|)^{b'}
+
\chi_{|\xi|< 1}(1+|\tau|)^{\theta'}
\Big]
|\widehat{h}(\xi,\tau)|.
\end{equation}
Then,  the modified Bourgain norm of $h$ is 
approximately  equal to
the $L^2$ norm of $c_{h,\alpha}$, that is
\begin{equation}
\label{eq:c_useful}
\|h\|_{X^{s,b',\theta'}_{\alpha}}
\simeq
\|c_{h,\alpha}(\xi,\tau)\|_{L^2_{\xi,\tau}}.
\end{equation}
Thus, to prove bilinear estimate  \eqref{bi-est-X-1},
 it suffices to show the following 
$L^2$ inequality (formulation)
\begin{align}
\label{bilinear-est-L2-form}
 \Big\| 
\int_{\rr^2}
 Q(\xi,\tau,\xi_1,\tau_1)
 c_{f,\alpha}(\xi-\xi_1,\tau-\tau_1) c_{g,\alpha}(\xi_1,\tau_1) 
 d\xi_1 d\tau_1 \Big\|_{L^2_{\xi,\tau}}
\lesssim
\big\|c_{f,\alpha}\big\|_{L^2_{\xi,\tau}} \big\|c_{g,\alpha}\big\|_{L^2_{\xi,\tau}},
\end{align}
where $Q$ is the  multiplier defined by
\begin{align*}
Q(\xi,\tau,\xi_1,\tau_1)
\doteq&
|\xi|
\cdot
\Big[
\frac{
(1+|\xi|)^s
}{
(1+|\tau-\xi^{3}|)^{b}
}
+
\frac{
\chi_{|\xi|< 1}
}
{
(1+|\tau|)^{1-\theta}
}
\Big]
\cdot
\frac{1}{
(1+|\xi_1|)^{s}
(1+|\tau_1-\alpha\xi_1^3|)^{b'}
+
\chi_{|\xi_1|< 1}(1+|\tau_1|)^{\theta'}
}
\\
\times&
\frac{1}{
(1+|\xi-\xi_1|)^{s}
(1+|\tau-\tau_1-\alpha(\xi-\xi_1)^3|)^{b'}
+
\chi_{|\xi-\xi_1|< 1}(1+|\tau-\tau_1|)^{\theta'}
}.
\end{align*}
For $b'\le b$ and $
1-\theta
\ge
b,
$
we have 
$
\big[
\frac{
(1+|\xi|)^s
}{
(1+|\tau-\xi^{3}|)^{b}
}
+
\frac{
\chi_{|\xi|< 1}
}
{
(1+|\tau|)^{1-\theta}
}
\big]
\lesssim
\frac{
(1+|\xi|)^s
}{
(1+|\tau-\xi^{3}|)^{b'}
}.
$
Using this
and
factoring out 
$
\frac{(1+|\xi|)^s}{
(1+|\xi_1|)^{s}
(1+|\xi-\xi_1|)^{s}
},
$
we get
$
Q(\xi,\tau,\xi_1,\tau_1)
\lesssim
Q'(\xi,\tau,\xi_1,\tau_1),
$
where 
\begin{align}
\label{Q'-def}
Q'(\xi,\tau,\xi_1,\tau_1)
=&
\frac{|\xi|}{(1+|\tau-\xi^{3}|)^{b'}}
\frac{(1+|\xi|)^s}{
(1+|\xi_1|)^{s}
(1+|\xi-\xi_1|)^{s}
}
\frac{1}{
(1+|\tau_1-\alpha\xi_1^3|)^{b'}
+
\chi_{|\xi_1|< 1}(1+|\tau_1|)^{\theta'}
}
\nn
\\
\times&
\frac{1}
{
(1+|\tau-\tau_1-\alpha(\xi-\xi_1)^3|)^{b'}
+
\chi_{|\xi-\xi_1|< 1}(1+|\tau-\tau_1|)^{\theta'}
}.
\end{align}
Hence for any $s\in\rr$,
to  prove the bilinear estimate \eqref{bilinear-est-L2-form-no1}, it suffices to prove that 
\begin{align}
\label{bilinear-est-L2-form-no1}
 \Big\| 
\int_{\rr^2}
 Q'(\xi,\tau,\xi_1,\tau_1)
c_{f,\alpha}(\xi-\xi_1,\tau-\tau_1) c_{g,\alpha}(\xi_1,\tau_1) 
d\xi_1 d\tau_1 \Big\|_{L^2_{\xi,\tau}}
\lesssim
\big\|c_{f,\alpha}\big\|_{L^2_{\xi,\tau}} \big\|c_{g,\alpha}\big\|_{L^2_{\xi,\tau}}.
\end{align}
%
%
%
%
%
%
%
%
%
%
%
%
\subsection{The case $s\ge 0$}
Then,
we have
$
\frac{(1+|\xi|)^s}{(1+|\xi-\xi_1|)^s(1+|\xi_1|)^s}
\lesssim
 1,
$
which implies that $Q'\leq Q_0$, where
\begin{align*}
Q_0(\xi,\tau,\xi_1,\tau_1)
\doteq&
\frac{|\xi|}{(1+|\tau-\xi^{3}|)^{b'}}
\frac{1}{
(1+|\tau_1-\alpha\xi_1^3|)^{b'}
+
\chi_{|\xi_1|< 1}(1+|\tau_1|)^{\theta'}
}
\\
\times&
\frac{1}
{
(1+|\tau-\tau_1-\alpha(\xi-\xi_1)^3|)^{b'}
+
\chi_{|\xi-\xi_1|< 1}(1+|\tau-\tau_1|)^{\theta'}
}.
\end{align*}
Therefore, when $s\geq 0$ the proof of (\ref{bilinear-est-L2-form-no1}) is 
reduced to proving
\begin{align}
\label{bilinear-est-s-0}
\Big\| 
\int_{\rr^2}
Q_0(\xi,\tau,\xi_1,\tau_1)
c_{f,\alpha}(\xi-\xi_1,\tau-\tau_1) c_{g,\alpha}(\xi_1,\tau_1) 
d\xi_1 d\tau_1 \Big\|_{L^2_{\xi,\tau}}
\lesssim
\big\|c_{f,\alpha}\big\|_{L^2_{\xi,\tau}} \big\|c_{g,\alpha}\big\|_{L^2_{\xi,\tau}},
\end{align}
which corresponds to proving the bilinear estimate when $s=0$.
By symmetry (in convolution writing), we may assume that 
\begin{equation}
\label{conv-symmetry-xi}
|\xi-\xi_1|
\leq
|\xi_1|.
\end{equation}
Next, for simplifying our proof,
we decompose the multiplier
$Q_0(\xi,\tau,\xi_1,\tau_1)$ as the sum of three multipliers
$
Q_0(\xi,\tau,\xi_1,\tau_1)
=
Q_1(\xi,\tau,\xi_1,\tau_1)
+
Q_2(\xi,\tau,\xi_1,\tau_1)
+
Q_3(\xi,\tau,\xi_1,\tau_1),
$
where
\begin{align}
\label{u-s-0-Q1}
&Q_1(\xi,\tau,\xi_1,\tau_1)
\doteq
\chi_{|\xi_1|< 1}
\cdot
\chi_{|\xi-\xi_1|< 1}
\cdot
Q_0(\xi,\tau,\xi_1,\tau_1)
\\
=&
\frac{|\xi|}{(1+|\tau-\xi^{3}|)^{b'}}
\frac{
\chi_{|\xi_1|< 1}
}{
(1+|\tau_1-\alpha\xi_1^3|)^{b'}
+
(1+|\tau_1|)^{\theta'}
}
\frac{
\chi_{|\xi-\xi_1|< 1}
}
{
(1+|\tau-\tau_1-\alpha(\xi-\xi_1)^3|)^{b'}
+
(1+|\tau-\tau_1|)^{\theta'}
},
\nn
\end{align}
\begin{align}
\label{u-s-0-Q2}
&Q_2(\xi,\tau,\xi_1,\tau_1)
\doteq
\chi_{|\xi_1|\ge 1}
\cdot
\chi_{|\xi-\xi_1|<1}
\cdot
Q_0(\xi,\tau,\xi_1,\tau_1)
\nn
\\
=&
\frac{|\xi|}{(1+|\tau-\xi^{3}|)^{b'}}
\frac{
\chi_{|\xi_1|\ge 1}
}{
(1+|\tau_1-\alpha\xi_1^3|)^{b'}
}
\frac{
\chi_{|\xi-\xi_1|< 1}
}
{
(1+|\tau-\tau_1-\alpha(\xi-\xi_1)^3|)^{b'}
+
(1+|\tau-\tau_1|)^{\theta'}
},
\end{align}
\vskip-0.1in
\nin
and
\begin{align}
\label{u-s-0-Q3}
Q_3(\xi,\tau,\xi_1,\tau_1)
\doteq&
\chi_{|\xi_1|\ge 1}
\cdot
\chi_{|\xi-\xi_1|\ge 1}
\cdot
Q_0(\xi,\tau,\xi_1,\tau_1)
\nn
\\
=&
\frac{|\xi|}{(1+|\tau-\xi^{3}|)^{b'}}
\frac{
\chi_{|\xi_1|\ge 1}
}{
(1+|\tau_1-\alpha\xi_1^3|)^{b'}
}
\frac{
\chi_{|\xi-\xi_1|\ge 1}
}
{
(1+|\tau-\tau_1-\alpha(\xi-\xi_1)^3|)^{b'}
}.
\end{align}
Now, in the multiplier $Q_3$ we recognize the 
familiar (from the Cauchy problem theory)
 Bourgain quantity 
\begin{align}
\label{d3-quantity}
(\tau-\xi^3)-(\tau_1-\alpha\xi_1^3)-[\tau-\tau_1-\alpha(\xi-\xi_1)^3]
=
-\xi^3+\alpha\xi_1^3+\alpha(\xi-\xi_1)^3
\doteq
d_{\alpha}(\xi,\xi_1),
\end{align}
which is related to the $u$-equation.
Observe that for the KdV equation ($\alpha=1$)
this quantity is 
\begin{equation}
\label{d-1}
d_1(\xi,\xi_1)
\doteq
-\xi^3+\xi_1^3+(\xi-\xi_1)^3
=
3\xi\xi_1(\xi_1-\xi)
.
\end{equation}
Below, we list several useful and elementary properties 
for this quantity.
%
%
%
%
%
%
%
%
%
%
%
\begin{lemma}
\label{d-alpha-property}
The  Bourgain quantity  $d_\alpha (\xi,\xi_1)$
satisfies the following  properties:
\begin{equation}
\label{d-comparison}
\max\{
|\tau-\xi^3|, |\tau_1-\alpha\xi_1^3|, |\tau-\tau_1-\alpha(\xi-\xi_1)^3|
\}
\ge
\frac13
|d_{\alpha}(\xi,\xi_1)|,
\end{equation}
\vskip-0.3in
\begin{align}
\label{d-le-4-xi1}
d_{\alpha}(\xi,\xi_1)
=&
3\alpha\xi
(\xi_1
-
r_1
\xi)
(
\xi_1-
r_2\xi
),
\,\,
\text{with}
\,\,
r_1
=
\frac12-\frac{\sqrt{-3+12\alpha^{-1}}}{6},
\,\,
r_2
=
\frac12+\frac{\sqrt{-3+12\alpha^{-1}}}{6},
\\
\label{d-4}
d_{4}(\xi,\xi_1)
=&
3\xi
(\xi-2\xi_1)^2
=
12
\xi
(\xi_1-\xi/2)^2,
\\
\label{d-bound-xi}
d_{\alpha}(\xi,\xi_1)
=&
3\alpha\xi
\big[
(\xi_1-\frac12\xi)^2
+
\frac{\alpha-4}{12\alpha}
\xi^2
\big]
\Longrightarrow
|d_{\alpha}(\xi,\xi_1)|
\ge
\frac{\alpha-4}{4}
|\xi|^3,
\,\,
\alpha>4,
\\
\label{d-bound-xi1}
d_{\alpha}(\xi,\xi_1)
=&
(\alpha-1)\xi
\Big[
\Big(
\xi-\frac{3\alpha\,\xi_1}{2(\alpha-1)}
\Big)^2
\hskip-0.04in
+
\frac{3\alpha(\alpha-4)}{4(\alpha-1)^2}
\xi_1^2
\Big]
\Longrightarrow
|d_{\alpha}(\xi,\xi_1)|
\ge
\frac{3\alpha(\alpha-4)}{4(\alpha-1)}
|\xi||\xi_1|^2,
\,\,
\alpha>4,
\end{align}
\vskip-0.2in
\begin{align}
\label{d-partial-xi1}
\frac{\p d_{\alpha}}{\p\xi_1}
=&
3\alpha\xi_1^2
-
3\alpha(\xi_1-\xi)^2
=
6\alpha
\xi(\xi_1-\xi/2),
\\
\label{d-partial-xi}
\frac{\p d_{\alpha}}{\p\xi}
=&
3(\alpha-1)(\xi-p_1)(\xi-p_2),
\,\,
\text{with}
\,\,
p_1
=
\frac{\sqrt{\alpha}}{\sqrt{\alpha}-1}\xi_1,
\,\,
p_2
=
\frac{\sqrt{\alpha}}{\sqrt{\alpha}+1}\xi_1,
\quad
\alpha\neq 1,
\\
\label{d-2partial-xi}
\frac{\p^2 d_{\alpha}}{\p\xi^2}
=&
-6\xi
+
6\alpha(\xi-\xi_1)
=
6(\alpha-1)
(\xi-q),
\,\,
\text{with}
\,\,
q
=
\frac{\alpha}{\alpha-1}\xi_1,
\quad
\alpha\neq 1.
\end{align}
\end{lemma}
\nin
To prove the bilinear estimate \eqref{bilinear-est-s-0},  it suffices to prove that for $\ell
=
1,2,3$ we have
\begin{align}
\label{bilinear-est-s-l}
\Big\| 
\int_{\rr^2}
Q_\ell(\xi,\tau,\xi_1,\tau_1)
c_{f,\alpha}(\xi-\xi_1,\tau-\tau_1) c_{g,\alpha}(\xi_1,\tau_1) 
d\xi_1 d\tau_1 \Big\|_{L^2_{\xi,\tau}}
\lesssim
\big\|c_{f,\alpha}\big\|_{L^2_{\xi,\tau}} \big\|c_{g,\alpha}\big\|_{L^2_{\xi,\tau}}.
\end{align}
\underline{\bf Proof for  multiplier $Q_1$.} 
Using the Cauchy--Schwarz inequality with respect to 
$(\xi_1,\tau_1)$  and taking the supremum in  $(\xi,\tau)$, 
we arrive at
\begin{align*}
\Big\| 
\hskip-0.05in
\int_{\rr^2}
\hskip-0.07in
Q_1(\xi,\tau,\xi_1,\tau_1)
 c_{f,\alpha}(\xi-\xi_1,\tau-\tau_1) c_{g,\alpha}(\xi_1,\tau_1) 
 d\xi_1 d\tau_1 \Big\|_{L^2}
 \hskip-0.06in
\lesssim
\hskip-0.02in
\big\|c_{f,\alpha}\big\|_{L^2} \big\|c_{g,\alpha}\big\|_{L^2}
\Big\|
\hskip-0.06in
\int_{\rr^2}
\hskip-0.06in
Q_1^2(\xi,\tau,\xi_1,\tau_1)
d\xi_1 d\tau_1
\Big\|_{L^{\infty}_{\xi,\tau}}^{\frac12}.
\end{align*}
Thus, to prove  bilinear estimate \eqref{bilinear-est-s-l} with $\ell=1$, it suffices to show 
 the following result for $Q_1$ in \eqref{u-s-0-Q1}.
\begin{lemma} 
\label{s-0-theta-1-lemma}
Let $\alpha>0$, $\theta'>\frac12$ and $b'\ge0$, for $\xi,\tau\in\rr$, we have
\begin{align}
\label{s-0-theta-1}
\Theta_1(\xi,\tau)
\doteq&
\frac{\xi^2}{(1+|\tau-\xi^{3}|)^{2b'}}
\int_{\rr^2}
\frac{
\chi_{|\xi_1|< 1}
}{
[
(1+|\tau_1-\alpha\xi_1^3|)^{b'}
+
(1+|\tau_1|)^{\theta'}
]^2
}
\nn
\\
\times&\frac{
\chi_{|\xi-\xi_1|< 1}
\quad
d\xi_1
d\tau_1
}
{
[
(1+|\tau-\tau_1-\alpha(\xi-\xi_1)^3|)^{b'}
+
(1+|\tau-\tau_1|)^{\theta'}
]^2
}
\lesssim 1.
\end{align}
\end{lemma}
\nin
{\bf Proof of Lemma \ref{s-0-theta-1-lemma}}. In \eqref{s-0-theta-1}, dropping  $(1+|\tau_1-\alpha\xi_1^3|)^{b'}$ and $(1+|\tau-\tau_1-\alpha(\xi-\xi_1)^3|)^{b'}$ from the denominator, 
applying  estimate \eqref{eq:calc_1a} with $x=\tau_1$,  $\ell=\ell'=\theta'>\frac12$, $a=0$ and $c=\tau$, we get
$
\Theta_1(\xi,\tau)
\leq
\frac{\xi^2}{(1+|\tau-\xi^{3}|)^{2b'}}
\int_\rr
\frac{
\chi_{|\xi_1|< 1}
\chi_{|\xi-\xi_1|<1}
}
{
(1+|\tau|)^{2\theta'}
}
d\xi_1\lesssim 1,
$
since $|\xi_1|<1$ and $|\xi-\xi_1|<1$.
\,\,
$\square$

\vskip0.05in
\noindent
\underline{\bf Proof for  multiplier $Q_2$.} We will  consider two possible  microlocalizations

\vskip0.05in
\noindent
{\bf $\bullet$ Microlocalization I: $|\tau-\xi^{3}|\ge |\tau_1-\alpha\xi_1^3|$}. In this case we define the domain 
\begin{align*}
B_{I}
\doteq
\big\{
(\xi,\tau,\xi_1,\tau_1) \in{\rr}^4: |\tau-\xi^{3}|\ge |\tau_1-\alpha\xi_1^3|,
|\xi_1|\ge 1,
|\xi-\xi_1|<1
\big\}.
\end{align*}
{\bf $\bullet$ Microlocalization II: $|\tau-\xi^{3}|< |\tau_1-\alpha\xi_1^3|$}. 
In this situation we define the domain 
\begin{align*}
B_{II}
\doteq
\big\{(\xi,\tau,\xi_1,\tau_1)\in{\rr}^4: |\tau-\xi^{3}|< |\tau_1-\alpha\xi_1^3|,
|\xi_1|\ge1,
|\xi-\xi_1|<1
\big\}.
\end{align*}
{\bf Estimate in microlocalization I:} $|\tau-\xi^{3}|\ge |\tau_1-\alpha\xi_1^3|$. As before, using the Cauchy--Schwarz inequality with respect to 
$(\xi_1,\tau_1)$, and taking the supremum over  $(\xi,\tau)$ we
arrive at
\begin{align*}
& \Big\| 
\int_{\rr^2}
 (\chi_{B_I} Q_2)(\xi,\tau,\xi_1,\tau_1)
 c_{f,\alpha}(\xi-\xi_1,\tau-\tau_1) c_{g,\alpha}(\xi_1,\tau_1) 
 d\xi_1 d\tau_1 \Big\|_{L^2_{\xi,\tau}}
 \nn
 \\
 &
\lesssim
\big\|c_{f,\alpha}\big\|_{L^2_{\xi,\tau}} \big\|c_{g,\alpha}\big\|_{L^2_{\xi,\tau}}
\Big\|
\int_{\rr^2}
 (\chi_{B_I} Q_2^2)(\xi,\tau,\xi_1,\tau_1)
d\xi_1 d\tau_1
\Big\|_{L^{\infty}_{\xi,\tau}}^{1/2}.
\end{align*}
Thus, to prove bilinear estimate \eqref{bilinear-est-s-l}, it suffices to show 
 the following result for $Q_2$ defined in \eqref{u-s-0-Q2}.
\begin{lemma} 
\label{s-0-theta-2-lemma}
Let $\alpha>0$.
If   $\frac13\leq b'<\frac12$ and 
$\theta'>1/2$, 
then for $\xi,\tau\in\rr$ we have
\begin{align}
\label{s-0-theta-2}
\hskip-0.12in
\Theta_2(\xi,\tau)
=&
\frac{
\xi^2
}{
(1\hskip-0.03in + \hskip-0.03in |\tau-\xi^{3}|)^{2b'}
}
\hskip-0.04in
\int_{\rr^2}
\hskip-0.04in
\frac{
\chi_{|\xi_1|\ge 1}
\,
\chi_{|\xi-\xi_1|<1}
\,
\chi_{B_{I}}(\xi,\tau,\xi_1,\tau_1)
\quad
d\xi_1
d\tau_1
}
{
(1\hskip-0.02in + \hskip-0.02in
|\tau_1-\alpha\xi_1^3|)^{2b'}
[
(1 \hskip-0.02in + \hskip-0.02in
|\tau-\tau_1-\alpha(\xi-\xi_1)^3|)^{b'}
\hskip-0.04in
+
(1  \hskip-0.02in + \hskip-0.02in
|\tau-\tau_1|)^{\theta'}
]^2
}
\lesssim 1.
\end{align}
\end{lemma}
\nin
{\bf Proof of Lemma \ref{s-0-theta-2-lemma}}.
Dropping $(1+|\tau-\tau_1-\alpha(\xi-\xi_1)^3|)^{b'}$ from the denominator, splitting 
$
\frac{1}{(1+|\tau-\xi^{3}|)^{2b'}}
=
\frac{1}{(1+|\tau-\xi^{3}|)^{6b'-2}}
\frac{1}{(1+|\tau-\xi^{3}|)^{2-4b'}}, 
$ 
putting $\frac{1}{(1+|\tau-\xi^{3}|)^{2-4b'}}$
inside the integral  in \eqref{s-0-theta-2} and using 
$
\frac{1}{(1+|\tau-\xi^{3}|)^{2-4b'}}
\lesssim
\frac{1}{(1+|\tau_1-\alpha\xi_1^3|)^{2-4b'}}
$ to replace $\frac{1}{(1+|\tau-\xi^{3}|)^{2-4b'}}$ by $\frac{1}{(1+|\tau_1-\alpha\xi_1^3|)^{2-4b'}}$,
we get
\begin{align}
\label{s-0-theta-2-est-no1}
\Theta_2(\xi,\tau)
\leq&
\frac{\xi^2}{(1+|\tau-\xi^{3}|)^{6b'-2}}
\int_{\rr^2}
\frac{
1
}{
(1+|\tau-\xi^3|)^{2-4b'}
}
\frac{
\chi_{|\xi_1|\ge1}
}{
(1+|\tau_1-\alpha\xi_1^{3}|)^{2b'}
}
\frac{
\chi_{|\xi-\xi_1|< 1}
}
{
(1+|\tau-\tau_1|)^{2\theta'}
}
d\xi_1
d\tau_1
\nn
\\
\leq&
\frac{\xi^2}{(1+|\tau-\xi^{3}|)^{6b'-2}}
\int_{\rr^2}
\frac{
\chi_{|\xi_1|\ge 1}
}{
(1+|\tau_1-\alpha\xi_1^{3}|)^{2-2b'}
}
\frac{
\chi_{|\xi-\xi_1|< 1}
}
{
(1+|\tau-\tau_1|)^{2\theta'}
}
d\xi_1
d\tau_1.
\end{align}
Now, applying  estimate \eqref{eq:calc_1a} with $x=\tau_1$, $\ell'=\theta'$, $\ell=1-b'$, $a=\tau$ and $c=\alpha\xi_1^3$, we get
\begin{equation}
\label{s-0-theta-2-est-no1}
\Theta_2(\xi,\tau)
\lesssim
\frac{
\xi^2
}{
(1+|\tau-\xi^3|)^{6b'-2}
}
\int_{\rr}
\frac{
\chi_{|\xi_1|\ge1}
\chi_{|\xi-\xi_1|< 1}
}{
(1+|\tau-\alpha\xi_1^{3}|)^{\min\{2\theta',2-2b'\}}
}
d\xi_1.
\end{equation}
Next we consider the following two possible cases

\vskip0.05in
\noindent
$\bullet $ Case 1: $|\xi|\leq 20$
\qquad
$\bullet $ Case 2: $|\xi|> 20$

\vskip0.05in
\noindent
{\bf Estimate in Case 1: $|\xi|\leq 20$.} Then, using the triangle inequality we get
$
|\xi_1|
=
|\xi_1-\xi+\xi|
\leq
|\xi-\xi_1|+|\xi|
\lesssim 1.
$
Combining this estimate with \eqref{s-0-theta-2-est-no1}, for $b'<1/2$ and $6b'-2\geq 0$ or $ 1/3\le b'<1/2$ we get
 the desired estimate \eqref{s-0-theta-2}.
This completes the proof in case 1.
\vskip0.05in
\noindent
{\bf Estimate in Case 2: $|\xi|> 20$.} Then, using the triangle inequality we have
$
|\xi_1|
\simeq
|\xi|.
$
Combining  this  with \eqref{s-0-theta-2-est-no1} and
making the change of variable, $\mu=\alpha\xi_1^3$, for $\min\{2\theta',2-b'\}>1$ we get
$
\Theta_2(\xi,\tau)
\lesssim
\frac{1}{(1+|\tau-\xi^{3}|)^{6b'-2}},
$
which is bounded if 
$6b'-2\ge 0$ or
$
b'
\geq
\frac13.
$
This completes the proof of Lemma \ref{s-0-theta-2-lemma}.
\,\,
$\square$

\vskip0.05in
\nin
{\bf Estimate in microlocalization II: $|\tau-\xi^{3}|< |\tau_1-\alpha\xi_1^3|$.} 
In this microlocalization, 
using duality and applying the Cauchy--Schwarz inequality  first in $(\xi_1, \tau_1)$ and then in $(\xi, \tau)$, we arrive at
\begin{align*}
 &\Big\| 
\int_{\rr^2}
(\chi_{B_{II}} Q_2)(\xi,\tau,\xi_1,\tau_1)
 c_{f,\alpha}(\xi-\xi_1,\tau-\tau_1) c_{g,\alpha}(\xi_1,\tau_1) 
 d\xi_1 d\tau_1 \Big\|_{L^2_{\xi,\tau}}
 \\
\leq&
\|c_{f,\alpha}\|_{L^2_{\xi,\tau}}\|c_{g,\alpha}\|_{L^2_{\xi,\tau}}  \Big\| \int_{\rr^2} (\chi_{B_{II}} Q_2^2)(\xi,\tau,\xi_1,\tau_1)d\xi d\tau \Big\|_{L^{\infty}_{\xi_1,\tau_1}}^{1/2}.
\end{align*}
Thus, to prove our bilinear estimate \eqref{bilinear-est-s-l} in this microlocalization, it suffices to show 
 the following result for $Q_2$ defined in \eqref{u-s-0-Q2}.
\begin{lemma} 
\label{s-0-theta-3-lemma}
Let $\alpha>0$.
If   $\frac13\leq b'<\frac12$ and $\theta'>\frac12$, then for $\xi_1,\tau_1\in\rr$ we have
\begin{align*}
\Theta_3(\xi_1,\tau_1)
=&
\frac{
\chi_{|\xi_1|\ge1}
}{
(1+|\tau_1-\alpha\xi_1^3|)^{2b'}
}
\int_{\rr^2}
\frac{
\xi^2
}{
(1+|\tau-\xi^{3}|)^{2b'}
}
\frac{
\chi_{|\xi-\xi_1|< 1}
\,
\chi_{B_{II}}(\xi,\tau,\xi_1,\tau_1)
\quad
d\xi
d\tau
}
{
[
(1+|\tau-\tau_1-\alpha(\xi-\xi_1)^3|)^{b'}
+
(1+|\tau-\tau_1|)^{\theta'}
]^2
}
\lesssim 1.
\end{align*}
\end{lemma}
\nin
{\bf Proof of Lemma \ref{s-0-theta-3-lemma}}. 
This proof is similar to the proof of Lemma  \ref{s-0-theta-2-lemma}. So, we omit it.
\,\,
$\square$

\vskip.05in
\noindent
\underline{\bf Proof for  multiplier  $Q_3$.} By symmetry (in convolution writing), we  assume that 
$
|\tau-\tau_1-\alpha(\xi-\xi_1)^{3}|
\leq
|\tau_1-\alpha\xi_1^{3}|.
$
Therefore, to prove  estimate \eqref{bilinear-est-s-l} for $Q_3$  we distinguish two  microlocalizations:
\vskip0.05in
\noindent
{\bf $\bullet$ Microlocalization III: $|\tau_1-\alpha\xi_1^{3}|<|\tau-\xi^{3}|$}. In this situation we define the domain 
\begin{align*}
B_{III}
\doteq
\big\{
(\xi,\tau,\xi_1,\tau_1) \in{\rr}^4: |\tau-\tau_1-\alpha(\xi-\xi)^{3}|\leq|\tau_1-\alpha\xi_1^{3}|<|\tau-\xi^{3}|,
|\xi_1|\ge1,
|\xi-\xi_1|\ge1
\big\}.
\end{align*}
{\bf $\bullet$ Microlocalization IV: $|\tau-\xi^{3}|\leq|\tau_1-\alpha\xi_1^{3}|$}. 
In this case we define the domain 
\begin{align*}
\hskip-0.1in
B_{IV}
\doteq
\big\{(\xi,\tau,\xi_1,\tau_1)\in{\rr}^4: 
\max\{|\tau-\tau_1-\alpha(\xi-\xi)^{3}|,
|\tau-\xi^{3}|
\}\leq|\tau_1-\alpha\xi_1^{3}|,
|\xi_1|
\hskip-0.01in
\ge1,
|\xi-\xi_1|
\hskip-0.01in
\ge1
\big\}.
\end{align*}
{\bf Proof of bilinear estimate in microlocalization III.} In this situation, using the Cauchy--Schwarz inequality with respect to 
$(\xi_1,\tau_1)$, and taking the supremum over  $(\xi,\tau)$ we
arrive at
\begin{align*}
& \Big\| 
\int_{\rr^2}
 (\chi_{B_{III}} Q_3)(\xi,\tau,\xi_1,\tau_1)
 c_{f,\alpha}(\xi-\xi_1,\tau-\tau_1) c_{g,\alpha}(\xi_1,\tau_1) 
 d\xi_1 d\tau_1 \Big\|_{L^2_{\xi,\tau}}
 \nn
 \\
 &
\lesssim
\big\|c_{f,\alpha}\big\|_{L^2_{\xi,\tau}} \big\|c_{g,\alpha}\big\|_{L^2_{\xi,\tau}}
\Big\|
\int_{\rr^2}
 (\chi_{B_{III}} Q_3^2)(\xi,\tau,\xi_1,\tau_1)
d\xi_1 d\tau_1
\Big\|_{L^{\infty}_{\xi,\tau}}^{1/2}.
\end{align*}
Thus, to prove our bilinear estimate \eqref{bilinear-est-s-l} in this microlocalization, it suffices to show 
 the following result for $Q_3$ defined in \eqref{u-s-0-Q3}.
\begin{lemma} 
\label{s-0-lem4}
Let $\alpha>0$ and $\alpha\neq 4$.
If $\frac{5}{12}\leq b'<1/2$, then  for $\xi,\tau\in\rr$
\begin{align}
\label{s-0-theta4}
\Theta_4(\xi,\tau)
\doteq
\frac{\xi^2}{(1+|\tau-\xi^{3}|)^{2b'}}
\int_{\rr^2}\frac{ \chi_{B_{III}}(\xi,\tau,\xi_1,\tau_1)\,\,\,
 \quad d\tau_1 d\xi_1 
}{(1+|\tau-\tau_1-\alpha(\xi-\xi_1)^{3}|)^{2b'}(1+|\tau_1-\alpha\xi_1^{3}|)^{2b'}}
\lesssim 1.
\end{align}
\end{lemma}

\noindent
{\bf Proof of bilinear estimate in microlocalization IV.}
In this case, 
using duality and applying the Cauchy--Schwarz inequality  first in $(\xi_1, \tau_1)$ and then in $(\xi, \tau)$, we get
\begin{align*}
 &\Big\| 
\int_{\rr^2}
(\chi_{B_{IV}} Q_3)(\xi,\tau,\xi_1,\tau_1)
 c_{f,\alpha}(\xi-\xi_1,\tau-\tau_1) c_{g,\alpha}(\xi_1,\tau_1) 
 d\xi_1 d\tau_1 \Big\|_{L^2_{\xi,\tau}}
\\
\lesssim&
\|c_{f,\alpha}\|_{L^2_{\xi,\tau}}\|c_{g,\alpha}\|_{L^2_{\xi,\tau}}  \Big\| \int_{\rr^2} (\chi_{B_{IV}} Q_3^2)(\xi,\tau,\xi_1,\tau_1)d\xi d\tau \Big\|_{L^{\infty}_{\xi_1,\tau_1}}^{1/2}.
\end{align*}
Thus, to prove our bilinear estimate \eqref{bilinear-est-s-l} in this microlocalization, it suffices to show 
 the following result for $Q_3$ as in \eqref{u-s-0-Q3}.
\begin{lemma}
\label{s-0-lem5}
Let $\alpha>0$, $\alpha\neq 1$ and $\alpha\neq 4$.
If $\frac{5}{12}\leq b'<1/2$, then  for $\xi_1,\tau_1\in\rr$
\begin{align}
\label{s-0-theta5}
\Theta_5(\xi_1, \tau_1)
\doteq
\frac{1}{(1+|\tau_1-\alpha\xi_1^{3}|)^{2b'}}
\int_{\rr^2}
\frac{
\chi_{B_{IV}}(\xi,\tau,\xi_1,\tau_1)
\,\,
\xi^2\,
\, d\tau d\xi}{(1+|\tau-\tau_1-\alpha(\xi-\xi_1)^{3}|)^{2b'}(1+|\tau-\xi^{3}|)^{2b'}}
\lesssim 1.
\end{align}
\end{lemma}
\nin
{\bf Proof of Lemma \ref{s-0-lem4}.} 
For the $\tau_1$-integral in \eqref{s-0-theta4}, applying estimate \eqref{eq:calc_5}  with $\ell=\ell'=b'$, $a=\tau-\alpha(\xi-\xi_1)^3$ and $c=\alpha\xi_1^3$, we get
\begin{equation}
\label{s-0-theta4-est}
\Theta_4(\xi,\tau)
\lesssim
\frac{\xi^2}{(1+|\tau-\xi^{3}|)^{2b'}}
\int_{\rr}
\frac{ 
\chi_{B_{III}}(\xi,\tau,\xi_1,\tau_1)
 d\xi_1 
}{(1+|\tau-\alpha(\xi-\xi_1)^{3}-\alpha\xi_1^3|)^{4b'-1}}
=
\frac{\xi^2}{(1+|\tau-\xi^{3}|)^{2b'}}
I(\xi,\tau),
\end{equation}
\vskip-0.08in
\nin
where 
\begin{equation}
\label{I-def-s-0}
I(\xi,\tau)
\doteq
\int_{\rr}
\frac{ 
\chi_{B_{III}}(\xi,\tau,\xi_1,\tau_1)
 d\xi_1 
}{(1+|\tau-\alpha(\xi-\xi_1)^{3}-\alpha\xi_1^3|)^{4b'-1}}
=
\int_{\rr}
\frac{ 
\chi_{B_{III}}(\xi,\tau,\xi_1,\tau_1)
 d\xi_1 
}{(1+|\tau-\xi^3-d_{\alpha}(\xi,\xi_1)|)^{4b'-1}}.
\end{equation}
Here we recall that  $d_{\alpha}(\xi,\xi_1)=-\xi^3+\alpha\xi_1^3+\alpha(\xi-\xi_1)^3$ is the Bourgain quantity of  $u$-equation.
Now  we need the following result.
\begin{lemma}
\label{l-lem-s-0}
Let $\alpha>0$.
If $\frac14<b'<\frac12$, then for all $\xi\neq 0$ and $\tau\in\rr$ we have 
\begin{align}
\label{I-est-s-0}
I(\xi, \tau)
\lesssim
\frac{|\xi|^{-\frac12}(1+|\tau-\xi^3|)^{2-4b'}}{(1+|\tau-\frac{\alpha}{4}\xi^3|)^{\frac12}}.
\end{align}
\end{lemma}
\nin
The proof of Lemma \ref{l-lem-s-0} comes later. 
Here combining estimate \eqref{I-est-s-0} with \eqref{s-0-theta4-est}, we get
\begin{align}
\label{s-0-theta4-est-no1}
\Theta_4(\xi,\tau)
\lesssim&
\frac{\xi^2}{(1+|\tau-\xi^{3}|)^{2b'}}
\frac{|\xi|^{-\frac12}(1+|\tau-\xi^3|)^{2-4b'}}{(1+|\tau-\frac{\alpha}{4}\xi^3|)^{\frac12}}
=
\frac{|\xi|^{3/2}}{(1+|\tau-\xi^{3}|)^{6b'-2}(1+|\tau-\frac{\alpha}{4}\xi^3|)^{\frac12}}.
\end{align}
For $|\xi|\le 1$ and $6b'-2\ge 0$ or $b'\ge \frac13$, from \eqref{s-0-theta4-est-no1} we get $\Theta_4\lesssim1$. Hence, we assume that $|\xi|>1$.
For $\alpha\ne 4$,  we  consider the  two cases possible.

\vskip0.05in
\noindent
$\bullet$ Case 1:  $|\tau-\xi^3|\leq \frac12\big|1-\frac{\alpha}{4}\big||\xi|^3$
\qquad
$\bullet$ Case 2:  $|\tau-\xi^{3}|> \frac12\big|1-\frac{\alpha}{4}\big||\xi|^3$

\vskip0.05in
\noindent
{\bf Estimate in Case 1:  $|\tau-\xi^3|\leq \frac12\big|1-\frac{\alpha}{4}\big||\xi|^3$.}
Then, using the triangle inequality, we have 
$
|\tau-\frac{\alpha}4\xi^3|
\ge
\frac12\big|1-\frac{\alpha}{4}\big||\xi|^3.
$
Hence, from \eqref{s-0-theta4-est-no1} for $6b'-2\ge 0$ or $b'\ge\frac13$ we get
$
\Theta_4(\xi,\tau)
\lesssim
\frac{|\xi|^{3/2}}{1}
\frac{1}{|\xi|^{3/2}}
=
1.
$

\vskip0.05in
\nin
{\bf Estimate in Case 2:  $|\tau-\xi^{3}|> \frac12\big|1-\frac{\alpha}{4}\big||\xi|^3$.} Then
$
|\tau-\xi^{3}|
\gtrsim
|\xi|^{3}
$
and therefore from \eqref{s-0-theta4-est-no1}  we get
$
\Theta_4(\xi,\tau)
\lesssim
\frac{|\xi|^{3/2}}{|\xi^3|^{6b'-2}}
\cdot
\frac{1}{1}
\leq
\frac{1}
{
|\xi|^{3(6b'-2)-3/2}
}.
$
Since $|\xi|> 1$, this quantity is bounded if and only if 
$
3(6b'-2)-3/2
\geq
0
$
or
$
b'
\geq
5/12.
$
This completes the proof of Lemma \ref{s-0-lem4}, once we prove Lemma \ref{l-lem-s-0}.
\,\,
$\square$

\vskip0.05in
\noindent
{\bf Proof of Lemma \ref{l-lem-s-0}.}
Our strategy is to make the change of variables, $
\mu=\mu(\xi_1)=d_{\alpha}(\xi,\xi_1)$. For this to be a good change, we must split the integration at the critical point of $\mu(\xi_1)$. Using the identity \eqref{d-partial-xi1}, i.e. 
$
\p d_{\alpha}/\p\xi_1
=
6\alpha
\xi(\xi_1-\xi/2),
$
we see that  $\xi_1={\xi}/{2}$ is the only  critical point.
Thus, assuming $\xi>0$ ($\xi<0$ is similar) we define the $\xi_1$-intervals
$
\mathcal{I}_1^{\xi_1}
\doteq
(-\infty, \xi/2)
$
and
$
\mathcal{I}_2^{\xi_1}
\doteq
(\xi/2, \infty).
$
Now, making the change of variables $\mu=\mu(\xi_1)=d_{\alpha}(\xi,\xi_1)$ in each one of these intervals
 and defining $\mathcal{I}_k^{\mu}$ be the range of $\mu(\xi_1)$ for $\xi_1\in \mathcal{I}_k^{\xi_1}$, from \eqref{I-def-s-0} we get
\begin{equation}
\label{s-0-xi1-Jk}
I(\xi,\tau)
=
J_1(\xi,\tau)+J_2(\xi,\tau),
\quad
\text{ where }
J_k(\xi,\tau)
=
\int_{\mathcal{I}_k^{\mu}}
\frac{
\chi_{B_{III}}(\xi,\tau,\xi_1,\tau_1)
}{(1+|\tau-\xi^3-\mu|)^{4b'-1}}
\frac{1}{|\mu'(\xi_1)|}
d\mu,
\quad
k=1,2.
\end{equation}
Next, we estimate $|\mu'(\xi_1)|$. Using 
$
|\mu(\xi_1)
-
\mu(\xi/2)|
=
\big|
[-\xi^{3}
+
\alpha\xi_1^{3}
+
\alpha(\xi-\xi_1)^{3}
]
-
[
-\xi^{3}
+
\alpha(\xi/2)^{3}
+
\alpha(\xi-\xi/2)^{3}
]
\big|
=
3\alpha|\xi|(\xi_1-\xi/2)^2,
$
we get
$
|\xi_1-\xi/2|
=
\sqrt{|\mu(\xi_1)-\mu(\xi/2)|/
(3\alpha|\xi|)}.
$
Combining this identity with formula \eqref{d-partial-xi1},  we have
$
|\mu'(\xi_1)|
=
6\alpha|\xi(\xi_1-\xi/2)|
\simeq
|\xi|^{1/2}
|\mu(\xi_1)-\mu(\xi/2)|^{1/2}.
$
Substituting this  into the integral in \eqref{s-0-xi1-Jk}, we obtain
\begin{equation}
\label{s-0-xi1-Jk-no1}
J_{k}(\xi,\tau)
\lesssim
|\xi|^{-\frac12}\int_{\mathcal{I}_k^\mu}\frac{
\chi_{B_{III}}(\xi,\tau,\xi_1,\tau_1)
}{(1+|\tau-\xi^3-\mu|)^{4b'-1}}
\frac{1}{|\mu-\mu(\xi/2)|^{\frac12}}
d\mu,
\quad
k=1,2.
\end{equation}
Moreover making the change of variables $\mu_1=\tau-\xi^3-\mu$ and using $\mu(\xi/2)=-(1-\frac14\alpha)\xi^3$, we get
\begin{align}
\label{I1-split-s-0-no3}
J_{k}(\xi,\tau)
\le
|\xi|^{-\frac12}\int_{\rr}^\infty
\frac{
\chi_{B_{III}}(\xi,\tau,\xi_1,\tau_1)
}{(1+|\mu_1|)^{4b'-1}}
\frac{1}{|\mu_1-(\tau-\frac{\alpha}4\xi^3)|^{\frac12}}
d\mu_1,
\quad
k=1,2.
\end{align}
In addition,  by comparison property  \eqref{d-comparison}  we have 
$
|\tau-\xi^3|
\ge
\frac13
|d_{\alpha}(\xi,\xi_1)|
=
\frac13|\mu|,
$
which implies that 
$
|\mu_1|
\leq
|\mu|
+
|\tau-\xi^3|
\lesssim
|\tau-\xi^3|.
$
Finally, for the $\mu_1$-integral 
in 
\eqref{I1-split-s-0-no3}, 
applying estimate \eqref{eq:calc_4} with $2(1-\ell)=4b'-1$, $a=\tau-\frac{\alpha}{4}\xi^3$ and $c=|\tau-\xi^3|$, we get
$
J_k(\xi, \tau)
\lesssim
\frac{|\xi|^{-\frac12}(1+|\tau-\xi^3|)^{2-4b'}}{(1+|\tau-\frac{\alpha}{4}\xi^3|)^{\frac12}}
$
($k=1,2$),
which combined with  \eqref{s-0-xi1-Jk} gives us
the desired estimate \eqref{I-est-s-0} and completes the proof of Lemma \ref{l-lem-s-0}.
\,\,
$\square$

\vskip0.05in
\nin
{\bf Proof of Lemma \ref{s-0-lem5}.} 
For the $\tau$-integral in \eqref{s-0-theta5}, applying estimate \eqref{eq:calc_5} with  $\ell=\ell'=b'$, $a=\xi^3$ and $c=\tau_1+\alpha(\xi-\xi_1)^3$, we get
\begin{align}
\label{s-0-theta5-est}
\Theta_5(\xi_1,\tau_1)
\lesssim
\frac{1}{(1+|\tau_1-\alpha\xi_1^{3}|)^{2b'}}
\int_{\rr} 
\frac{ 
\chi_{B_{IV}}(\xi,\tau,\xi_1,\tau_1)
\, \xi^2 
}{(1+|\tau_1-\alpha\xi_1^3+d_{\alpha}(\xi,\xi_1)|)^{4b'-1}}
d\xi ,
\end{align}
\begin{minipage}{0.65\linewidth}
where we recall that
$
d_{\alpha}(\xi,\xi_1)
=
-\xi^{3}
+
\alpha\xi_1^{3}
+
\alpha
(\xi-\xi_1)^{3}.
$
Now, our strategy for estimating the $\xi$-integral in \eqref{s-0-theta5-est} is to make the change of variables 
$\mu=\mu(\xi)=d_{\alpha}(\xi,\xi_1)$. For this to be a good change,
we must split the integration at the critical points  of 
$d_{\alpha}(\xi,\xi_1)$, which are
$
p_1
=
\frac{\sqrt{\alpha}}{\sqrt{\alpha}-1}\xi_1,
$
$
p_2
=
\frac{\sqrt{\alpha}}{\sqrt{\alpha}+1}\xi_1.
$
Also, the  inflection point of $\mu(\xi)=d_{\alpha}(\xi,\xi_1)$ is 
$q=\frac{\alpha}{\alpha-1}\xi_1$.
In the case that  $0<\alpha<1$ and $\xi_1>0$,
the graph of $\mu=\mu(\xi)=d_{\alpha}(\xi,\xi_1)$ looks as in Figure \ref{fig:s-0-xi}. 
\end{minipage}
\hskip-0.4in
\begin{minipage}{0.5\linewidth}
\begin{center}
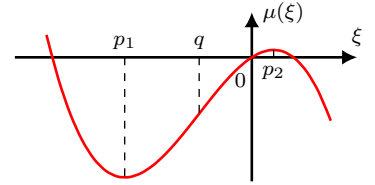

\begin{tikzpicture}[yscale=0.3, xscale=0.7]
%
%
\newcommand\X{0};
\newcommand\Y{0};
\newcommand\FX{8};
\newcommand\FY{8};
\newcommand\R{0.6};
\newcommand\Alpha{0.5};
\newcommand*{\TickSize}{2pt};
%
%
\draw[black,line width=1pt,-{Latex[black,length=2mm,width=2mm]}]
(-4.5,0)
--
(2,0)
node[above]
{\fontsize{\FX}{\FY}\bf \textcolor{black}{$\xi$}};

\draw[black,line width=1pt,-{Latex[black,length=2mm,width=2mm]}]
(0,-5.5)
--
(0,2)
node[right]
{\fontsize{\FX}{\FY}\bf \textcolor{black}{$\mu(\xi)$}};

\draw[line width=1pt, yscale=1,domain=-3.9:1.5,smooth,variable=\x,red]  
plot ({\x},{-(\x)^3+\Alpha+\Alpha*(\x-1)^3});

\draw[dashed,black,line width=0.5pt]
(0,0)
node[xshift=-0.15cm, yshift=-0.3cm]
{\fontsize{\FX}{\FY}\bf $0$}

({sqrt(\Alpha)/(sqrt(\Alpha)-1)},0)
node[yshift=0.2cm]
{\fontsize{\FX}{\FY}\bf $p_1$}
--
({sqrt(\Alpha)/(sqrt(\Alpha)-1)},{(\Alpha*\Alpha-2*(sqrt(\Alpha))^(3))/(1-sqrt(\Alpha))^2})

({sqrt(\Alpha)/(sqrt(\Alpha)+1)},0)
node[yshift=-0.2cm]
{\fontsize{\FX}{\FY}\bf $p_2$}
--
({sqrt(\Alpha)/(sqrt(\Alpha)+1)},{(\Alpha*\Alpha+2*(sqrt(\Alpha))^(3))/(1+sqrt(\Alpha))^2})



({\Alpha/(\Alpha-1)},0)
node[xshift=0cm, yshift=0.2cm]
{\fontsize{\FX}{\FY}\bf $q$}
--
({\Alpha/(\Alpha-1)},{-(\Alpha/(\Alpha-1))^3+\Alpha+\Alpha*(\Alpha/(\Alpha-1)-1)^3})
;

\end{tikzpicture}

\vskip-0.05in
\captionof{figure}{$0<\alpha<1$, $\xi_1>0$}
\label{fig:s-0-xi}
\end{center}
\end{minipage}

\vskip0.04in
Also, we would like to split the $\xi$-integral
  in a way such that $\xi$ is close to the critical points $p_1$, $p_2$ and away from the inflection points $q$. Furthermore, we notice that $|p_1-q|=|p_2-q|=\frac{\sqrt{\alpha}}{|1-\alpha|}|\xi_1|$.
Thus, choosing a small number $\delta=\delta(\alpha)>0$ as follows
\begin{equation}
\label{delta-def}
\delta
=
\delta(\alpha)
\doteq
\frac{1}{10}
\min\big\{
\frac{\sqrt{\alpha}}{|1-\alpha|},
\frac{1}{1000}
\big\},
\end{equation}
\vskip-0.03in
\nin
and   assuming that $\xi_1>0$ and $0<\alpha<1$ (the other cases are similar), we define the $\xi$-intervals
$
\mathcal{I}_{1}^\xi
\doteq
\big(
-\infty
,
p_1-\delta\xi_1
\big),
$
$
\mathcal{I}_{2}^\xi
\doteq
\big(
p_1-\delta\xi_1
,
p_1
\big),
$
$
\mathcal{I}_{3}^\xi
\doteq
\big(
p_1
,
p_1+\delta\xi_1
\big),
$
$
\mathcal{I}_{4}^\xi
\doteq
\big(
p_1+\delta\xi_1
,
p_2-\delta\xi_1
\big),
$
$
\mathcal{I}_{5}^\xi
\doteq
\big(
p_2-\delta\xi_1
,
p_2
\big),
$
$
\mathcal{I}_{6}^\xi
\doteq
\big(
p_2
,
p_2+\delta\xi_1
\big)
$
and
$
\mathcal{I}_{7}^\xi
\doteq
\big(
p_2+\delta\xi_1
,
\infty
\big),
$
which are visualized  in Figure \ref{fig:s-0-xi-axis}

\noindent
\begin{minipage}{1\linewidth}
\begin{center}
\begin{tikzpicture}[yscale=1, xscale=0.8]
%
%
\newcommand\X{0};
\newcommand\Y{0};
\newcommand\FX{11};
\newcommand\FY{11};
\newcommand\R{0.6};
\newcommand*{\TickSize}{2pt};
%
%
\draw[black,line width=1pt,-{Latex[black,length=2mm,width=2mm]}]
(-9,0)
--
(9,0)
node[above]
{\fontsize{\FX}{\FY} $\xi$};

\draw[red]

(-10,0)
node[yshift=0cm,black,line width=1pt]
{\bf \fontsize{\FX}{\FY}$\xi_1> 0:$}
node[xshift=-0.3cm,yshift=-0.5cm,black,line width=1pt]
{\bf \fontsize{\FX}{\FY}$0<\alpha<1$}


(0,0)
node[yshift=-0.4cm]
{\fontsize{\FX}{\FY}$q$}
node[yshift=0cm]
{$\bullet$}
node[xshift=0cm,yshift=0.6cm, black]
{\fontsize{\FX}{\FY}$\mathcal{I}_4^\xi$}

(-7,0)
node[xshift=-1cm,yshift=0.6cm, black]
{\fontsize{\FX}{\FY}$\mathcal{I}_1^\xi$}

(-5,0)
node[yshift=-0.4cm]
{\fontsize{\FX}{\FY}$p_1$}
node[yshift=0cm]
{$\bullet$}
node[xshift=-1cm,yshift=0.6cm, black]
{\fontsize{\FX}{\FY}$\mathcal{I}_2^\xi$}
node[xshift=1cm,yshift=0.6cm, black]
{\fontsize{\FX}{\FY}$\mathcal{I}_3^\xi$}

(-2.5,0)
node[yshift=-0.4cm]
{\fontsize{\FX}{\FY}$p_1+\delta\xi_1$}
node[yshift=0cm]
{$\bullet$}

(-7.5,0)
node[yshift=-0.4cm]
{\fontsize{\FX}{\FY}$p_1-\delta\xi_1$}
node[yshift=0cm]
{$\bullet$}

(5,0)
node[yshift=-0.4cm]
{\fontsize{\FX}{\FY}$p_2$}
node[yshift=0cm]
{$\bullet$}
node[xshift=-1cm,yshift=0.6cm, black]
{\fontsize{\FX}{\FY}$\mathcal{I}_5^\xi$}
node[xshift=1cm,yshift=0.6cm, black]
{\fontsize{\FX}{\FY}$\mathcal{I}_6^\xi$}

(2.5,0)
node[yshift=-0.4cm]
{\fontsize{\FX}{\FY}$p_2-\delta\xi_1$}
node[yshift=0cm]
{$\bullet$}

(7.5,0)
node[yshift=-0.4cm]
{\fontsize{\FX}{\FY}$p_2+\delta\xi_1$}
node[yshift=0cm]
{$\bullet$}

(8,0)
node[xshift=0cm,yshift=0.6cm, black]
{\fontsize{\FX}{\FY}$\mathcal{I}_7^\xi$}

;

\draw [decorate,decoration={brace,mirror,raise=5pt},thick] (-2.5,0) --  (-5,0);

\draw [decorate,decoration={brace,mirror,raise=5pt},thick] (-5,0) --  (-7.5,0);

\draw [decorate,decoration={brace,mirror,raise=5pt},thick] (7.5,0) --  (5,0);

\draw [decorate,decoration={brace,mirror,raise=5pt},thick] (5,0) --  (2.5,0);

\end{tikzpicture}

\vskip-0.15in
\captionof{figure}{$\xi$-intervals}
\label{fig:s-0-xi-axis}
\end{center}
\end{minipage}

\vskip0.05in
\nin
Then, making the change of variables $\mu=\mu(\xi)=d_{\alpha}(\xi,\xi_1)$ in each one of these intervals and defining $\mathcal{I}_k^\mu$ be the range of $\mu(\xi)$ for $\xi\in \mathcal{I}_k^\xi$, from \eqref{s-0-theta5-est} we obtain
\begin{align*}
\Theta_5
\lesssim
\sum\limits_{k=1}^7
J_k(\xi_1,\tau_1)
\text{  with  }
J_k(\xi_1,\tau_1)
\doteq
\frac{1}{(1+|\tau_1-\alpha\xi_1^{3}|)^{2b'}}
\int_{\mathcal{I}_{k}^{\mu}}
\frac{ 
\chi_{B_{IV}}(\xi,\tau,\xi_1,\tau_1) 
\,
\xi^2
}{(1+|\tau_1-\alpha\xi_1^3+\mu|)^{4b'-1}}
\frac{1}{|\mu'(\xi)|}
d\mu.
\end{align*}
For our estimation below, we need the next inequality.
\vskip-0.1in
\begin{lemma}
\label{dist-lem}
If $|a-b|\ge \varepsilon |b|$ with $\varepsilon>0$, then 
$
\displaystyle
|a-b|
\ge
\min\{
\frac12, \frac{\varepsilon}{2}
\}
|a|.
$
\end{lemma}
\nin
{\bf Estimation away from critical points.} 
If $\xi\in \mathcal{I}_k^\xi$,  $k=1,4,7$, then using formula \eqref{d-partial-xi} and the elementary inequality of Lemma \ref{dist-lem}, we get 
$
|\mu'(\xi)|
=
3|\alpha-1|
|
\xi-p_1
|
|
\xi-p_2
|
 \gtrsim \xi^2,
$
which gives us
\begin{align}
\label{s-0-theta5-J1-est}
J_k(\xi_1,\tau_1)
\lesssim
\frac{1}{(1+|\tau_1-\alpha\xi_1^{3}|)^{2b'}}
\int_{\rr}
\frac{ 
\chi_{B_{IV}}(\xi,\tau,\xi_1,\tau_1)
}{(1+|\mu_1|)^{4b'-1}}
d\mu_1,
\quad
 k=1,4,7.
\end{align}
Now, using comparison property \eqref{d-comparison}, i.e. $|\tau_1-\alpha\xi_1^3|\ge \frac13|d_\alpha|=\frac13|\mu|$,  we have $|\mu_1|\le|\mu|+|\tau_1-\alpha\xi_1^3|\lesssim |\tau_1-\alpha\xi_1^3|$, which helps us  estimate the $\mu_1$-integral in \eqref{s-0-theta5-J1-est} as follows
\begin{align*}
J_k(\xi_1,\tau_1)
\lesssim&
\frac{1}{(1+|\tau_1-\alpha\xi_1^{3}|)^{2b'}}
\cdot
(1+|\tau_1-\alpha\xi_1^{3}|)^{2-4b'}
=
\frac{1}{(1+|\tau_1-\alpha\xi_1^{3}|)^{6b'-2}},
\quad
 k=1,4,7,
\end{align*}
and we see that  $J_k$ is bounded if $6b'-2\ge 0$ or $b'\ge 1/3$. 

\vskip0.05in
\nin
{\bf Estimation near critical points.}  Since $p_1, p_2$ are like 
$\xi_1$, and $\xi$ is near $p_1$ and $p_2$, we have 
$|\xi|\simeq |\xi_1|$. Using this and taking  $|\xi|^2$ outside the integral,   we get 
\begin{align}
\label{s-0-theta5-J2}
J_k(\xi_1,\tau_1)
\lesssim
\frac{|\xi_1|^2}{(1+|\tau_1-\alpha\xi_1^{3}|)^{2b'}}
\int_{\mathcal{I}_{k}^{\mu}}
\frac{ 
\chi_{B_{IV}}(\xi,\tau,\xi_1,\tau_1) 
}{(1+|\tau_1-\alpha\xi_1^3+\mu|)^{4b'-1}}
\frac{1}{|\mu'(\xi)|}
d\mu,
\quad
k=2,3,5,6.
\end{align}
\underline{Estimation for $J_2$ and $J_3$.}  We begin by deriving a lower  bound for $|\mu'(\xi)|=|\p d_\alpha/\p \xi|$. Using linear 
approximation  at  $p_1$ with remainder, we have
$
\mu(\xi)
=
\mu(p_1)
+
\mu'(p_1)
(\xi-p_1)
+
\frac{\mu''(\eta_1)}2
(\xi-p_1)^2,
$
where 
$\eta_1$
is some number between
$\xi$
and
$p_1$. Since $p_1$ is a critical point,  we have
$
\mu(\xi)
=
\mu(p_1)
+
\frac{\mu''(\eta_1)}2
(\xi-p_1)^2,
$
which gives
\begin{equation}
\label{xi-p1-dis}
|\xi-p_1|
=
\sqrt{
2[\mu(\xi)-\mu(p_1)]/\mu''(\eta_1)
}.
\end{equation}
Furthermore, since $\xi\in \mathcal{I}_{2}^\xi\cup \mathcal{I}_{3}^{\xi}$,  $\eta_1$ is also in $\mathcal{I}_{2}^\xi\cup \mathcal{I}_{3}^{\xi}$. Moreover, by the choosing of $\delta$, we get $\delta |\xi_1|\le \frac1{10}|p_1-q|=\frac1{10}\frac{\sqrt{\alpha}}{|\alpha-1|}|\xi_1|$, which gives us  $|\eta_1-q|\le \frac{11}{10}|p_1-q|=\frac{11}{10}\frac{\sqrt{\alpha}}{|\alpha-1|}|\xi_1|$. Now, using formula \eqref{d-2partial-xi}, i.e. 
$
\mu''(\eta_1)
=
6(1-\alpha)
\big(
\eta_1-q
\big),
$
we get $|\mu''(\eta_1)|\lesssim |\xi_1|$,
which combined with \eqref{xi-p1-dis} implies 
\begin{equation}
\label{s-0-xi-p1-dist}
|\xi-p_1|
\gtrsim
\sqrt{
|\mu(\xi)-\mu(p_1)|/|\xi_1|
}.
\end{equation}
Also, for $\xi\in \mathcal{I}_2^\xi\cup \mathcal{I}_3^\xi$, we have $|\xi-p_2|\ge \delta|\xi_1|$. Combining this with estimate \eqref{s-0-xi-p1-dist}, by \eqref{d-partial-xi} we get 
\begin{equation}
\label{s-0-mu-der-est}
|\mu'(\xi)|
=
3|\alpha-1|
|
\xi-p_1
|
|
\xi-p_2
|
\gtrsim
\sqrt{
|\mu(\xi)-\mu(p_1)|/|\xi_1|
}
|\xi_1|
=
|\xi_1|^{1/2} \sqrt{
|\mu(\xi)-\mu(p_1)|
}.
\end{equation}
Using estimate \eqref{s-0-mu-der-est}  and making the change of variables $\mu_1=\mu+\tau_1-\alpha\xi_1^3$, from \eqref{s-0-theta5-J2} we get
\begin{align}
\label{s-0-case2-theta-3-est2}
J_k(\xi_1,\tau_1)
\lesssim
\frac{|\xi_1|^{3/2}}{(1+|\tau_1-\alpha\xi_1^{3}|)^{2b'}}
\int_{\rr}
\frac{ 
\chi_{B_{IV}}(\xi,\tau,\xi_1,\tau_1) 
}{(1+|\mu_1|)^{4b'-1}}
\frac{1}{\sqrt{
|\mu_1-(\tau_1-\alpha\xi_1^3)-\mu(p_1)|
}}
d\mu_1,
\quad
k=2,3.
\end{align}
Moreover,  using comparison property \eqref{d-comparison},  we have $|\mu_1|\lesssim |\tau_1-\alpha\xi_1^3|$. For the $\mu_1$-integral in \eqref{s-0-case2-theta-3-est2},
applying estimate \eqref{eq:calc_4} with $2(1-\ell)=4b'-1$, $a=\tau_1-\alpha\xi_1^3+\mu(p_1)$ and $c=|\tau_1-\alpha\xi_1^3|$, 
we get
\begin{align}
\label{s-0-case1-theta-3-est-no2}
J_k(\xi_1,\tau_1)
\lesssim
\frac{|\xi_1|^{3/2}}{
(1+|\tau_1-\alpha\xi_1^{3}|)^{6b'-2}
(1+|\tau_1-\alpha\xi_1^3+\mu(p_1)|)^{\frac12}},
\quad
k=2,3.
\end{align}
We observe that the difference between the two factors in the denominator is $\mu(p_1)$,  which is
$
\mu(p_1)
=
-p_1^3
+
\alpha\xi_1^3
+
\alpha(p_1-\xi_1)^3
=
\frac{\alpha^2-2\alpha^{3/2}}{(1-\sqrt{\alpha})^2}
\xi_1^3.
$
Noticing that for $|\xi_1|\ge 1$,  $\mu(p_1)=0$ if and only if $\alpha=4$. So  if $\alpha\neq 4$, then like the estimate for the quantity in \eqref{s-0-theta4-est-no1},   by considering the following  two cases:
$|\tau_1-\alpha\xi_1^3|\leq \frac12|\mu(p_1)|$
and
$|\tau_1-\alpha\xi_1^{3}|> \frac12|\mu(p_1)|$, 
we get $J_k\lesssim 1$ ($k=2,3$) for $b'\ge \frac5{12}$.

\vskip0.05in
\nin
\underline{Estimation for $J_5$ and $J_6$.}  Using the linear approximation at $p_2$, we get  (similar to estimate \eqref{s-0-case1-theta-3-est-no2})
\begin{align*}
J_k(\xi_1,\tau_1)
\lesssim
\frac{|\xi_1|^{3/2}}{
(1+|\tau_1-\alpha\xi_1^{3}|)^{6b'-2}
(1+|\tau_1-\alpha\xi_1^3+\mu(p_2)|)^{\frac12}},
\quad
k=5,6.
\end{align*}
Also, since $\mu(p_2)=\frac{\alpha^2+2\alpha^{3/2}}{(1+\sqrt{\alpha})^2}
\xi_1^3\neq 0$ for $|\xi_1|>1$, by considering the two cases $|\tau_1-\alpha\xi_1^3|\leq \frac12|\mu(p_2)|$ and  $|\tau_1-\alpha\xi_1^{3}|> \frac12|\mu(p_2)|$, we get $J_5\lesssim 1$ and $J_6\lesssim 1$  for $b'\ge \frac{5}{12}$.
This completes the proof of Lemma \ref{s-0-lem5}.
\,\,
$\square$

\vskip0.05in
\nin
{\bf Conclusion.} In the case $s\ge 0$,  Lemmas \ref{s-0-theta-1-lemma}--\ref{s-0-lem5} imply bilinear estimate \eqref{bilinear-est-s-0},
if  $\alpha\neq 1$ and $\alpha\neq 4$.

%
%
%
%
%
%
%
%
%
%
\subsection{$s<0$ with $\alpha>4$}
%
%
First, we observe that  to prove the bilinear estimate \eqref{bi-est-X-1}, 
it suffices to show it for $|\xi_1|$ and  $|\xi-\xi_1|$ away from 0,
that is
\begin{equation}
\label{above-1-cond}
|\xi_1|> 1
\qquad
\text{and}
\qquad
|\xi-\xi_1|>1.
\end{equation}
In fact, if $|\xi_1|\leq 1$ or $|\xi-\xi_1|\leq 1$ then
our estimate \eqref{bilinear-est-L2-form-no1} is reduced to the case that $s=0$ by
\begin{equation}
\label{sob-mult}
 \frac{(1+|\xi|)^s}{
(1+|\xi_1|)^{s}
(1+|\xi-\xi_1|)^{s}
}
\lesssim 
1,
\end{equation}
which is  the combination of Sobolev multipliers appearing in $Q'$ defined in \eqref{Q'-def}.
For $|\xi_1|< 1$ inequality \eqref{sob-mult} follows from
$
(1+|\xi-\xi_1|)(1+|\xi_1|)\leq (1+|\xi|+1)\cdot2 
= 
4(1+|\xi|/2)\leq 4(1+|\xi|)
$
and $s<0$.
For $|\xi-\xi_1|< 1$ this inequality \eqref{sob-mult} follows from
$
(1+|\xi-\xi_1|)(1+|\xi_1|)\leq 2\cdot(1+|\xi|+1) 
= 
4(1+|\xi|/2)\leq 4(1+|\xi|)
$
and $s<0$.
Then  $Q'$ in estimate \eqref{bilinear-est-L2-form-no1} is replaced by $Q_0$ and so we are reduced to $s=0$.

\vskip0.05in
Now, using assumption \eqref{above-1-cond} we can replace 
$
(1+|\xi_1|)^{-s}
(1+|\xi-\xi_1|)^{-s}
$
by
$
|\xi_1(\xi-\xi_1)|^{-s}
$.
Also, dropping 
$
\chi_{|\xi_1|< 1}(1+|\tau_1|)^{\theta'}
$
and
$
\chi_{|\xi-\xi_1|< 1}(1+|\tau-\tau_1|)^{\theta'},
$
we see that  to prove estimate \eqref{bilinear-est-L2-form-no1} it suffices to  do it with  $Q'$ being replaced with
\begin{align}
\label{Q1-modify-no1-s-neg}
Q_1(\xi,\tau,\xi_1,\tau_1)
\doteq
\frac{|\xi|}{(1+|\tau-\xi^{3}|)^{b'}}
\frac{(1+|\xi|)^s|\xi_1(\xi-\xi_1)|^{-s}}{
(1+|\tau_1-\alpha\xi_1^3|)^{b'}
}
\frac{1}
{
(1+|\tau-\tau_1-\alpha(\xi-\xi_1)^3|)^{b'}
}.
\end{align}
Furthermore, by symmetry (in convolution writing), we assume that 
$
|\tau-\tau_1-\alpha(\xi-\xi_1)^{3}|
\leq
|\tau_1-\alpha\xi_1^{3}|.
$
Also, 
we consider the following two microlocalizations:
\vskip0.05in
\noindent
{\bf $\bullet$ Microlocalization I: $|\tau_1-\alpha\xi_1^{3}|\leq|\tau-\xi^{3}|$}. In this situation we define the domain 
\begin{align*}
E_I
\doteq
\big\{
(\xi,\tau,\xi_1,\tau_1) \in{\rr}^4: |\tau-\tau_1-\alpha(\xi-\xi)^{3}|\leq|\tau_1-\alpha\xi_1^{3}|<|\tau-\xi^{3}|,
|\xi_1|\ge1,
|\xi-\xi_1|\ge1
\big\}.
\end{align*}
{\bf $\bullet$ Microlocalization II: $|\tau-\xi^{3}|\leq|\tau_1-\alpha\xi_1^{3}|$.} 
In this case we define the domain 
\begin{align*}
E_{II}
\doteq
\big\{
(\xi,\tau,\xi_1,\tau_1)\in{\rr}^4: 
\max\{|\tau-\tau_1-\alpha(\xi-\xi)^{3}|,
|\tau-\xi^{3}|
\}\leq|\tau_1-\alpha\xi_1^{3}|,
|\xi_1|
\ge1,
|\xi-\xi_1|
\ge1
\big\}.
\end{align*}
{\bf Proof of bilinear estimate in  microlocalization I.} 
Here our bilinear estimate \eqref{bilinear-est-L2-form-no1} is reduced to
\begin{align}
\label{bilinear-est-neg-EI}
 \Big\| 
 \int_{\rr^2} 
 (\chi_{E_I} Q_1)(\xi,\tau,\xi_1,\tau_1)
 c_{f,\alpha}(\xi-\xi_1,\tau-\tau_1) c_{g,\alpha}(\xi_1,\tau_1) 
 d\xi_1 d\tau_1 \Big\|_{L^2_{\xi,\tau}}
\lesssim
\big\|c_{f,\alpha}\big\|_{L^2_{\xi,\tau}} \big\|c_{g,\alpha}\big\|_{L^2_{\xi,\tau}}.
\end{align}
Using the Cauchy--Schwarz inequality with respect to 
$(\xi_1,\tau_1)$, and taking the supremum over  $(\xi,\tau)$ we
obtain the following estimate for the left-hand side of \eqref{bilinear-est-neg-EI}
\begin{align*}
&\Big\| 
 \int_{\rr^2} 
 (\chi_{E_I} Q_1)(\xi,\tau,\xi_1,\tau_1)
 c_{f,\alpha}(\xi-\xi_1,\tau-\tau_1) c_{g,\alpha}(\xi_1,\tau_1) 
 d\xi_1 d\tau_1 \Big\|_{L^2_{\xi,\tau}}
 \nn
 \\
 &
\lesssim
\big\|c_{f,\alpha}\big\|_{L^2_{\xi,\tau}} \big\|c_{g,\alpha}\big\|_{L^2_{\xi,\tau}}
\Big\|
\int_{\rr^2}
 (\chi_{E_{I}} Q_1^2)(\xi,\tau,\xi_1,\tau_1)
d\xi_1 d\tau_1
\Big\|_{L^{\infty}_{\xi,\tau}}^{\frac12}.
\end{align*}
So,  to prove bilinear estimate \eqref{bilinear-est-neg-EI},  it suffices to show 
 the following result for $Q_1$ defined by  \eqref{Q1-modify-no1-s-neg}.
%
%
%
%
%
%
%
%
%
%
%
%
\begin{lemma} 
\label{s-neg-lem1} 
Let   $\alpha>4$.
If $\max\{\frac{15-4s}{36},\frac{7}{18}\}<  b'<\frac12$  and $s>-\frac34$
then  for $\xi,\tau\in\rr$
\begin{align}
\label{s-neg-theta1}
\Theta_1(\xi,\tau)
\doteq
\frac{\xi^2(1+|\xi|)^{2s}}{(1+|\tau-\xi^{3}|)^{2b'}}
\int_{\rr^2} \frac{ \chi_{E_I}(\xi,\tau,\xi_1,\tau_1)|\xi_1(\xi-\xi_1)|^{-2s}\,\,\,
 \quad d\tau_1 d\xi_1 
}{(1+|\tau-\tau_1-\alpha(\xi-\xi_1)^{3}|)^{2b'}(1+|\tau_1-\alpha\xi_1^{3}|)^{2b'}}
\lesssim 1.
\end{align}
\end{lemma}
\nin
{\bf Proof of bilinear estimate in  microlocalization II.}
In this situation,  using duality and applying the Cauchy--Schwarz inequality  first in $(\xi_1, \tau_1)$ and then in $(\xi, \tau)$, we get
\begin{align*}
 & \Big\| 
 \int_{\rr^2} 
 (\chi_{E_{II}} Q_1)(\xi,\tau,\xi_1,\tau_1)
 c_{f,\alpha}(\xi-\xi_1,\tau-\tau_1) c_{g,\alpha}(\xi_1,\tau_1) 
 d\xi_1 d\tau_1 \Big\|_{L^2_{\xi,\tau}}
 \\
 \lesssim&
\|c_{f,\alpha}\|_{L^2_{\xi,\tau}}\|c_{g,\alpha}\|_{L^2_{\xi,\tau}} \Big\| 
\int_{\rr^2} (\chi_{E_{II}} Q_1^2)(\xi,\tau,\xi_1,\tau_1)d\xi d\tau 
\Big\|_{L^{\infty}_{\xi_1,\tau_1}}^{\frac12}.
\end{align*}
Thus, to prove estimate  \eqref{bilinear-est-neg-EI} in this case it suffices to prove the following lemma
with $Q_1$ as in \eqref{Q1-modify-no1-s-neg}.
\begin{lemma}
\label{s-neg-lem2} 
Let   $\alpha>4$.
If $\max\{\frac{15-4s}{36},\frac{7}{18}\}< b'<1/2$ and $s>-\frac34$, then  for $\xi_1,\tau_1\in\rr$
\begin{align}
\label{s-neg-theta2}
\Theta_2(\xi_1, \tau_1)
\doteq
\frac{1}{(1+|\tau_1-\alpha\xi_1^{3}|)^{2b'}}
\int_{\rr^2}
\frac{
\chi_{E_{II}}(\xi,\tau,\xi_1,\tau_1)
\,\,
\xi^2
(1+|\xi|)^{2s}
|\xi_1(\xi-\xi_1)|^{-2s}
\,
\, d\tau d\xi}{(1+|\tau-\tau_1-\alpha(\xi-\xi_1)^{3}|)^{2b'}(1+|\tau-\xi^{3}|)^{2b'}}
\lesssim 1.
\end{align}
\end{lemma}
\nin
{\bf Proof of Lemma \ref{s-neg-lem1} .} For the $\tau_1$-integral in \eqref{s-neg-theta1},
applying calculus estimate \eqref{eq:calc_5} with  $\ell=\ell'=b'$, $a=\tau-\alpha(\xi-\xi_1)^3$ and $c=\alpha\xi_1^3$, we get
\begin{align}
\label{s-neg-theta1-est}
\Theta_1(\xi,\tau)
\lesssim
\frac{\xi^2(1+|\xi|)^{2s}}{(1+|\tau-\xi^{3}|)^{2b'}}
\int_\rr
\frac{ 
\chi_{E_I}(\xi,\tau,\xi_1,\tau_1)
|\xi_1(\xi-\xi_1)|^{-2s}
 d\xi_1 
}{(1+|\tau-\xi^3-d_{\alpha}(\xi,\xi_1)|)^{4b'-1}},
\end{align}
\begin{minipage}{0.64 \linewidth}
where we recall that
$
d_{\alpha}(\xi,\xi_1)
=
-\xi^{3}
+
\alpha\xi_1^{3}
+
\alpha
(\xi-\xi_1)^{3}.
$
Now, our strategy for estimating the above $\xi_1$-integral is to make the change of variables,
$
\mu=\mu(\xi_1)=d_{\alpha}(\xi,\xi_1).
$
Since 
by  formula \eqref{d-partial-xi1},  i.e.
$
\p d_{\alpha}/\p\xi_1
=
6\alpha
\xi(\xi_1-\xi/2),
$
we need to split  the $\xi_1$-integral  at the critical point $\xi_1=1/2\xi$, so that this change of variables makes sense.
In the case that  $\alpha>4$ and $\xi>0$,
the graph of $\mu=\mu(\xi_1)=d_{\alpha}(\xi,\xi_1)$ looks as in Figure \ref{fig:s-neg-xi1}.
\end{minipage}
\hskip-0.2in
\begin{minipage}{0.45\linewidth}
\begin{center}
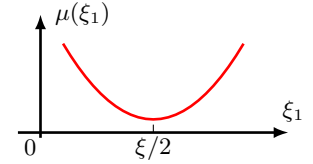

\begin{tikzpicture}[yscale=0.35, xscale=3,every node/.style={scale=1}]
%
%
\newcommand\X{0};
\newcommand\Y{0};
\newcommand\FX{9};
\newcommand\FY{9};
\newcommand\R{0.6};
\newcommand*{\TickSize}{2pt};
\newcommand*{\Num}{4};
\newcommand*{\Alphalpha}{6};
\newcommand*{\Dist}{2};
\newcommand*{\Step}{\Dist/\Num};
%
%
%
\draw[line width=1pt,-{Latex[length=2mm]}]
(-0.1,0)
--
(1.1,0)
node[above]
{\fontsize{\FX}{\FY} $\xi_1$}
;

\draw[line width=1pt,-{Latex[length=2mm]}]
(0,-0.2)
--
(0,4.5)
node[xshift=0.5cm]
{\fontsize{\FX}{\FY} $\mu (\xi_1)$};

\draw[red, line width=1pt, domain=0.1:0.9,variable=\x]
plot
({\x},{-1+\Alphalpha*(\x)^3+\Alphalpha*(1-\x)^3})
;

\draw[black,dashed]
(0.5,0)
node[yshift=-0.2cm]
{\fontsize{\FX}{\FY}$ \xi/2$}
--
(0.5,{-1+\Alphalpha*(0.5)^3+\Alphalpha*(1-0.5)^3})
;

\draw[black]
(0,0)
node[xshift=-0.2cm,yshift=-0.2cm]
{\fontsize{\FX}{\FY} $0$}
;

\end{tikzpicture}

\vskip-0.15in
\captionof{figure}{$\alpha>4$ and $\xi>0$}
\label{fig:s-neg-xi1}
\end{center}
\end{minipage}

\vskip0.05in
\noindent
Now, for $\xi>0$ we define the $\xi_1$-intervals 
($\xi<0$  is similar)
$
\mathcal{I}_1^{\xi_1}
\doteq
(-\infty, -\frac34\xi),
$
$
\mathcal{I}_2^{\xi_1}
\doteq
(-\frac34\xi, \frac34\xi)
$
and
$
\mathcal{I}_3^{\xi_1}
\doteq
(\frac34\xi, \infty).
$
Thus, from estimate \eqref{s-neg-theta1-est} we obtain
\begin{align}
\label{s-neg-theta1-est-deco}
\Theta_1(\xi,\tau)
\lesssim
\sum\limits_{k=1}^3
J_k
\text{ with }
J_k(\xi,\tau)
\doteq
\frac{\xi^2(1+|\xi|)^{2s}}{(1+|\tau-\xi^{3}|)^{2b'}}
\int_{\mathcal{I}_k^{\xi_1}}
\frac{ 
\chi_{E_I}(\xi,\tau,\xi_1,\tau_1)
|\xi_1(\xi-\xi_1)|^{-2s}
 d\xi_1 
}{(1+|\tau-\xi^3-d_{\alpha}(\xi,\xi_1)|)^{4b'-1}}.
\end{align}
{\bf Estimation near critical points (estimate for $J_2$).} Since $\xi_1\in \mathcal{I}_2^{\xi_1}$, we have $1\le |\xi_1|\lesssim |\xi|$ and $|\xi-\xi_1|\lesssim |\xi|$, which gives us $|\xi_1(\xi-\xi_1)|^{-2s}\lesssim |\xi|^{-4s}$. Using this estimate,  from \eqref{s-neg-theta1-est-deco} we get 
\begin{align}
\label{s-neg-theta1-est-no1}
J_2(\xi,\tau)
\lesssim
\frac{|\xi|^{2-2s}}{(1+|\tau-\xi^{3}|)^{2b'}}
\int_{-3/4\xi}^{3/4\xi}
\frac{ 
\chi_{E_I}(\xi,\tau,\xi_1,\tau_1)
 d\xi_1
}{(1+|\tau-\xi^3-d_{\alpha}(\xi,\xi_1)|)^{4b'-1}}
\lesssim
\frac{|\xi|^{2-2s}}{(1+|\tau-\xi^{3}|)^{2b'}}
I(\xi,\tau),
\end{align}
where $I(\xi,\tau)$ is  defined by \eqref{I-def-s-0},
which was estimated in Lemma \ref{l-lem-s-0}.
Combining estimate \eqref{I-est-s-0} with \eqref{s-neg-theta1-est-no1}, we obtain
\begin{align}
\label{s-neg-theta1-est-no2}
J_2(\xi,\tau)
\lesssim&
\frac{|\xi|^{2-2s}}{(1+|\tau-\xi^{3}|)^{2b'}}
\frac{|\xi|^{-\frac12}(1+|\tau-\xi^3|)^{2-4b'}}{(1+|\tau-\frac{\alpha}{4}\xi^3|)^{\frac12}}
=
\frac{|\xi|^{3/2-2s}}{(1+|\tau-\xi^{3}|)^{6b'-2}}
\frac{1}{(1+|\tau-\frac{\alpha}{4}\xi^3|)^{\frac12}}.
\end{align}
Furthermore, by property  \eqref{d-comparison}, i.e. 
$
|\tau-\xi^3|
\ge
\frac13
|d_{\alpha}(\xi,\xi_1)|,
$
and
using  estimate  \eqref{d-bound-xi} for $d_{\alpha}$, we get
  $|\tau-\xi^3|\gtrsim |\xi|^3$.
Combining this  with \eqref{s-neg-theta1-est-no2} and dropping $(1+|\tau-\frac{\alpha}{4}\xi^3|)^{\frac12}$ we get
$
J_2(\xi,\tau)
\lesssim
\frac{|\xi|^{3/2-2s}}{|\xi^{3}|^{6b'-2}},
$
which is bounded if $3(6b'-2)\ge(3/2-2s)$ or 
$
b'
\geq
\frac{15-4s}{36}.
$
For $b'<\frac12$, it suffices to have 
$
s>-3/4.
$
This completes the estimate for $J_2(\xi,\tau)$. 

\vskip0.05in
\noindent
{\bf Estimation away from critical points (estimate for $J_1$ and $J_3$).} Making the change of variables $\mu=\mu(\xi_1)=d_{\alpha}(\xi,\xi_1)$ and  
defining $\mathcal{I}_k^\mu$ be the range of $\mu(\xi_1)$ for $\xi_1\in \mathcal{I}_k^{\xi_1}$, we get 
\begin{align}
\label{s-neg-theta1-case2-no1}
J_k(\xi,\tau)
\lesssim
\frac{\xi^2
(1+|\xi|)^{2s}}{(1+|\tau-\xi^{3}|)^{2b'}}
\int_{\mathcal{I}_k^{\mu}}
\frac{ 
\chi_{E_I}(\xi,\tau,\xi_1,\tau_1)\, |\xi_1(\xi-\xi_1)|^{-2s}
}{(1+|\tau-\xi^3-\mu|)^{4b'-1}}
\frac{1}{|\mu'|}
 d\mu,
 \,\,
 k=1,3.
\end{align}
For  $\xi_1<-\frac34\xi$ or $\xi_1>\frac34\xi$, we have $|\xi-\xi_1|\lesssim |\xi_1|$.
Also,
using the formula \eqref{d-partial-xi1} and using inequality of Lemma \ref{dist-lem}   we have 
$
|\mu'(\xi_1)|
=
6\alpha
|\xi(\xi_1-\xi/2)|
\gtrsim
|\xi|
|\xi_1|,
$
which combined  with   \eqref{s-neg-theta1-case2-no1} implies that
\begin{align}
\label{s-neg-theta1-case2-no3}
J_k(\xi,\tau)
\lesssim
 \frac{
|\xi|
(1+|\xi|)^{2s}}{(1+|\tau-\xi^{3}|)^{2b'}}
 \int_{\rr}
\frac{ 
|\xi_1|^{-4s-1}
\chi_{E_I}(\xi,\tau,\xi_1,\tau_1)\, 
}{(1+|\tau-\xi^3-\mu|)^{4b'-1}}
 d{\mu},
 \quad
 k=1,3.
\end{align}
Next, we consider the following two cases

\vskip0.05in
\noindent
$\bullet$
Case 1:
$0>s\ge -\frac14$
\qquad
$\bullet$
Case 2:
$s< -\frac14$
 
\vskip0.05in
\noindent
\underline{Proof in Case 1.} Using $-4s-1\le 0$ and $|\xi_1|\ge 1$ we have $|\xi_1|^{-4s-1}\le 1$. Combining this estimate with  
 \eqref{s-neg-theta1-case2-no3} and making the change of variables $\mu_1=\tau-\xi^3-\mu$, we get 
\begin{align}
\label{s-neg-theta1-case2-no5}
J_k(\xi,\tau)
\lesssim
\frac{
|\xi|
(1+|\xi|)^{2s}}{(1+|\tau-\xi^{3}|)^{2b'}}
 \int_{\rr}
\frac{ 
\chi_{E_I}(\xi,\tau,\xi_1,\tau_1)\, 
}{(1+|\mu_1|)^{4b'-1}}
 d{\mu_1},
 \quad
 k=1,3.
 \end{align}
Also, in this microlocalization using comparison property 
 \eqref{d-comparison} we get $|\tau-\xi^3|\ge \frac13|d_\alpha|=\frac13|\mu|$, which implies that
$|\mu_1|\le |\tau-\xi^3|+|\mu|\lesssim |\tau-\xi^3|$. This helps us integrate $\mu_1$ in \eqref{s-neg-theta1-case2-no5}
\begin{align}
\label{s-neg-theta1-case2-no6}
J_k(\xi,\tau)
\lesssim&
  \frac{
|\xi|
(1+|\xi|)^{2s}}{(1+|\tau-\xi^{3}|)^{2b'}}
\cdot
(1+|\tau-\xi^{3}|)^{2-4b'}
=
  \frac{
|\xi|
(1+|\xi|)^{2s}}{(1+|\tau-\xi^{3}|)^{6b'-2}},
\quad
k=1,3.
 \end{align}
Furthermore, using $|\tau-\xi^3|\ge \frac13|d_\alpha|$ and  estimate \eqref{d-bound-xi} i.e. $|d_\alpha|\ge \frac{\alpha-4}{4}
|\xi|^3$, we get  $|\tau-\xi^3|\gtrsim |\xi|^3$. Combining  this with estimate 
\eqref{s-neg-theta1-case2-no6}, we get
$
\Theta_1(\xi,\tau)
\lesssim
\frac{
|\xi|
(1+|\xi|)^{2s}}{(1+|\xi^3|)^{6b'-2}}
\lesssim
  \frac{
1}{(1+|\xi|)^{3(6b'-2)-1-2s}}, 
$
 which is bounded if $3(6b'-2)-1-2s\ge 0$ or
$
 b'\ge \frac{7}{18}+\frac19s.
$
Since
$
0>
s
\ge
-
\frac14,
$
it suffices to have $b'\ge 7/18$.
This completes the proof for $J_1(\xi,\tau)$ and $J_3(\xi,\tau)$ in case 1.

\vskip0.05in
\noindent
\underline{Proof in Case 2.} Using property  \eqref{d-comparison}, i.e. $|\tau-\xi^3|\ge \frac13|d_\alpha|$, and by estimate \eqref{d-bound-xi1}, i.e. $|d_\alpha|\ge\frac{3\alpha(\alpha-4)}{4(\alpha-1)}
|\xi||\xi_1|^2$, we get $|\tau-\xi^3|\gtrsim |\xi|\xi_1^2$ or $\xi_1^2\lesssim \frac{|\tau-\xi^3|}{|\xi|}
$. Thus, 
for $-4s-1>0$  we have 
$
|\xi_1|^{-4s-1}
\lesssim
(|\tau-\xi^3|/|\xi|
)^{-2s-1/2}.
$
Using this estimate  and  the change of variables $\mu_1=\tau-\xi^3-\mu$, from  \eqref{s-neg-theta1-case2-no3} we get 
\begin{align}
\label{s-neg-theta1-case2-2-no1}
J_k(\xi,\tau)
\lesssim
\frac{
|\xi|^{3/2+2s}
(1+|\xi|)^{2s}}{(1+|\tau-\xi^{3}|)^{2b'+2s+1/2}}
 \int_{\rr}
\frac{ 
\chi_{E_I}(\xi,\tau,\xi_1,\tau_1)\, 
}{(1+|\mu_1|)^{4b'-1}}
 d{\mu_1},
 \quad
 k=1,3.
 \end{align}
Also,  by property  \eqref{d-comparison}, we have 
$|\mu_1|\lesssim |\tau-\xi^3|$, which helps us integrate $\mu_1$ in \eqref{s-neg-theta1-case2-2-no1}
\begin{align}
\label{s-neg-theta1-case2-2-no2}
J_k(\xi,\tau)
\lesssim&
\frac{
|\xi|^{3/2+2s}
(1+|\xi|)^{2s}}{(1+|\tau-\xi^{3}|)^{2b'+2s+1/2}}
\cdot
(1+|\tau-\xi^{3}|)^{2-4b'}
=
\frac{
|\xi|^{3/2+2s}
(1+|\xi|)^{2s}}{(1+|\tau-\xi^{3}|)^{6b'-3/2+2s}},
\quad
k=1,3.
 \end{align}
Furthermore, using $|\tau-\xi^3|\ge \frac13|d_\alpha|$ and  estimate \eqref{d-bound-xi}, we get  $|\tau-\xi^3|\gtrsim |\xi|^3$. Combining this estimate with \eqref{s-neg-theta1-case2-2-no2}, for
$
6b'-3/2+2s
\ge 
0
$
or 
$
b'
\ge
1/4
-
1/3s,
$ 
we get
\begin{equation}
\label{s-neg-theta1-case2-2-no3}
J_k(\xi,\tau)
\lesssim
  \frac{
|\xi|^{3/2+2s}
(1+|\xi|)^{2s}}{(1+|\xi^{3}|)^{6b'-3/2+2s}},
\quad
k=1,3.
\end{equation} 
In addition, for $3/2+2s\ge0$ or $s\ge-3/4$, we have $|\xi|^{3/2+2s}\lesssim (1+|\xi|)^{3/2+2s}$. Combining this estimate with \eqref{s-neg-theta1-case2-2-no3}, we get 
$
\Theta_1(\xi,\tau)
\lesssim
  \frac{
1}{(1+|\xi|)^{3(6b'-3/2+2s)-3/2-4s}}
\lesssim
  \frac{
1}{(1+|\xi|)^{18b'-6+2s}}, 
$
 which is bounded if $18b'-6+2s\ge 0$ or
 $
b'\ge \frac13-\frac19s.
$
For $b'<\frac12$, it suffices to have $
s+3/2>0$ or $s>-3/2$. %
This completes  the  proof of Lemma \ref{s-neg-lem1} .
\,\,
$\square$ 

\vskip0.05in
\nin
{\bf Proof of Lemma \ref{s-neg-lem2}.} For the $\tau$-integral in  \eqref{s-neg-theta2},  applying estimate \eqref{eq:calc_5}  with $\ell=\ell'=b'$, $a=\xi^3$ and $c=\tau_1+\alpha(\xi-\xi_1^3)$, we get
\begin{align}
\label{s-neg-theta2-est-no1}
\Theta_2(\xi_1, \tau_1)
\lesssim
\frac{1}{(1+|\tau_1-\alpha\xi_1^{3}|)^{2b'}}
\int_\rr
\frac{
\chi_{E_{II}}(\xi,\tau,\xi_1,\tau_1)
|\xi_1(\xi-\xi_1)|^{-2s}
\xi^2
(1+|\xi|)^{2s}
\,
\, d\xi}{
(1+|\tau_1-\alpha\xi_1^3+d_{\alpha}(\xi,\xi_1)|)^{4b'-1}
}.
\end{align}
\begin{minipage}{0.7\linewidth}
Like before, we split the $\xi$-integral near critical points and away from critical points. Here we recall that the two critical points are given by 
$
p_1
=
\frac{\sqrt{\alpha}}{\sqrt{\alpha}-1}\xi_1,
$
$
p_2
=
\frac{\sqrt{\alpha}}{\sqrt{\alpha}+1}\xi_1,
$
and the inflection point is $q=\frac{\alpha}{\alpha-1}\xi_1$.
In the case that   $\xi_1>0$,
the graph of $\mu=\mu(\xi)=d_{\alpha}(\xi,\xi_1)$ looks as in Figure \ref{fig:s-neg-xi}.  Also, using the small number $\delta$ given by \eqref{delta-def}, for $\xi_1>0$ (for $\xi_1<0$, it is similar) we define the $\xi$-intervals %
$
\mathcal{I}_{1}^\xi
\doteq
\big(
-\infty
,
p_2-\delta\xi_1
\big),
$
$
\mathcal{I}_{2}^\xi
\doteq
\big(
p_2-\delta\xi_1
,
p_2
\big),
$
$
\mathcal{I}_{3}^\xi
\doteq
\big(
p_2
,
p_2+\delta\xi_1
\big),
$
$
\mathcal{I}_{4}^\xi
\doteq
\big(
p_2+\delta\xi_1
,
p_1-\delta\xi_1
\big),
$
$
\mathcal{I}_{5}^\xi
\doteq
\big(
p_1-\delta\xi_1
,
p_1
\big),
$
$
\mathcal{I}_{6}^\xi
\doteq
\big(
p_1
,
p_1+\delta\xi_1
\big),
$
$
\mathcal{I}_{7}^\xi
\doteq
\big(
p_2+\delta\xi_1
,
\infty
\big),
$
which are visualized  in Figure \ref{fig:s-neg-xi-axis}.
\end{minipage}
\hskip-0.3in
\begin{minipage}{0.4\linewidth}
\begin{center}
\begin{tikzpicture}[yscale=0.35, xscale=0.7]
%
%
\newcommand\X{0};
\newcommand\Y{0};
\newcommand\FX{8};
\newcommand\FY{8};
\newcommand\R{0.6};
\newcommand\Alpha{5};
\newcommand*{\TickSize}{2pt};
%
%
\draw[black,line width=1pt,-{Latex[black,length=2mm,width=2mm]}]
(-0.5,0)
--
(3,0)
node[above]
{\fontsize{\FX}{\FY}\bf \textcolor{black}{$\xi$}};

\draw[black,line width=1pt,-{Latex[black,length=2mm,width=2mm]}]
(0,-1.5)
--
(0,5)
node[right]
{\fontsize{\FX}{\FY}\bf \textcolor{black}{$\mu(\xi)$}};

\draw[line width=1pt, yscale=1,domain=-0.1:2.3,smooth,variable=\x,red]  
plot ({\x},{-(\x)^3+\Alpha+\Alpha*(\x-1)^3});

\draw[dashed,black,line width=0.5pt]
(0,0)
node[xshift=-0.2cm, yshift=-0.2cm]
{\fontsize{\FX}{\FY}\bf $0$}

({sqrt(\Alpha)/(sqrt(\Alpha)-1)},0)
node[xshift=0cm,yshift=-0.2cm]
{\fontsize{\FX}{\FY}\bf $p_1$}
--
({sqrt(\Alpha)/(sqrt(\Alpha)-1)},{(\Alpha*\Alpha-2*(sqrt(\Alpha))^(3))/(1-sqrt(\Alpha))^2})

({sqrt(\Alpha)/(sqrt(\Alpha)+1)},0)
node[yshift=-0.2cm]
{\fontsize{\FX}{\FY}\bf $p_2$}
--
({sqrt(\Alpha)/(sqrt(\Alpha)+1)},{(\Alpha*\Alpha+2*(sqrt(\Alpha))^(3))/(1+sqrt(\Alpha))^2})

%
%

({\Alpha/(\Alpha-1)},0)
node[xshift=0cm, yshift=-0.2cm]
{\fontsize{\FX}{\FY}\bf $q$}
--
({\Alpha/(\Alpha-1)},{-(\Alpha/(\Alpha-1))^3+\Alpha+\Alpha*(\Alpha/(\Alpha-1)-1)^3})
;

\end{tikzpicture}

\captionof{figure}{$\xi_1>0$}

\label{fig:s-neg-xi}
\end{center}
\end{minipage}

\vskip0.05in
\noindent
\begin{minipage}{1\linewidth}
\begin{center}
\begin{tikzpicture}[yscale=1, xscale=0.8]
%
%
\newcommand\X{0};
\newcommand\Y{0};
\newcommand\FX{11};
\newcommand\FY{11};
\newcommand\R{0.6};
\newcommand*{\TickSize}{2pt};
%
%
\draw[black,line width=1pt,-{Latex[black,length=2mm,width=2mm]}]
(-9,0)
--
(9,0)
node[above]
{\fontsize{\FX}{\FY} $\xi$};

\draw[red]

(-10,0)
node[yshift=0cm,black,line width=1pt]
{\bf \fontsize{\FX}{\FY}$\xi_1> 0:$}


(0,0)
node[yshift=-0.4cm]
{\fontsize{\FX}{\FY}$q$}
node[yshift=0cm]
{$\bullet$}
node[xshift=0cm,yshift=0.6cm, black]
{\fontsize{\FX}{\FY}$\mathcal{I}_4^\xi$}

(-7,0)
node[xshift=-1cm,yshift=0.6cm, black]
{\fontsize{\FX}{\FY}$\mathcal{I}_1^\xi$}

(-5,0)
node[yshift=-0.4cm]
{\fontsize{\FX}{\FY}$p_2$}
node[yshift=0cm]
{$\bullet$}
node[xshift=-1cm,yshift=0.6cm, black]
{\fontsize{\FX}{\FY}$\mathcal{I}_2^\xi$}
node[xshift=1cm,yshift=0.6cm, black]
{\fontsize{\FX}{\FY}$\mathcal{I}_3^\xi$}

(-2.5,0)
node[yshift=-0.4cm]
{\fontsize{\FX}{\FY}$p_2+\delta\xi_1$}
node[yshift=0cm]
{$\bullet$}

(-7.5,0)
node[yshift=-0.4cm]
{\fontsize{\FX}{\FY}$p_2-\delta\xi_1$}
node[yshift=0cm]
{$\bullet$}

(5,0)
node[yshift=-0.4cm]
{\fontsize{\FX}{\FY}$p_1$}
node[yshift=0cm]
{$\bullet$}
node[xshift=-1cm,yshift=0.6cm, black]
{\fontsize{\FX}{\FY}$\mathcal{I}_5^\xi$}
node[xshift=1cm,yshift=0.6cm, black]
{\fontsize{\FX}{\FY}$\mathcal{I}_6^\xi$}

(2.5,0)
node[yshift=-0.4cm]
{\fontsize{\FX}{\FY}$p_1-\delta\xi_1$}
node[yshift=0cm]
{$\bullet$}

(7.5,0)
node[yshift=-0.4cm]
{\fontsize{\FX}{\FY}$p_1+\delta\xi_1$}
node[yshift=0cm]
{$\bullet$}

(8,0)
node[xshift=0cm,yshift=0.6cm, black]
{\fontsize{\FX}{\FY}$\mathcal{I}_7^\xi$}

;

\draw [decorate,decoration={brace,mirror,raise=5pt},thick] (-2.5,0) --  (-5,0);

\draw [decorate,decoration={brace,mirror,raise=5pt},thick] (-5,0) --  (-7.5,0);

\draw [decorate,decoration={brace,mirror,raise=5pt},thick] (7.5,0) --  (5,0);

\draw [decorate,decoration={brace,mirror,raise=5pt},thick] (5,0) --  (2.5,0);

\end{tikzpicture}

\vskip-0.15in
\captionof{figure}{$\xi$-intervals}
\label{fig:s-neg-xi-axis}
\end{center}
\end{minipage}

\vskip0.05in
\nin
Now,
we see that our change of variables makes sense in each one of these intervals. Thus, making the change of variables $\mu=\mu(\xi)=d_{\alpha}(\xi,\xi_1)$ and defining $\mathcal{I}_k^\mu$ be the range of $\mu(\xi)$ for $\xi\in \mathcal{I}_k^\xi$,  we obtain
\begin{align}
\label{s-neg-theta2-est-no2}
\Theta_2
\lesssim
\sum\limits_{k=1}^7
J_k
\text{  with  }
J_k
\doteq
\frac{1}{(1+|\tau_1-\alpha\xi_1^{3}|)^{2b'}}
\int_{\mathcal{I}_{k}^{\mu}}
\frac{ 
\chi_{E_{II}}(\xi,\tau,\xi_1,\tau_1)
\, |\xi_1(\xi-\xi_1)|^{-2s}
\xi^2
(1+|\xi|)^{2s}
d\mu
}{(1+|\tau_1-\alpha\xi_1^3+\mu|)^{4b'-1}
\quad
|\mu'(\xi)|}.
\end{align}
{\bf Estimation near critical points.} 
Since $p_1, p_2$ are like 
$\xi_1$, and $\xi$ is near $p_1$ or $p_2$, we have 
$|\xi|\simeq |\xi_1|$ and $|\xi-\xi_1|\lesssim |\xi_1|$, 
which give us that  $|\xi_1(\xi-\xi_1)|^{-2s}
\xi^2
(1+|\xi|)^{2s}\lesssim |\xi_1|^{2-2s}$. Using this, we get
\begin{align*}
J_k(\xi_1,\tau_1)
\lesssim
\frac{|\xi_1|^{2-2s}}{(1+|\tau_1-\alpha\xi_1^{3}|)^{2b'}}
\int_{\mathcal{I}_{k}^{\mu}}
\frac{ 
\chi_{E_{II}}(\xi,\tau,\xi_1,\tau_1)
}{(1+|\tau_1-\alpha\xi_1^3+\mu|)^{4b'-1}}
\frac{1}{|\mu'(\xi)|}
d\mu,
\quad
k=2,3,5,6.
\end{align*}
\underline{Estimation for $J_2$ and $J_3$.} Like  estimate \eqref{s-0-mu-der-est}, 
for $\xi\in I_2^{\xi}\cup I_3^{\xi}$,
we get 
$
|\mu'(\xi)|
\gtrsim
|\xi_1|^{1/2} \sqrt{
|\mu(\xi)-\mu(p_2)|
}.
$
Using 
this estimate and making the change of variables $\mu_1=\tau_1-\alpha\xi_1^3+\mu$, we obtain
\begin{align}
\label{s-neg-theta2-est-no3}
J_k(\xi_1,\tau_1)
\lesssim
\frac{|\xi_1|^{3/2-2s}}{(1+|\tau_1-\alpha\xi_1^{3}|)^{2b'}}
\int_{\rr}
\frac{ 
\chi_{E_{II}}(\xi,\tau,\xi_1,\tau_1)
}{(1+|\mu_1|)^{4b'-1}}
\frac{1}{\sqrt{
|\mu_1-(\tau_1-\alpha\xi_1^3)-\mu(p_2)|
}}
d\mu_1,
\quad
k=2,3.
\end{align}
Furthermore,  using comparison property \eqref{d-comparison},  we have $|\mu_1|\lesssim |\tau_1-\alpha\xi_1^3|$. For the $\mu_1$-integral in \eqref{s-neg-theta2-est-no3},
applying estimate \eqref{eq:calc_4} with $2(1-\ell)=4b'-1$, $a=\tau_1-\alpha\xi_1^3+\mu(p_2)$ and $c=|\tau_1-\alpha\xi_1^3|$, 
we get
\begin{align}
\label{s-neg-theta2-est-no4}
J_k(\xi_1,\tau_1)
\lesssim
\frac{|\xi_1|^{3/2-2s}}{(1+|\tau_1-\alpha\xi_1^{3}|)^{2b'}}
\frac{(1+|\tau_1-\alpha\xi_1^{3}|)^{2-4b'}}{(1+|\tau_1-\alpha\xi_1^3+\mu(p_2)|)^{\frac12}}
\lesssim
\frac{|\xi_1|^{3/2-2s}}{
(1+|\tau_1-\alpha\xi_1^{3}|)^{6b'-2}
},
\quad
k=2,3.
\end{align}
Moreover, using $|\tau_1-\alpha\xi_1^3|\ge \frac13|d_\alpha|$, estimate \eqref{d-bound-xi1} and $|\xi|\simeq |\xi_1|$ we get 
$
|\tau_1-\alpha\xi_1^3|
\gtrsim
|\xi|
\xi_1^2
\simeq
|\xi_1|^3.
$
Combining this estimate with \eqref{s-neg-theta2-est-no4},  we get 
$
J_k(\xi_1,\tau_1)
\lesssim
|\xi_1|^{15/2-2s-18b'}
$
($k=2,3$),
which is bounded if 
$
15/2-2s-18b'
\le 
0
$
or
$b'\ge \frac{15-4s}{36}$.
For $b'<\frac12$, it suffices to have 
$
s>-\frac34.
$
This completes the estimation near critical points.

\vskip0.05in
\nin
\underline{Estimation for $J_5$ and $J_6$.} Since this estimate is like that for $J_2$ and $J_3$, we omit it here.

\vskip0.05in
\noindent
{\bf Estimation away from critical points.} 
Using property \eqref{d-comparison}, we get
 $|\tau_1-\alpha\xi_1^3+\mu|\lesssim |\tau_1-\alpha\xi_1^3|$, 
which implies that
 $
 |\tau_1-\alpha\xi_1^3+\mu|/|\tau_1-\alpha\xi_1^3|
\lesssim
1.
$
Moving $\frac{1}{(1+|\tau_1-\alpha\xi_1^{3}|)^{2b'}}$ inside $\tau$-integral and also dividing the integrand by  
$ 
(|\tau_1-\alpha\xi_1^3+\mu|/|\tau_1-\alpha\xi_1^3|)^{\varepsilon}
$
for $\varepsilon>0$, to be chosen later,  from \eqref{s-neg-theta2-est-no2} we get 
\begin{align}
\label{s-neg-case2-theta-2-fin-1}
J_k(\xi_1,\tau_1)
\lesssim&
\int_{\mathcal{I}_k^\mu}
\frac{1
}{(1+|\tau_1-\alpha\xi_1^3|)^{2b'-\varepsilon}}
\frac{ 
\chi_{E_{II}}(\xi,\tau,\xi_1,\tau_1)
\, |\xi_1(\xi-\xi_1)|^{-2s}
\xi^2
(1+|\xi|)^{2s}
d\mu
}{(1+|\tau_1-\alpha\xi_1^3+\mu|)^{4b'-1+\varepsilon}
\quad
|\mu'(\xi)|},
\quad
k=1,4,7.
\end{align}
Using comparison property \eqref{d-comparison} and adding estimates \eqref{d-bound-xi}, \eqref{d-bound-xi1} up we get $|\tau_1-\alpha\xi_1^3|\gtrsim |\xi|(\xi^2+\xi_1^2)$.  Also, when $\xi$ is away from the critical points $p_1$ and $p_2$, by using formula \eqref{d-partial-xi}
and by
using elementary inequality of Lemma \ref{dist-lem},
we have a lower bound for $|\mu'(\xi)|$, that is $|\mu'(\xi)|\gtrsim \xi^2+\xi_1^2$.  Furthermore, for $s<0$, we have $(1+|\xi|)^{2s}\lesssim1$.
Combining these estimates with \eqref{s-neg-case2-theta-2-fin-1}, we get 
\begin{align}
\label{s-neg-case2-theta-2-fin-2}
J_k(\xi_1,\tau_1)
\lesssim
\int_{\mathcal{I}_k^\mu}
\frac{|\xi_1(\xi-\xi_1)|^{-2s}|\xi|^{2+\varepsilon-2b'}}{(\xi^2+\xi_1^2)^{2b'-\varepsilon+1}}
\frac{ 
\chi_{E_{II}}(\xi,\tau,\xi_1,\tau_1)
}{(1+|\tau_1-\alpha\xi_1^3+\mu|)^{4b'-1+\varepsilon}}
d\mu,
\quad
k=1,4,7.
\end{align}
Finally, combining 
 the elementary facts $|\xi|\le(\xi^2+\xi_1^2)^{1/2}$, $|\xi(\xi-\xi_1)|\lesssim (\xi^2+\xi_1^2)$ with  \eqref{s-neg-case2-theta-2-fin-2}, we get
\begin{align*}
J_k(\xi_1,\tau_1)
\lesssim
\int_{\rr}
\frac{1}{(\xi^2+\xi_1^2)^{2b'-\varepsilon+1-(1+\varepsilon/2-b')+2s}}
\frac{ 
\chi_{E_{II}}(\xi,\tau,\xi_1,\tau_1)
}{(1+|\tau_1-\alpha\xi_1^3+\mu|)^{4b'-1+\varepsilon}}
d\mu,
\quad
k=1,4,7,
\end{align*}
which is bounded if $2b'-\varepsilon+1-(1+\varepsilon/2-b')+2s\ge 0$ and $4b'-1+\varepsilon\ge 1$.
Combining these two estimates, we get 
$
b'
\ge
(3-2s)/9.
$
For $b'<1/2$, it suffices to have
$
s>-3/4.
$
This completes the estimate for all $J_k$ and  the proof of Lemma \ref{s-neg-lem2}.
\,\,
$\square$

%
%
%
%
%
%
%
%
%
%
%
%
\subsection{$s\ge \frac34$ with $\alpha=4$}
%
%
First, we observe that to prove the  bilinear estimate \eqref{bi-est-X-1} with $\alpha=4$, i.e.
$
\|
\p_x(fg)
\|_{X^{s,-b,\theta-1}_{1}}
\le
c_2
\|
f
\|_{X^{s,b',\theta'}_{4}}
\|
g
\|_{X^{s,b',\theta'}_{4}},
$
it suffices to prove its $L^2$ formulation \eqref{bilinear-est-L2-form-no1}.
By symmetry (in convolution writing), we  assume that 
$
|\xi-\xi_1|
\leq
|\xi_1|.
$
Also, like the case $\alpha\neq 4$, we decompose $Q'(\xi,\tau,\xi_1,\tau_1)$ as the sum of three functions
$
Q'(\xi,\tau,\xi_1,\tau_1)
=
Q_1(\xi,\tau,\xi_1,\tau_1)
+
Q_2(\xi,\tau,\xi_1,\tau_1)
+
Q_3(\xi,\tau,\xi_1,\tau_1),
$
where
$
Q_1(\xi,\tau,\xi_1,\tau_1)
\doteq
\chi_{|\xi_1|< 1}
\cdot
\chi_{|\xi-\xi_1|< 1}
\cdot
Q'(\xi,\tau,\xi_1,\tau_1),
$
$
Q_2(\xi,\tau,\xi_1,\tau_1)
\doteq
\chi_{|\xi_1|\ge 1}
\cdot
\chi_{|\xi-\xi_1|<1}
\cdot
Q'(\xi,\tau,\xi_1,\tau_1)
$
and
\begin{align}
\label{u4-def-Q3}
&Q_3(\xi,\tau,\xi_1,\tau_1)
\doteq
\chi_{|\xi_1|\ge 1}
\cdot
\chi_{|\xi-\xi_1|\ge 1}
\cdot
Q'(\xi,\tau,\xi_1,\tau_1)
\nn
\\
=&
\frac{|\xi|}{(1+|\tau-\xi^{3}|)^{b'}}
\frac{(1+|\xi|)^s}{
(1+|\xi_1|)^{s}
(1+|\xi-\xi_1|)^{s}
}
\frac{
\chi_{|\xi_1|\ge 1}
}{
(1+|\tau_1-4\xi_1^3|)^{b'}
}
\frac{
\chi_{|\xi-\xi_1|\ge 1}
}
{
(1+|\tau-\tau_1-4(\xi-\xi_1)^3|)^{b'}
}.
\end{align}
Now, to prove the bilinear estimate \eqref{bilinear-est-L2-form-no1} with $\alpha=4$,  it suffices to show that for $\ell
=
1,2,3$ we have
\begin{align}
\label{u4-bilinear-est-l}
\Big\| 
\int_{\rr^2}
Q_\ell(\xi,\tau,\xi_1,\tau_1)
c_{f,4}(\xi-\xi_1,\tau-\tau_1) c_{g,4}(\xi_1,\tau_1) 
d\xi_1 d\tau_1 \Big\|_{L^2_{\xi,\tau}}
\lesssim
\big\|c_{f,4}\big\|_{L^2_{\xi,\tau}} \big\|c_{g,4}\big\|_{L^2_{\xi,\tau}}.
\end{align}
\vskip0.05in
\noindent
\underline{\bf Proof for  multipliers $Q_1$  and $Q_2$.} For $s\ge 0$, we have 
$
\frac{(1+|\xi|)^s}{
(1+|\xi_1|)^{s}
(1+|\xi-\xi_1|)^{s}
}
\lesssim
1.
$
Using this estimate, $Q_1$ and $Q_2$ are reduced to the $Q_1$ and $Q_2$ defined by \eqref{u-s-0-Q1} and \eqref{u-s-0-Q2} respectively. Also, we notice  that 
 Lemmas \ref{s-0-theta-1-lemma}--\ref{s-0-theta-3-lemma} hold for $\alpha=4$. Therefore,  we get the desired estimate \eqref{u4-bilinear-est-l} for $Q_1$ and $Q_2$.

\vskip.1in
\noindent
\underline{\bf  Proof for  multiplier  $Q_3$.}
By symmetry (in convolution writing), we  assume that 
$
|\tau-\tau_1-4(\xi-\xi_1)^{3}|
\leq
|\tau_1-4\xi_1^{3}|.
$
Also, 
we consider the following two microlocalizations:
\vskip0.05in
\noindent
{\bf $\bullet$ Microlocalization I: $|\tau_1-4\xi_1^{3}|\leq|\tau-\xi^{3}|$}. In this situation we define the domain $B_I$ to be
\begin{align*}
B_I
\doteq
\big\{
(\xi,\tau,\xi_1,\tau_1) \in{\rr}^4: |\tau-\tau_1-4(\xi-\xi)^{3}|\leq|\tau_1-4\xi_1^{3}|\leq|\tau-\xi^{3}|,
|\xi_1|\ge1,
|\xi-\xi_1|\ge1
\big\}.
\end{align*}
\vskip0.05in
\noindent
{\bf $\bullet$ Microlocalization II: $|\tau-\xi^{3}|\leq|\tau_1-4\xi_1^{3}|$}. 
In this case we define the domain $B_{II}$ to be
\begin{align*}
\hskip-0.1in
B_{II}
\doteq
\big\{(\xi,\tau,\xi_1,\tau_1)\in{\rr}^4: 
\max\{|\tau-\tau_1-4(\xi-\xi)^{3}|,
|\tau-\xi^{3}|
\}\leq|\tau_1-4\xi_1^{3}|,
|\xi_1|
\hskip-0.01in
\ge1,
|\xi-\xi_1|
\hskip-0.01in
\ge1
\big\}.
\end{align*}
{\bf Proof of bilinear estimate in microlocalization I.} In this situation, using the Cauchy--Schwarz inequality with respect to 
$(\xi_1,\tau_1)$, and taking the supremum over  $(\xi,\tau)$ we
arrive at
\begin{align*}
&  \Big\| 
 \int_{\rr^2} 
 (\chi_{B_I} Q_3)(\xi,\tau,\xi_1,\tau_1)
 c_{f,4}(\xi-\xi_1,\tau-\tau_1) c_{g,4}(\xi_1,\tau_1) 
 d\xi_1 d\tau_1 \Big\|_{L^2_{\xi,\tau}}
 \nn
 \\
 &
\lesssim
\big\|c_{f,4}\big\|_{L^2_{\xi,\tau}} \big\|c_{g,4}\big\|_{L^2_{\xi,\tau}}
\Big\|
\int_{\rr^2}
 (\chi_{B_{I}} Q_3^2)(\xi,\tau,\xi_1,\tau_1)
d\xi_1 d\tau_1
\Big\|_{L^{\infty}_{\xi,\tau}}^{1/2}.
\end{align*}
Now, to prove bilinear estimate \eqref{u4-bilinear-est-l} for $Q_3$,  it suffices to show 
 the following result for $Q_3$ as in \eqref{u4-def-Q3}.
\begin{lemma} 
\label{u4-bi-lem-1}
If $\frac{1}{3}< b'<1/2$ and $s\ge \frac34$, then  for $\xi,\tau\in\rr$
\begin{align}
\label{u4-theta1}
\hskip-0.05in
\Theta_1(\xi,\tau)
\hskip-0.05in
\doteq&
\frac{\xi^2}{(1+|\tau-\xi^{3}|)^{2b'}}
\int_{\rr^2}
\hskip-0.05in
\frac{(1+|\xi|)^{2s}}{
(1+|\xi_1|)^{2s}
(1+|\xi-\xi_1|)^{2s}
}
\frac{ \chi_{B_I}(\xi,\tau,\xi_1,\tau_1) \quad d\tau_1 d\xi_1 
}{(1+|\tau-\tau_1-4(\xi-\xi_1)^{3}|)^{2b'}(1+|\tau_1-4\xi_1^{3}|)^{2b'}}
\nn
\\
\lesssim&
 1.
\end{align}
\end{lemma}
\vskip0.05in
\noindent
{\bf Proof of bilinear estimate in microlocalization II.} In this case,  using duality and applying the Cauchy--Schwarz inequality  first in $(\xi_1, \tau_1)$ and then in $(\xi, \tau)$, we get
\begin{align*}
 & \Big\| 
 \int_{\rr^2} 
 (\chi_{B_{II}} Q_3)(\xi,\tau,\xi_1,\tau_1)
 c_{f,4}(\xi-\xi_1,\tau-\tau_1) c_{g,4}(\xi_1,\tau_1) 
 d\xi_1 d\tau_1 \Big\|_{L^2_{\xi,\tau}}
 \\
 \lesssim&
\|c_{f,4}\|_{L^2_{\xi,\tau}}\|c_{g,4}\|_{L^2_{\xi,\tau}} \Big\| 
\int_{\rr^2} (\chi_{B_{II}} Q_3^2)(\xi,\tau,\xi_1,\tau_1)d\xi d\tau 
\Big\|_{L^{\infty}_{\xi_1,\tau_1}}^{1/2}.
\end{align*}
Now,  to prove bilinear estimate \eqref{u4-bilinear-est-l} for $Q_3$, it suffices to prove next result
with $Q_3$ as in \eqref{u4-def-Q3}.
\begin{lemma} 
\label{u4-bi-lem-2}
If $\frac{1}{3}< b'<1/2$ and $s\ge \frac34$, then  for $\xi_1,\tau_1\in\rr$
\begin{align}
\label{u4-theta2}
\Theta_2(\xi_1,\tau_1)
\doteq&
\frac{1}{(1+|\tau_1-4\xi_1^{3}|)^{2b'}}
\int_{\rr^2}
\frac{\xi^2(1+|\xi|)^{2s}}{
(1+|\xi_1|)^{2s}
(1+|\xi-\xi_1|)^{2s}
}
\frac{ 
\chi_{B_{II}}(\xi,\tau,\xi_1,\tau_1) 
\quad d\tau d\xi 
}{(1+|\tau-\tau_1-4(\xi-\xi_1)^{3}|)^{2b'}(1+|\tau-\xi^{3}|)^{2b'}}
\nn
\\
\lesssim& 1.
\end{align}
\end{lemma}
\noindent
{\bf Proof of Lemma \ref{u4-bi-lem-1}.} For the $\tau_1$-integral in \eqref{u4-theta1}, applying estimate \eqref{eq:calc_5}  with  $\ell=\ell'=b'>\frac14$, $a=\tau-4(\xi-\xi_1)^3$ and $c=4\xi_1^3$, we get
\begin{align}
\label{u4-theta1-est}
\Theta_1(\xi,\tau)
\lesssim
\frac{\xi^2}{(1+|\tau-\xi^{3}|)^{2b'}}
\int_{\rr}
\frac{(1+|\xi|)^{2s}}{
(1+|\xi_1|)^{2s}
(1+|\xi-\xi_1|)^{2s}
}
\frac{ 
\chi_{B_I}(\xi,\tau,\xi_1,\tau_1)
\,
 d\xi_1 
}{(1+|\tau-\xi^3-d_{4}(\xi,\xi_1)|)^{4b'-1}},
\end{align}
where we recall that  
$
d_{4}(\xi,\xi_1)
=
-\xi^3+4\xi_1^3+4(\xi-\xi_1)^3
=
12\xi(\xi_1-\xi/2)^2.
$
Furthermore, in this microlocalization, using property \eqref{d-comparison} we get $|\tau-\xi^3|\ge\frac13|d_{4}(\xi,\xi_1)|$, which implies that 
$
|\tau-\xi^3-d_{4}(\xi,\xi_1)|
\lesssim 
|\tau-\xi^3|.
$ 
This gives us
$
|\tau-\xi^3-d_{4}(\xi,\xi_1)|/
|\tau-\xi^3|
\lesssim 1.
$ 
Moving $\frac{1}{(1+|\tau-\xi^{3}|)^{2b'
}}$ inside $\xi_1$-integral in \eqref{u4-theta1-est} and also dividing the integrand by  
$ 
(|\tau-\xi^3-d_4(\xi,\xi_1)|/|\tau-\xi^3|)^{2b'}
$,  we have 
\begin{align}
\label{u4-theta1-est-no1}
\Theta_1(\xi,\tau)
\lesssim
\xi^2
\int_{\rr}
\frac{(1+|\xi|)^{2s}}{
(1+|\xi_1|)^{2s}
(1+|\xi-\xi_1|)^{2s}
}
\frac{ 
\chi_{B_I}
(\xi,\tau,\xi_1,\tau_1)
 d\xi_1 
}{(1+|\tau-\xi^3-d_{4}(\xi,\xi_1)|)^{6b'-1}}.
\end{align}
\begin{minipage}{0.65\linewidth}
Now, our strategy for estimating the above $\xi_1$-integral is to make the change of variables,
$
\mu=\mu(\xi_1)=d_{4}(\xi,\xi_1).
$
Using  formula \eqref{d-partial-xi1} with $\alpha=4$, i.e.
$
\p d_{4}/\p\xi_1
=
24
\xi(\xi_1-\xi/2),
$
we find the critical point $\xi/2$ and
we need to split  the $\xi_1$-integral  at this point, so that this change of variables makes sense.
In the case that   $\xi>0$,
the graph of $\mu=\mu(\xi_1)=d_{4}(\xi,\xi_1)$ looks as in Figure \ref{fig:a-4-xi1}. 
\end{minipage}
\hskip-0.2in
\begin{minipage}{0.45\linewidth}
\begin{center}
\begin{tikzpicture}[yscale=0.5, xscale=0.8]
%
%
\newcommand\X{0};
\newcommand\Y{0};
\newcommand\FX{8};
\newcommand\FY{8};
\newcommand\R{0.6};
\newcommand*{\TickSize}{2pt};
%
%
\draw[line width=1pt,-{Latex[length=2mm]}]
(-0.3,0)
--
(3,0)
node[above]
{\fontsize{\FX}{\FY} $\xi_1$}
;

\draw[line width=1pt,-{Latex[length=2mm]}]
(0,-0.2)
--
(0,2.9)
node[xshift=0.5cm]
{\fontsize{\FX}{\FY} $\mu(\xi_1)$};

\draw[red, line width=1pt, domain=-0.5:2.5,variable=\x]
plot
({\x},{(\x-1)^2})
;


\draw[black,dashed]
(1,0)
node[yshift=-0.3cm]
{\fontsize{\FX}{\FY}$ \xi/2$};

\draw[black]
(0,0)
node[xshift=-0.2cm,yshift=-0.2cm]
{\fontsize{\FX}{\FY} $0$}

;

\end{tikzpicture}

\vskip-0.15in
\captionof{figure}{$\xi>0$}
\label{fig:a-4-xi1}

\end{center}
\end{minipage}

\vskip0.05in
\nin
Now, using the  $\xi_1$-intervals 
$
\mathcal{I}_1^{\xi_1}
\doteq
(-\infty, \frac14\xi),
$
$
\mathcal{I}_2^{\xi_1}
\doteq
(\frac14\xi, \frac12\xi),
$
$
\mathcal{I}_3^{\xi_1}
\doteq
(\frac12\xi, \frac34\xi)
$
and
$
\mathcal{I}_4^{\xi_1}
\doteq
(\frac34\xi, \infty),
$
we obtain
\begin{align}
\label{a-4-theta1-est-deco}
\Theta_1
\lesssim
\sum\limits_{k=1}^4
J_k
\text{ with }
J_k(\xi,\tau)
\doteq
\xi^2
\int_{I_k^{\xi_1}}
\frac{(1+|\xi|)^{2s}}{
(1+|\xi_1|)^{2s}
(1+|\xi-\xi_1|)^{2s}
}
\frac{ 
\chi_{B_I}
(\xi,\tau,\xi_1,\tau_1)
 d\xi_1 
}{(1+|\tau-\xi^3-d_{4}(\xi,\xi_1)|)^{6b'-1}}.
\end{align}
{\bf Estimation away from critical points.} On each one of the intervals $\mathcal{I}_1^{\xi_1}$ and $\mathcal{I}_4^{\xi_1}$, making the change of variables $\mu=\mu(\xi_1)=\tau-\xi^3-d_{4}(\xi,\xi_1)$ and using formula \eqref{d-partial-xi1} with $\alpha=4$,
 we get 
\begin{align}
\label{u4-theta1-est-no2}
J_k(\xi,\tau)
\lesssim
\xi^2
\int_{\mathcal{I}_k^\mu}
\frac{(1+|\xi|)^{2s}}{
(1+|\xi_1|)^{2s}
(1+|\xi-\xi_1|)^{2s}
}
\frac{ 
\chi_{B_I}
(\xi,\tau,\xi_1,\tau_1)
}{(1+|\mu|)^{6b'-1}}
\frac{1}{|\xi||\xi_1-\xi/2|}
d\mu,
\quad
k=1,4.
\end{align}
For $\xi_1\le \frac14\xi$ or $\xi_1\ge 3/4\xi$, we get $|\xi_1-\xi/2|\ge \frac14|\xi|\gtrsim |\xi|$. Also, using the triangle inequality, we get $1+|\xi|\le (1+|\xi_1|)(1+|\xi-\xi_1|)$, which implies that $\frac{(1+|\xi|)^{2s}}{
(1+|\xi_1|)^{2s}
(1+|\xi-\xi_1|)^{2s}
}\lesssim 1$ for $s\ge 0$. Combining  these estimates with \eqref{u4-theta1-est-no2},  we obtain
$
J_k(\xi,\tau)
\lesssim
\int_{\mathcal{I}_k^{\mu}}
\frac{ 
\chi_{B_I}
(\xi,\tau,\xi_1,\tau_1)
}{(1+|\mu|)^{6b'-1}}
d\mu
$
($k=1,4$),
which is bounded if  $6b'-1>1$ or $b'>1/3$. This completes the estimate for $J_1$ and $J_4$.

\vskip0.05in
\noindent
{\bf Estimation near critical points.}
Since the estimate for $J_3$ is same to the estimate of $J_2$, here we only provide the estimate for $J_2$.
For $1/4\xi\le \xi_1\le 1/2\xi$, we have $1\le |\xi_1|\simeq |\xi|\simeq |\xi-\xi_1|$, which implies that  $\frac{(1+|\xi|)^{2s}}{
(1+|\xi_1|)^{2s}
(1+|\xi-\xi_1|)^{2s}
}\lesssim |\xi|^{-2s}$.
Using this estimate, from \eqref{a-4-theta1-est-deco} we obtain
$
J_2(\xi,\tau)
\lesssim
|\xi|^{2-2s}
\int_{1/4\xi}^{1/2\xi}
\frac{ 
1
}{(1+|\tau-\xi^3-d_{4}(\xi,\xi_1)|)^{6b'-1}} d\xi_1.
$
In addition, we split $J_2$ as
$
J_2
=
J_{21}
+
J_{22},
$
with
$$
J_{2k}
\doteq
|\xi|^{2-2s}
\int_{\mathcal{I}_{2k}^{\xi_1}}
\frac{ 
1
}{(1+|\tau-\xi^3-d_{4}(\xi,\xi_1)|)^{6b'-1}}
 d\xi_1,
 \quad
 k=1,2,
$$
where $\mathcal{I}_{21}^{\xi_1}=(\frac14\xi,\frac12\xi -\frac{1}{\sqrt{|\xi|}})$ and $\mathcal{I}_{22}^{\xi_1}=(\frac12\xi -\frac{1}{\sqrt{|\xi|}}, \frac12\xi)$ .

\vskip0.05in
\nin
\underline{Estimate for $J_{21}$.}  Making the change of variables $\mu=\mu(\xi_1)=\tau-\xi^3-d_{4}(\xi,\xi_1)$, we get
\begin{align}
\label{u4-theta1-B2-est}
&J_{21}
\lesssim
|\xi|^{2-2s}
\int_{\mathcal{I}_{21}^{\mu}}
\frac{ 
1 
}{(1+|\tau-\xi^3-\mu|)^{6b'-1}}
\frac{1}{|\mu'|}
d\mu.
\end{align}
Now,
using formula \eqref{d-partial-xi1} with $\alpha=4$,  we get the lower bound for $|\mu'(\xi_1)|$, that is  $|\mu'(\xi_1)|=|24\xi(\xi_1-\xi/2)|\gtrsim |\xi| \frac{1}{\sqrt{|\xi|}}=\sqrt{|\xi|}$. Combining this estimate with \eqref{u4-theta1-B2-est}, for $6b'-1>1$ or $b'>1/3$ we get 
\begin{align}
\label{u4-theta1-J21-est}
J_{21}(\xi,\tau)
\lesssim
|\xi|^{3/2-2s}
\int_{\mathcal{I}_{21}^{\mu}}
\frac{ 
1 
}{(1+|\tau-\xi^3-\mu|)^{6b'-1}}
d\mu
\lesssim
|\xi|^{3/2-2s}.
\end{align}
\underline{Estimate for $J_{22}$.}
For $6b'-1\ge 0$ or $b'\ge \frac16$ we get 
\begin{align}
\label{u4-theta1-J22-est}
J_{22}(\xi,\tau)
\lesssim
|\xi|^{2-2s}
\int_{\frac12\xi -\frac{1}{\sqrt{|\xi|}}}^{\frac12\xi}
\frac{1}{1}
d\xi_1
\lesssim
|\xi|^{2-2s}
\frac{1}{\sqrt{|\xi|}}
=
|\xi|^{3/2-2s}.
\end{align}
Combining estimates \eqref{u4-theta1-J21-est} with \eqref{u4-theta1-J22-est}, we get
$
J_2(\xi,\tau)
\lesssim
|\xi|^{3/2-2s}.
$
Since $|\xi|\gtrsim 1$, $J_{2}$ is bounded if $3/2-2s\le 0$ or $s\ge \frac34$.
This completes the proof of Lemma \ref{u4-bi-lem-1}.\,\,
$\square$

\vskip0.05in
\noindent
{\bf Proof of Lemma \ref{u4-bi-lem-2}.} 
For the $\tau$-integral in \eqref{u4-theta2},
applying estimate \eqref{eq:calc_5}  with $\ell=\ell'=b'>\frac14$, $a=\tau_1+4(\xi-\xi_1)^3$ and $c=\xi^3$, we get
\begin{align}
\label{u4-theta2-est}
\Theta_2(\xi_1,\tau_1)
\lesssim
\frac{1}{(1+|\tau_1-4\xi_1^{3}|)^{2b'}}
\int_{\rr}
\frac{\xi^2(1+|\xi|)^{2s}}{
(1+|\xi_1|)^{2s}
(1+|\xi-\xi_1|)^{2s}
}
\frac{ 
\chi_{B_{II}}(\xi,\tau,\xi_1,\tau_1) 
 d\xi 
}{(1+|\tau_1-4\xi_1^3+d_{4}(\xi,\xi_1)|)^{4b'-1}}
d\xi,
\end{align}
where $d_{4}=-\xi^3+4\xi_1^3+4(\xi-\xi_1)^3$.
Like  \eqref{u4-theta1-est-no1}, using $|\tau_1-4\xi_1^{3}|\gtrsim |\tau_1-4\xi_1^3+d_{4}(\xi,\xi_1)|$, we obtain
\begin{align}
\label{u4-theta2-est-no1}
\Theta_2(\xi_1,\tau_1)
\lesssim
\int_{\rr}
\frac{\xi^2(1+|\xi|)^{2s}}{
(1+|\xi_1|)^{2s}
(1+|\xi-\xi_1|)^{2s}
}
\frac{ 
\chi_{B_{II}}
}{(1+|\tau_1-4\xi_1^3+d_{4}(\xi,\xi_1)|)^{6b'-1}}
d\xi.
\end{align}
Next, we  estimate  the $\xi$-integral for $\xi$ near the critical point of $d_{4}(\xi,\xi_1)$ and for $\xi$ away from the critical point of $d_{4}(\xi,\xi_1)$. Since these estimates are similar to the ones that we present in Lemma  \ref{u4-bi-lem-1}, we omit them here completing the proof for Lemma \ref{u4-bi-lem-2}.
\,\,
$\square$

%
%
%
%
%
%
%
%
%
%
%
%
\subsection{Proof of $v$-equation Bilinear Estimates 
\eqref{bi-est-X-2}
}
\setcounter{equation}{0}
In its $L^2$ formulation  this  estimate reads
\begin{align}
\label{v-bilinear-est-L2-form}
 \Big\| 
\int_{\rr^2}
 Q(\xi,\tau,\xi_1,\tau_1)
 c_{f,\alpha}(\xi-\xi_1,\tau-\tau_1) c_{g,1}(\xi_1,\tau_1) 
 d\xi_1 d\tau_1 \Big\|_{L^2_{\xi,\tau}}
\lesssim
\big\|c_{f,\alpha}\big\|_{L^2_{\xi,\tau}} \big\|c_{g,1}\big\|_{L^2_{\xi,\tau}},
\end{align}
where 
$Q$ is the multiplier defined by
\begin{align*}
Q(\xi,\tau,\xi_1,\tau_1)
\doteq&
|\xi|
\cdot
\Big[
\frac{
(1+|\xi|)^s
}{
(1+|\tau-\alpha\xi^{3}|)^{b}
}
+
\frac{
\chi_{|\xi|< 1}
}
{
(1+|\tau|)^{1-\theta}
}
\Big]
\cdot
\frac{1}{
(1+|\xi_1|)^{s}
(1+|\tau_1-\xi_1^3|)^{b'}
+
\chi_{|\xi_1|< 1}(1+|\tau_1|)^{\theta'}
}
\\
\times&
\frac{1}{
(1+|\xi-\xi_1|)^{s}
(1+|\tau-\tau_1-\alpha(\xi-\xi_1)^3|)^{b'}
+
\chi_{|\xi-\xi_1|< 1}(1+|\tau-\tau_1|)^{\theta'}
}.
\end{align*}
The corresponding Bourgain quantity is 
\begin{align}
\label{v-d3-quantity}
\widetilde{d}_{\alpha}(\xi,\xi_1)
\doteq
(\tau-\alpha\xi^3)-(\tau_1-\xi_1^3)-[\tau-\tau_1-\alpha(\xi-\xi_1)^3]
=
-\alpha\xi^3+\xi_1^3+\alpha(\xi-\xi_1)^3.
\end{align}
Since $\widetilde{d}_{\alpha}(\xi,\xi_1)=-d_{\alpha}(\xi_1,\xi)$,
where $d_{\alpha}(\xi,\xi_1)$ is the Bourgain quantity for $u$-equation, the proof for the bilinear estimate \eqref{bi-est-X-2} is similar to the proof for the bilinear estimate \eqref{bi-est-X-1}. In fact, we have similar microlocalizations, which follows from the new comparison property:  
$
\max\{
|\tau-\alpha\xi^3|, |\tau_1-\xi_1^3|, |\tau-\tau_1-\alpha(\xi-\xi_1)^3|
\}
\ge
\frac13
|\widetilde{d}_{\alpha}(\xi,\xi_1)|.
$
The main difference between the proof of bilinear estimate \eqref{bi-est-X-1} and the proof of bilinear estimate \eqref{bi-est-X-2} is that
the Bourgain quantity $\widetilde{d}_{\alpha}(\xi,\xi_1)$
 is a quadratic function of $\xi$ and a cubic function of $\xi_1$.

%
%
%
%
%
%
%
%
%
%
%
%
\section{ Proof of   Bilinear Estimates in $Y^{s,b}$}
\label{sec:Y-bilinear-est}
\setcounter{equation}{0}
The temporal bilinear estimate  \eqref{bi-est-Y-1} for the nonlinearity of the
 $u$-equation can be reduced to the 
 corresponding estimates for the KdV equation,
 whose proof can be found  in \cite{h2006}.
 In fact, for $s\ge 0$, in the $L^2$ formulation, if $|\tau|\lesssim |\xi|^3$, then we have $\|\p_x(f,g)\|_{Y^{s,-b}_1}\lesssim \|\p_x(f,g)\|_{X^{s,-b}_1}$, which will be reduced to the bilinear estimate in Bourgain spaces. Also, for $|\tau|\gg|\xi|^3$, we get $|\tau-\xi^3|\simeq |\tau-\alpha\xi^3|$, which reduces  
estimate  \eqref{bi-est-Y-1} to the corresponding estimates for the KdV equation. For $s<0$, working similarly we can reduce it to the corresponding estimates for the KdV equation.

Now,
we do the proof of the temporal bilinear estimate for $v$, i.e.  estimate \eqref{bi-est-Y-2}.
Here we consider the case  $\max\{0, s_{c}(\alpha)\}\le s<3$
and $\alpha\ne 1$, which is most relevant for the Majda-Biello system.
The case  $-3/4<s<0$ and $\alpha>4$  is similar for 
$\alpha=1$ we refer to \cite{h2006}.
We begin by  defining a domain $A$ in which $\tau\gg |\xi|^3$,
that is
\begin{equation}
\label{Y-bi-domain-A}
A
\doteq
\Big\{
(\xi,\tau)\in\rr^2
:
|\tau|>
(10+10\alpha)^7\Big[2\big|\frac{\sqrt{\alpha}}{\sqrt{\alpha}-1}\big|
+
2
\big|
\frac{3\alpha}{\alpha-1}
\big|
+
2\Big]^3|\xi|^3
\Big\}.
\end{equation}
Then, for  $s\ge 0$ we have
\begin{align}
\label{v-Y-est-1}
\|
w_{fg}
\|_{Y^{s,-b}_{\alpha}}^2
\lesssim
\int_{\rr^2}
\chi_{A}
(1+|\tau|)^{\frac{2s}{3}}
(1+|\tau-\alpha\xi^3|)^{-2b}
|\widehat{w}_{fg}(\xi,\tau)|^2
d\xi d\tau
+
\|w_{fg}\|_{X^{s,-b}_\alpha}^2.
\end{align}
Thus, to prove bilinear estimate \eqref{bi-est-Y-2}, it suffices to prove that
\begin{align}
\label{v-bi-est-Y-reduced}
\Big(
\int_{\rr^2}
\chi_{A}
(1+|\tau|)^{\frac{2s}{3}}
(1+|\tau-\alpha\xi^3|)^{-2b}
|\widehat{w}_{fg}(\xi,\tau)|^2
d\xi d\tau
\Big)^{1/2}
\lesssim
\|f\|_{X^{s,b'}_{\alpha}}
\|g\|_{X^{s,b'}_{1}}.
\end{align}
Using 
$
c_{h,\alpha}(\xi,\tau)
\doteq
(1+|\xi|)^{s}
(1+|\tau-\alpha\xi^3|)^{b'}
|\widehat{h}(\xi,\tau)|
$,
the bilinear estimate  \eqref{v-bi-est-Y-reduced} reads as follows
\begin{align}
\label{v-bilinear-est-L2-form-Y}
 \Big\| 
\int_{\rr^2}
 Q(\xi,\tau,\xi_1,\tau_1)
 c_{f,\alpha}(\xi-\xi_1,\tau-\tau_1) c_{g,1}(\xi_1,\tau_1) 
 d\xi_1 d\tau_1 \Big\|_{L^2_{\xi,\tau}}
\lesssim
\|c_{f,\alpha}\|_{L^2_{\xi,\tau}} \|c_{g,1}\|_{L^2_{\xi,\tau}},
\end{align}
where the multiplier $Q=Q(\xi,\tau,\xi_1,\tau_1)$ is given by
\begin{align*}
Q
\doteq&
\frac{ |\xi|\, \chi_{A}(\xi,\tau) }
{
(1+|\tau-\alpha\xi^3|)^{b}
}
\frac{
(1+|\tau|)^{\frac{s}{3}}
}
{(1+|\xi_1|)^{s}
(1+|\xi-\xi_1|)^{s}
}
\frac{1}{
(1+|\tau_1-\xi_1^3|)^{b'}
}
\,
\frac{1}{
(1+|\tau-\tau_1-\alpha(\xi-\xi_1)^3|)^{b'}
}.
\end{align*}
To make further simplification, we consider the following two microlocalizations:

\vskip0.05in
\noindent
$\bullet$ Microlocalization I:
$|\tau|< (10+10\alpha)|\xi_1|^3$
\qquad
$\bullet$ Microlocalization II:
$|\tau|\ge (10+10\alpha)|\xi_1|^3$.

\vskip0.05in
\noindent
{\bf Estimate in microlocalization I.} Since  $|\tau|\le (10+10\alpha)|\xi_1|^3$ and $s\ge 0$, we have 
$
\frac{(1+|\tau|)^{s/3}}{
(1+|\xi_1|)^{s}(1+|\xi-\xi_1|)^{s}
}
\lesssim
1,
$
which gives us 
$
Q(\xi,\tau,\xi_1,\tau_1)
\lesssim
Q_1(\xi,\tau,\xi_1,\tau_1),
$
where the multiplier  $Q_1$ is defined as follows 
\begin{align*}
Q_1(\xi,\tau,\xi_1,\tau_1)
\doteq
\frac{
|\xi|
\,
\chi_{A}(\xi,\tau)
}
{
(1+|\tau-\alpha\xi^3|)^{b}
}
\,
\frac{\chi_{|\tau|<(10+10\alpha)|\xi_1|^3}}{
(1+|\tau_1-\xi_1^3|)^{b'}
}
\,
\frac{1}{
(1+|\tau-\tau_1-\alpha(\xi-\xi_1)^3|)^{b'}
}.
\end{align*}
To show bilinear  estimate \eqref{v-bilinear-est-L2-form-Y} in this situation, using the Cauchy--Schwarz inequality with respect to 
$(\xi_1,\tau_1)$, and taking the supremum over  $(\xi,\tau)$, we get
\begin{align*}
\Big\| 
\int_{\rr^2}
Q_1(\xi,\tau,\xi_1,\tau_1)
c_{f,\alpha}(\xi  -  \xi_1,\tau  -  \tau_1) c_{g,1}(\xi_1,\tau_1) 
d\xi_1 d\tau_1 
\Big\|_{L^2_{\xi,\tau}}
\le
\|
\Theta_1
\|_{L^{\infty}_{\xi,\tau}}^{1/2}
\|c_{f,\alpha}\|_{L^2_{\xi,\tau}}
\|c_{g,1}\|_{L^2_{\xi,\tau}} ,
\end{align*}
which shows that the proof of the bilinear estimate \eqref{v-bilinear-est-L2-form-Y} follows from the next result.
\begin{lemma}
\label{Y-lem-1}
Let $\alpha>0$ and $\alpha\ne 1$.
If $\frac13<b'\le b<\frac12$, then we have
\begin{align}
\label{Y-theta1}
\Theta_1(\xi,\tau)
\doteq
\frac{
\xi^2
\,
\chi_{A}(\xi,\tau)
}
{
(1+|\tau-\alpha\xi^3|)^{2b}
}
\int_{\rr^2}
\frac{\chi_{|\tau|<(10+10\alpha)|\xi_1|^3}}{
(1+|\tau_1-\xi_1^3|)^{2b'}
}
\,
\frac{d\tau_1 d\xi_1}{
(1+|\tau-\tau_1-\alpha(\xi-\xi_1)^3|)^{2b'}
}
\lesssim
1.
\end{align}
\end{lemma}
\vskip0.05in
\noindent
{\bf Estimate in microlocalization II.}
The multiplier  $Q(\xi,\tau,\xi_1,\tau_1)$ is reduced to $Q_2$, where
\begin{align*}
Q_2(\xi,\tau,\xi_1,\tau_1)
\doteq
\frac{|\xi|
\,
\chi_{A}(\xi,\tau)
}{(1+|\tau-\alpha\xi^{3}|)^{b}}
\frac{\chi_{|\tau|>(10+10\alpha)|\xi_1|^3}(1+|\tau|)^{s/3}}{
(1+|\xi_1|)^{s}(1+|\xi-\xi_1|)^{s}
}
\frac{1}
{
(1+|\tau-\tau_1-\alpha(\xi-\xi_1)^{3}|)^{b'}
(1+|\tau_1-\xi_1^{3}|)^{b'}
}.
\end{align*}
Again, applying the Cauchy--Schwarz inequality with respect to 
$(\xi_1,\tau_1)$, and taking the supremum over  $(\xi,\tau)$, we arrive at
\begin{align*}
\Big\| 
\int_{\rr^2}
Q_2(\xi,\tau,\xi_1,\tau_1)
c_{f,\alpha}(\xi  -  \xi_1,\tau  - \tau_1) c_{g,1}(\xi_1,\tau_1) 
d\xi_1 d\tau_1 
\Big\|_{L^2_{\xi,\tau}}
\le
\|
\Theta_2
\|_{L^{\infty}_{\xi,\tau}}^{1/2}
\|c_{f,\alpha}\|_{L^2_{\xi,\tau}}
\|c_{g,1}\|_{L^2_{\xi,\tau}},
\end{align*}
which shows that the proof of the bilinear estimate \eqref{v-bilinear-est-L2-form-Y}  follows from the next result.
\begin{lemma}
\label{Y-lem-2}
Let $\alpha>0$ and $\alpha\ne 1$.
If $0\le s< 3$ and $\max\{\frac{2s+3}{18},\frac{4}{9}\}<b'\le b<\frac12$, then we have
\begin{align}
\label{Y-theta2}
\Theta_2(\xi,\tau)
\doteq&
\frac{ |\xi|^2(1+|\tau|)^{2s/3}
\,
\chi_{A}(\xi,\tau)
}{(1+|\tau-\alpha\xi^{3}|)^{2b}}
\\
\times
&\int_{\rr^2}  
\frac{\chi_{|\tau|\ge(10+10\alpha)|\xi_1|^3}}{
(1+|\xi_1|)^{2s}(1+|\xi-\xi_1|)^{2s}
}
\frac{1}
{
(1+|\tau-\tau_1-(\xi-\xi_1)^{3}|)^{2b'}
(1+|\tau_1-\xi_1^{3}|)^{2b'}
}
d\xi_1 d\tau_1 
\lesssim
1.
\nn
\end{align}
\end{lemma}

\vskip.05in
\noindent
{\bf Proof of Lemma \ref{Y-lem-1}.} For the $\tau_1$-integral in \eqref{Y-theta1}, applying calculus estimate \eqref{eq:calc_5} with $\ell=\ell'=b'$, $a=\tau-\alpha(\xi-\xi_1)^3$, $c=\xi_1^3$
we arrive at the following estimate
\begin{align}
\label{theta1-est-Y}
\Theta_1(\xi,\tau)
\lesssim
\frac{ |\xi|^2 \chi_{A}(\xi,\tau)}{(1+|\tau-\alpha\xi^{3}|)^{2b}}
\int_{\rr}
\frac{\chi_{|\tau|<(10+10\alpha)|\xi_1|^3}}
{
(1+|\tau-\alpha\xi^3-\widetilde{d}_{\alpha}(\xi,\xi_1)|)^{4b'-1}
}
d\xi_1,
\end{align}
where $\widetilde{d}_{\alpha}(\xi,\xi_1)$ is  the Bourgain quantity 
defined in \eqref{v-d3-quantity}, that is 
$
\widetilde{d}_{\alpha}(\xi,\xi_1)
\doteq
-\alpha\xi^3+\xi_1^3+\alpha(\xi-\xi_1)^3.
$
Now, our strategy is to split the $\xi_1$ integral at the critical points of $\mu(\xi_1)=\widetilde{d}_{\alpha}(\xi,\xi_1)$. Differentiating we get
\begin{align}
\label{v-d-partial-xi}
\frac{\p \widetilde{d}_{\alpha}(\xi,\xi_1)}{\p\xi_1}
=
3(1-\alpha)(\xi_1-\widetilde{p}_1)(\xi_1-\widetilde{p}_2),
\,\,
\text{with}
\,\,
\widetilde{p}_1
=
\frac{\sqrt{\alpha}}{\sqrt{\alpha}-1}\xi,
\,\,
\widetilde{p}_2
=
\frac{\sqrt{\alpha}}{\sqrt{\alpha}+1}\xi,
\,\,
\alpha>0
\,\,
\text{and}
\,\,
\alpha\neq 1.
\end{align}
Furthermore, in this microlocalization $\xi_1$ is away from the critical points.
In fact,  we have 
$
|\tau|
>
(10+10\alpha)^7\Big[2\big|\frac{\sqrt{\alpha}}{\sqrt{\alpha}-1}\big|
+
2
\big|
\frac{3\alpha}{\alpha-1}
\big|
+
2\Big]^3|\xi|^3
$
and
$
|\tau|<(10+10\alpha)|\xi_1|^3$, which implies that
$
(10+10\alpha)|\xi_1|^3
>
(10+10\alpha)^7\Big[2\big|\frac{\sqrt{\alpha}}{\sqrt{\alpha}-1}\big|
+
2
\big|
\frac{3\alpha}{\alpha-1}
\big|
+
2\Big]^3|\xi|^3,
$
or
\begin{equation}
\label{Y-s-0-xi1-away}
|\xi_1|
>
N|\xi|,
\quad
\text{with}
\quad
N
\doteq
(10+10\alpha)^2\Big[2\big|\frac{\sqrt{\alpha}}{\sqrt{\alpha}-1}\big|
+
2
\big|
\frac{3\alpha}{\alpha-1}
\big|
+
2\Big]|\xi|.
\end{equation}
Thus, using the two $\xi_1$-intervals 
$
\mathcal{I}_1^{\xi_1}
\doteq
(-\infty,-N|\xi|) 
$
and
$
\mathcal{I}_2^{\xi_1}
\doteq
(N|\xi|, \infty),
$
from \eqref{theta1-est-Y} we have 
\begin{equation*}
\Theta_1
=
J_1+J_2
\text{ with }
J_k(\xi,\tau)
\doteq
\frac{|\xi|^2 \chi_{A}(\xi,\tau)}{(1+|\tau-\alpha\xi^{3}|)^{2b}}
\int_{\mathcal{I}_k^{\xi_1}}
\frac{\chi_{|\tau|<(10+10\alpha)|\xi_1|^3}}
{
(1+|\tau-\alpha\xi^3-\widetilde{d}_{\alpha}(\xi,\xi_1)|)^{4b'-1}
}
d\xi_1.
\end{equation*}
In addition, for each one of the intervals $\mathcal{I}_1^{\xi_1}$, $\mathcal{I}_2^{\xi_1}$ we do further split. Using $\xi_1$-intervals
$
\mathcal{I}_{k,1}^{\xi_1}
\doteq
\big\{
\xi_1\in \mathcal{I}_{k}^{\xi_1};
\,\,
|\tau-\alpha\xi^3|\ge \frac12 |\widetilde{d}_{\alpha}(\xi,\xi_1)|
\big\}
$
and
$
\mathcal{I}_{k,2}^{\xi_1}
\doteq
\big\{
\xi_1\in \mathcal{I}_{k}^{\xi_1};
\,\,
|\tau-\alpha\xi^3|< \frac12 |\widetilde{d}_{\alpha}(\xi,\xi_1)|
\bigg\},
$
for 
$k=1,2$ we have
\begin{equation}
\label{Y-theta1-est-Jk}
J_k
=
J_{k1}+J_{k2}
\text{ with }
J_{kj}(\xi,\tau)
\doteq
\frac{|\xi|^2 \chi_{A}(\xi,\tau)}{(1+|\tau-\alpha\xi^{3}|)^{2b}}
\int_{\mathcal{I}_{kj}^{\xi_1}}
\frac{\chi_{|\tau|<(10+10\alpha)|\xi_1|^3}}
{
(1+|\tau-\alpha\xi^3-\widetilde{d}_{\alpha}(\xi,\xi_1)|)^{4b'-1}
}
d\xi_1.
\end{equation}
{\bf Estimate for $J_{k1}$.} 
Using $|\tau-\alpha\xi^3|\ge \frac12 |\widetilde{d}_{\alpha}(\xi,\xi_1)|$, we get 
$
|\tau-\alpha\xi^3|
\gtrsim
|\tau-\alpha\xi^3-\widetilde{d}_{\alpha}(\xi,\xi_1)|,
$
which gives
\begin{align}
\label{Y-theta1-est-Jk-no1}
J_{k1}(\xi,\tau)
\lesssim
\int_{\mathcal{I}_{k1}^{\xi_1}}
\frac{\chi_{A}(\xi,\tau)\, \chi_{|\tau|<(10+10\alpha)|\xi_1|^3}\,\,\, |\xi|^2}
{
(1+|\tau-\alpha\xi^3-\widetilde{d}_{\alpha}(\xi,\xi_1)|)^{2b+4b'-1}
}
d\xi_1
,
\quad
k=1,2.
\end{align}
Now making the change of variables  $\mu=\mu(\xi_1)=\tau-\alpha\xi^3-\widetilde{d}_{\alpha}(\xi,\xi_1)$ and using the property  
\eqref{v-d-partial-xi}, 
we get $|\mu'|\gtrsim \big|
\xi_1-\frac{\sqrt{\alpha}}{\sqrt{\alpha}-1}\xi
\big|
\big|
\xi_1-\frac{\sqrt{\alpha}}{\sqrt{\alpha}+1}\xi
\big|$, which combined with
estimate \eqref{Y-theta1-est-Jk-no1} implies that
\begin{align}
\label{Y-theta1-est-Jk-no2}
J_{k1}(\xi,\tau)
\lesssim
\int_{\mathcal{I}_{k1}^{\mu}}
\frac{\chi_{A}(\xi,\tau)  \,\chi_{|\tau|<(10+10\alpha)|\xi_1|^3}}
{
(1+|\mu|)^{2b+4b'-1}
}
\frac{|\xi|^2}{
\big|
\xi_1-\frac{\sqrt{\alpha}}{\sqrt{\alpha}-1}\xi
\big|
\big|
\xi_1-\frac{\sqrt{\alpha}}{\sqrt{\alpha}+1}\xi
\big|}
d\mu,
\quad
k=1,2.
\end{align}
Now, using estimate \eqref{Y-s-0-xi1-away}, i.e. $|\xi_1|\gg |\xi|$, we get $\big|
\xi_1-\frac{\sqrt{\alpha}}{\sqrt{\alpha}-1}\xi
\big|
\big|
\xi_1-\frac{\sqrt{\alpha}}{\sqrt{\alpha}+1}\xi
\big|\gtrsim \xi^2$, which combined with \eqref{Y-theta1-est-Jk-no2} gives
$
J_{k1}(\xi,\tau)(\xi,\tau)
\lesssim
\int_{\rr}
\frac{1}
{
(1+|\mu|)^{2b+4b'-1}
}
d\mu
$
($
k=1,2
$
).
This integral is bounded if 
$2b+4b'-1>1$. It suffices to have
$
1/3<b'\le b.
$

\vskip.05in
\noindent
{\bf Estimate for $J_{k2}$.} Using that 
$
|\tau-\alpha\xi^3|< \frac12 |\widetilde{d}_{\alpha}(\xi,\xi_1)|
$
we get
$
|\tau-\alpha\xi^3-\widetilde{d}_{\alpha}(\xi,\xi_1)|\gtrsim |\widetilde{d}_{\alpha}(\xi,\xi_1)|
$, which combined with  \eqref{Y-theta1-est-Jk} implies that  
\begin{align}
\label{Y-s-0-Jk2}
J_{k2}(\xi,\tau)
\lesssim
\frac{\chi_{A}(\xi,\tau) |\xi|^2}{(1+|\tau-\alpha\xi^{3}|)^{2b}}
\int_{\mathcal{I}_{k2}^{\xi_1}}
\frac{\chi_{|\tau|<(10+10\alpha)|\xi_1|^3}}
{
(1+|\widetilde{d}_{\alpha}(\xi,\xi_1)|)^{4b'-1}
}
d\xi_1,
\quad
k=1,2.
\end{align}
Furthermore, using estimate \eqref{Y-s-0-xi1-away}, i.e. $|\xi_1|\gg |\xi|$, and  the following factorization
\begin{equation}
\label{v-d-le-4-xi}
\widetilde{d}_{\alpha}(\xi,\xi_1)
=
(1-\alpha)\xi_1
\Big[\xi_1-\frac{3\alpha+\sqrt{3\alpha(4-\alpha)}}{2(\alpha-1)}\xi\Big]
\Big[\xi_1-\frac{3\alpha-\sqrt{3\alpha(4-\alpha)}}{2(\alpha-1)}\xi\Big]
,
\quad
0<\alpha\le 4
\,\,
\text{and}
\,\,
\alpha\neq 1,
\end{equation}
we get 
$
|\widetilde{d}_{\alpha}(\xi,\xi_1)|
\gtrsim
|\xi_1|^3
$ 
for $0<\alpha\le 4$ and $\alpha\neq 1$.
Also for
$
\alpha>4,
$ we have 
$
|\widetilde{d}_{\alpha}(\xi,\xi_1)|
\ge
\frac{\alpha-4}{4}
|\xi_1|^3.
$
In addition, for $(\xi,\tau)\in A$, we have 
$|\tau-\alpha\xi^3|\gtrsim |\tau|\gtrsim |\xi|^3$.
Combining these estimates with \eqref{Y-s-0-Jk2}, we get 
\begin{align*}
J_{k2}(\xi,\tau)
\lesssim
\frac{|\xi|^2}{(1+|\xi^{3}|)^{2b}}
\int_{\rr}
\frac{1}
{
(1+|\xi_1^3|)^{4b'-1}
}
d\xi_1
\lesssim
\frac{1}{(1+|\xi|)^{6b-2}}
\int_{\rr}
\frac{1}
{
(1+|\xi_1|)^{12b'-3}
}
d\xi_1,
\quad
k=1,2.
\end{align*}
For $\frac{1}{(1+|\xi|)^{6b-2}}\lesssim1 $, we need $6b-2\ge 0$ or $b\ge 1/3$. For  the integral of $\xi_1$ to be bounded, we need $12b'-3>1$ or $b'>1/3$. This completes the proof of Lemma \ref{Y-lem-1}.\,\, $\square$

\vskip.05in
\nin
{\bf Proof of Lemma \ref{Y-lem-2}.} 
For the $\tau_1$-integral in \eqref{Y-theta2}, applying calculus estimate \eqref{eq:calc_5} with $\ell=\ell'=b'$, $a=\tau-\alpha(\xi-\xi_1)^3$, $c=\xi_1^3$,
we arrive at the following estimate
\begin{align}
\label{Y-theta2-est}
\Theta_2(\xi,\tau)
\lesssim&
\frac{(1+|\tau|)^{2s/3}|\xi|^2}{(1+|\tau-\alpha\xi^{3}|)^{2b}}
\int_{\rr}
\frac{
\chi_{A}(\xi,\tau)
\chi_{|\tau|\ge(10+10\alpha)|\xi_1|^3}}{
(1+|\xi_1|)^{2s}(1+|\xi-\xi_1|)^{2s}
}
\cdot
\frac{1}
{
(1+|\tau-\alpha(\xi-\xi_1)^{3}-\xi_1^3|)^{4b'-1}
}
d\xi_1.
\end{align}
Using  $|\tau|>(10+10\alpha)|\xi_1|^3$ and $(\xi,\tau)\in A$, i.e. 
$
|\tau|>
(10+10\alpha)^7\big[2\big|\frac{\sqrt{\alpha}}{\sqrt{\alpha}-1}\big|
+
2
\big|
\frac{3\alpha}{\alpha-1}
\big|
+
2\big]^3|\xi|^3,
$
we get
\begin{equation}
\label{Y-theat2-xi-est}
|\tau-\alpha\xi^3|
\simeq
|\tau-\alpha(\xi-\xi_1)^{3}-\xi_1^3|
\simeq
|\tau|
\gtrsim
|\xi|^3+|\xi_1|^3,
\end{equation}
which combined with \eqref{Y-theta2-est} implies that 
\begin{align}
\label{Y-theta2-est-no1}
\Theta_2(\xi,\tau)
\lesssim
\frac{(1+|\tau|)^{2s/3}|\xi|^2}{(1+|\tau|)^{2b+4b'-1}}
\int_{\rr}
\frac{
\chi_{A}(\xi,\tau)
\chi_{|\tau|\ge(10+10\alpha)|\xi_1|^3}}{
(1+|\xi_1|)^{2s}(1+|\xi-\xi_1|)^{2s}
}
d\xi_1.
\end{align}
Next, we consider the following two cases.

\vskip.05in
\noindent
$\bullet$ Case 1: $0\le s\le 1$
\quad
$\bullet$ Case 2: $1< s<3$

\vskip0.05in
\noindent
{\bf Proof in Case 1.} Using estimate  \eqref{Y-theat2-xi-est}, we get $|\xi|\lesssim |\tau|^{1/3}$ and $|\xi_1|\lesssim |\tau|^{1/3}$, which helps us do the $\xi_1$-integration in \eqref{Y-theta2-est-no1}. We have
$
\Theta_2
\lesssim
\frac{(1+|\tau|)^{2s/3}|\xi|^2}{(1+|\tau|)^{2b+4b'-1}}
\int_{|\xi_1|\lesssim |\tau|^{1/3}}
1
d\xi_1
\lesssim
\frac{1}{(1+|\tau|)^{2b+4b'-2s/3-2}},
$
which is bounded if 
$
2b+4b'-2s/3-2
\ge
0
$
or
$
s
\le
3(b+2b'-1).
$
For $s\le 1$, it suffices to have $4/9\le b'\le b$.
This completes the proof in Case 1.

\vskip0.05in
\noindent
{\bf Proof in Case 2.} 
Using $|\xi|\lesssim\max\{|\xi_1|, |\xi-\xi_1|\}$ and $s>1$,
we get
$
\frac{|\xi|^2 }{
(1+|\xi_1|)^{2s}(1+|\xi-\xi_1|)^{2s}
}
$
$
\lesssim
$
$
\max\{
\frac{1}{(1+|\xi_1|)^{2s}},
$
$
\frac{1}{(1+|\xi-\xi_1|)^{2s}}
\},
$
which implies that 
$
\int_\rr
\frac{|\xi|^2 }{
(1+|\xi_1|)^{2s}(1+|\xi-\xi_1|)^{2s}
}
d\xi_1
\lesssim
1.
$
Combining this with  \eqref{Y-theta2-est-no1}, we get
\begin{align*}
\Theta_2(\xi,\tau)
\lesssim
\frac{(1+|\tau|)^{2s/3}}{(1+|\tau|)^{2b+4b'-1}}
\int_\rr
\frac{|\xi|^2 }{
(1+|\xi_1|)^{2s}(1+|\xi-\xi_1|)^{2s}
}
d\xi_1
\lesssim
\frac{(1+|\tau|)^{2s/3}}{(1+|\tau|)^{2b+4b'-1}}
=
\frac{1}{(1+|\tau|)^{2b+4b'-1-2s/3}},
\end{align*}
which is bounded if $2b+4b'-1-2s/3\ge 0$. For this, it suffices to have
$
b
\ge
b'>
\frac{3+2s}{18}.
$
For $b'\le b<\frac12$, we need $s< 3$.
This completes  the proof Lemma \ref{Y-lem-2}.
\,\,$\square$

%
%
%
%
%
%
%
%
%
%
\section{Counterexamples for bilinear estimates}
\label{sec:u-bi-counter}
\setcounter{equation}{0}
In this section we  provide counterexamples for
estimate \eqref{bi-est-X-1}, which corresponds 
to the nonlinearity of the $u$-equation.
We do this in the three cases:
$\bullet$  $\alpha=4$;
\quad
$\bullet$  $0<\alpha<4$ and $\alpha\neq 1$;
\quad
$\bullet$  $\alpha>4$ or $\alpha=1$.

%
%
%
%
%
%
%
%
%
%
\subsection{Counterexample for $\alpha=4$} 
We will show that  estimate \eqref{bi-est-X-1} fails for $s<\frac34$, if $\alpha=4$. In fact, we have the following result.
\begin{lemma}
\label{u4-counter-example-b-small}
Let $\alpha=4$.
If $s<\frac34$, then the bilinear estimate \eqref{bi-est-X-1} (with $b'=b$) fails for any  $b\in\rr$.
\end{lemma}
\nin
{\bf Proof  of Lemma \ref{u4-counter-example-b-small}.}
Using the notation given by \eqref{eq:c_u}, i.e.
$$
c_{f,\alpha}(\xi,\tau)
\doteq
\Big[
(1+|\xi|)^{s}
(1+|\tau-\alpha\xi^3|)^{b'}
+
\chi_{|\xi|< 1}(1+|\tau|)^{\theta'}
\Big]
|\widehat{f}(\xi,\tau)|,
$$
we see that
the $L^2$ formulation \eqref{bilinear-est-L2-form} of  bilinear estimate \eqref{bi-est-X-1}   when  $f=g$,
$\supp\, c_{f,\alpha}\subset\{(\xi_1,\tau_1)\in\rr^2: |\xi_1|\ge 1\}$ and dropping $\chi_{|\xi|<1}(1+|\tau|)^{\theta-1}$,
 reads as follows
\begin{align}
\label{u4-counter-1st-estimate}
\hskip-0.15in
\Big\|
\frac{\xi(1+|\xi|)^s}{(1+|\tau-\xi^3|)^{b}}
\hskip-0.05in
\int_{\rr^2}
\hskip-0.09in
\frac{
(1+|\xi-\xi_1|)^{-s}
(1+|\xi_1|)^{-s}
c_{f,4}(\xi-\xi_1,\tau-\tau_1) c_{f,4}(\xi_1,\tau_1)
d\xi_1 d\tau_1
}{
(1+|\tau-\tau_1-4(\xi-\xi_1)^{3}|)^{b}(1+|\tau_1-4\xi_1^{3}|)^{b}}
\Big\|_{L^2_{\xi,\tau}}
\hskip-0.05in
\hskip-0.07in
\lesssim
\hskip-0.03in
\big\|c_{f,4}\big\|_{L^2_{\xi_1,\tau_1}}^2.
\end{align}
We will show  that if 
bilinear estimate \eqref{u4-counter-1st-estimate} holds, then we must have  $s\ge\frac 34$. 
Since the Bourgain quantity 
$d_{4}(\xi,\xi_1)=12\xi(\xi_1-\xi/2)^2$ has a double zero at $\xi_1=\xi/2$, in our counterexample $\xi_1\simeq \xi/2$. 
In fact, for $N\in \zz^+$ we define  
$
c_{f,4}(\xi_1,\tau_1)
\doteq
\chi_{A^+_4}(\xi_1,\tau_1),
$
where  
$
A^+_4
\doteq
\big\{
(\xi_1,\tau_1)
\in
\rr^2
:
N-N^{-\frac{1}{2}}
\le
\xi_1
\le
N+N^{-\frac{1}{2}},
\,\,
|\tau_1-4\xi_1^3|
\le
10\cdot 5^3
\big\}.
$
Also, we define the $(\xi,\tau)$-domain $\widetilde{A}_1^+$ as follows
$
\widetilde{A}^+_1
\doteq
\big\{
(\xi,\tau)
\in
\rr^2
:
2N-(2N)^{-\frac{1}{2}}
\le
\xi
\le
2N+(2N)^{-\frac{1}{2}},
\,\,
|\tau-\xi^3|
\le
10\cdot 5^3
\big\}.
$
Therefore,  we get
\begin{align}
\label{u4-L2-norm-cf}
\|c_{f,4}\|_{L^2(\rr^2)}^2
=&
\int_{A^+_4} 1^2 d\xi_1d\tau_1
=
\int_{\xi_1=N-N^{-\frac{1}{2}}}^{N+N^{-\frac{1}{2}}}\int_{\tau_1=4\xi_1^3-10\cdot 5^3}^{4\xi_1^3+10\cdot 5^3}
1
d\tau_1d\xi_1
=
2\cdot10\cdot 5^3
\cdot
N^{-\frac{1}{2}}.
\end{align}
Next, defining the quantity in the $L^2$-norm of \eqref{u4-counter-1st-estimate}
\begin{align}
\label{u4-def-theta}
\Theta(\xi,\tau)
\doteq
\frac{\xi(1+|\xi|)^s}{(1+|\tau-\xi^3|)^{b}}
\int_{\rr^2}
\frac{
(1+|\xi-\xi_1|)^{-s}
(1+|\xi_1|)^{-s}
c_{f,4}(\xi-\xi_1,\tau-\tau_1) c_{f,4}(\xi_1,\tau_1)
d\xi_1 d\tau_1
}{
(1+|\tau-\tau_1-4(\xi-\xi_1)^{3}|)^{b}(1+|\tau_1-4\xi_1^{3}|)^{b}},
\end{align}
we derive  the following key estimate
\begin{align}
\label{u4-counter-example-bounded-from-below}
\| 
\Theta
\|_{L^2_{\xi,\tau}}
\gtrsim
N^{-s}N^{1/4}.
\end{align}
We prove estimate \eqref{u4-counter-example-bounded-from-below}  later.
Now, combining it with \eqref{u4-L2-norm-cf}, we see that if  the bilinear estimate \eqref{u4-counter-1st-estimate} holds, then  we would have
$
N^{-s}N^{1/4}
\lesssim
N^{-\frac{1}{2}}.
$
Since $N\gg 1$, we must have
$
-s+\frac14
\le
-\frac12
$
or
$
s
\ge
\frac34.
$
This completes the proof of  Lemma \ref{u4-counter-example-b-small}.
\,\,
$\square$

\nin
\begin{minipage}{0.5\linewidth}
\center
\begin{tikzpicture}[yscale=1, xscale=3,every node/.style={scale=1}]
%
%
\newcommand\X{0};
\newcommand\Y{0};
\newcommand\FX{9};
\newcommand\FY{9};
\newcommand\R{0.6};
\newcommand*{\TickSize}{2pt};
\newcommand*{\Num}{4};
\newcommand*{\Dist}{2}
\newcommand*{\Step}{\Dist/\Num}
%
%

\draw[line width=1pt,-{Latex[length=2mm]}]
(0.35,0)
--
(2.3,0)
node[above]
{\fontsize{\FX}{\FY} $\xi_1$}
;

\draw[line width=1pt,-{Latex[length=2mm]}]
(1,0)
--
(1,3.5)
node[right]
{\fontsize{\FX}{\FY} $\tau_1$};

\draw[red, line width=1pt, domain=1:1.5,variable=\x]
plot
({\x},{(\x)^3})
;

\draw[DarkGreen, line width=1pt, domain=1:1.5,variable=\x]
plot
({\x},{3*\x-2})
;

\draw[DarkGreen, line width=1pt, domain=1:1.5,variable=\x]
plot
({\x},{3*\x+1.5^3-4.5})
;


\draw[black, line width=0.4pt,dashed]
(1,0)
node[xshift=-0.4cm, yshift=-0.2cm]
{\fontsize{\FX}{\FY} $N$}

(1.5,0)
node[xshift=0.4cm,yshift=-0.3cm]
{\fontsize{\FX}{\FY} $N+\frac{1}{10^3}N^{-\frac{1}{2}}$}

(1.5,1.5^3)
node[xshift=-0.5cm]
{\fontsize{\FX}{\FY}\color{red} $C^*$}

(1,1)
node[below left]
{\fontsize{\FX}{\FY}\color{red} $E^*$}


%

;

\draw[blue, line width=1pt,dashed]
(1.5,0)
--
(1.5,1.5^3)

;

\draw[DarkGreen, line width=1pt, domain=1.23:1.5,variable=\x]
plot
({\x},{-1/3*\x+3})
;

\draw[DarkGreen, line width=1pt, domain=1:1.27,variable=\x]
plot
({\x},{-1/3*\x+2.2})
;

\draw[]
(1.5,1.35^3)
node[right]
{\fontsize{\FX}{\FY}\color{DarkGreen} $B^+$}

(1.25,1.35^3+0.1)
node[left]
{\fontsize{\FX}{\FY}\color{DarkGreen} $C^+$}

(1,{2.2-1/3})
node[left]
{\fontsize{\FX}{\FY}\color{DarkGreen} $D^+$}

(1.27,{2.2-1.27/3})
node[right]
{\fontsize{\FX}{\FY}\color{DarkGreen} $E^+$}

(1.7,1)
node[right]
{\fontsize{\FX}{\FY}\color{red} $\boxed{\tau_1=4\xi_1^3}$};

\draw [decorate,decoration={brace,amplitude=10pt},xshift=1pt,yshift=0pt]
(1.5,1.5^3)--(1.5,2.5)  node [black,midway,xshift=0.7cm] 
{\footnotesize $\le 1$};

\draw [decorate,decoration={brace,amplitude=10pt},xshift=-2pt,yshift=0pt]
(1,1)--(1,{2.2-1/3})  node [black,midway,xshift=-0.7cm] 
{\footnotesize $\le 1$};

\end{tikzpicture}

\vskip-0.15in

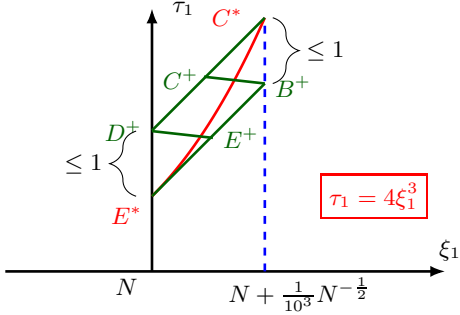
\captionof{figure}{Rectangle in $A_4^+$}
\label{fig:rec-A4}
\end{minipage}
\begin{minipage}{0.5\linewidth}
\center
\begin{tikzpicture}[yscale=1, xscale=3,every node/.style={scale=1}]
%
%
\newcommand\X{0};
\newcommand\Y{0};
\newcommand\FX{9};
\newcommand\FY{9};
\newcommand\R{0.6};
\newcommand*{\TickSize}{2pt};
\newcommand*{\Num}{4};
\newcommand*{\Dist}{2}
\newcommand*{\Step}{\Dist/\Num}
%
%

\draw[line width=1pt,-{Latex[length=2mm]}]
(0.35,0)
--
(2.3,0)
node[above]
{\fontsize{\FX}{\FY} $\xi$}
;

\draw[line width=1pt,-{Latex[length=2mm]}]
(1,0)
--
(1,3.5)
node[right]
{\fontsize{\FX}{\FY} $\tau$};

\draw[red, line width=1pt, domain=1:1.5,variable=\x]
plot
({\x},{(\x)^3})
;

\draw[DarkGreen, line width=1pt, domain=1:1.5,variable=\x]
plot
({\x},{3*\x-2})
;

\draw[DarkGreen, line width=1pt, domain=1:1.5,variable=\x]
plot
({\x},{3*\x+1.5^3-4.5})
;


\draw[black, line width=0.4pt,dashed]
(1,0)
node[xshift=-0.4cm, yshift=-0.2cm]
{\fontsize{\FX}{\FY} $2N$}

(1.5,0)
node[xshift=0.4cm,yshift=-0.3cm]
{\fontsize{\FX}{\FY} $2N+\frac{1}{10^3}(2N)^{-\frac{1}{2}}$}

(1.5,1.5^3)
node[xshift=-0.5cm]
{\fontsize{\FX}{\FY}\color{red} $\widetilde{C}^*$}

(1,1)
node[below left]
{\fontsize{\FX}{\FY}\color{red} $\widetilde{E}^*$}


%

;

\draw[blue, line width=1pt,dashed]
(1.5,0)
--
(1.5,1.5^3)

;

\draw[DarkGreen, line width=1pt, domain=1.23:1.5,variable=\x]
plot
({\x},{-1/3*\x+3})
;

\draw[DarkGreen, line width=1pt, domain=1:1.27,variable=\x]
plot
({\x},{-1/3*\x+2.2})
;

\draw[]
(1.5,1.35^3)
node[right]
{\fontsize{\FX}{\FY}\color{DarkGreen} $\widetilde{B}^+$}

(1.25,1.35^3+0.1)
node[left]
{\fontsize{\FX}{\FY}\color{DarkGreen} $\widetilde{C}^+$}

(1,{2.2-1/3})
node[left]
{\fontsize{\FX}{\FY}\color{DarkGreen} $\widetilde{D}^+$}

(1.27,{2.2-1.27/3})
node[right]
{\fontsize{\FX}{\FY}\color{DarkGreen} $\widetilde{E}^+$}

(1.7,1)
node[right]
{\fontsize{\FX}{\FY}\color{red} $\boxed{\tau=\xi^3}$};

\draw [decorate,decoration={brace,amplitude=10pt},xshift=1pt,yshift=0pt]
(1.5,1.5^3)--(1.5,2.5)  node [black,midway,xshift=0.7cm] 
{\footnotesize $\le 1$};

\draw [decorate,decoration={brace,amplitude=10pt},xshift=-2pt,yshift=0pt]
(1,1)--(1,{2.2-1/3})  node [black,midway,xshift=-0.7cm] 
{\footnotesize $\le 1$};

\end{tikzpicture}

\vskip-0.15in
\captionof{figure}{Rectangle in $\widetilde{A}_1^+$}
\label{fig:rec-A1}
\end{minipage}

\vskip.1in
\noindent
{\bf Proof of \eqref{u4-counter-example-bounded-from-below}.} 
Following \cite{kpv1996},
we can show that there is rectangle with dimension $c_1N^{-2}\times c_2N^{3/2}$, where $c_1$,$c_2$ are appropriate constants, in the domain $A^+_4$. In fact, as shown in  Figure \ref{fig:rec-A4}, we denote this rectangle by $B^+C^+D^+E^+$, where 
$B^+E^+$  is on the tangent line of $\tau=4\xi^3$ at $E^*=(N,4N^3)$. The coordinate for $B^+$ is  $(N+\frac{1}{10^3}N^{-1/2}, 4 N^3+\frac{12}{10^3} N^{\frac32})$
and the coordinate  for $D^+$ is 
$(N, 4(N+\frac{N^{-\frac12}}{10^3})^3
-
\frac{12}{10^3}N^{\frac32})$. 
This rectangle in in $A_4^+$.
Similarly, as shown in Figure \ref{fig:rec-A1},  we have a rectangle with  dimension $c_1N^{-2}\times c_2N^{3/2}$ in the domain $\widetilde{A}^+_1$ and we denote it by $\widetilde{B}^+\widetilde{C}^+\widetilde{D}^+\widetilde{E}^+$.
Next, we prove that if  $(\xi,\tau)\in  \widetilde{B}^+\widetilde{C}^+\widetilde{D}^+\widetilde{E}^+$ and $(\xi_1,\tau_1)\in B^+C^+D^+E^+$,  then $(\xi-\xi_1,\tau-\tau_1)\in A^+_4$. 
For $(\xi,\tau)\in \widetilde{B}^+\widetilde{C}^+\widetilde{D}^+\widetilde{E}^+$, we get
\begin{align} 
\label{u4-counter-xi-tau-est}
2N\le \xi\le 2N+\frac{1}{10^3}(2N)^{-\frac{1}{2}}
\quad
\text{and}
\quad
|\tau-\xi^3|\le 1.
\end{align}
Also, for $(\xi_1,\tau_1)\in B^+C^+D^+E^+$, we have
\begin{align} 
\label{u4-counter-xi1-tau1-est}
N\le \xi_1\le N+\frac{1}{10^3}N^{-\frac{1}{2}}
\quad
\text{and}
\quad
|\tau_1-4\xi_1^3|\le 1.
\end{align}
Using estimates \eqref{u4-counter-xi-tau-est} and \eqref{u4-counter-xi1-tau1-est}, for $(\xi,\tau)\in \widetilde{B}^+\widetilde{C}^+\widetilde{D}^+\widetilde{E}^+$ and $(\xi_1,\tau_1)\in B^+C^+D^+E^+$ we obtain
\begin{align}
\label{u4-counter-xi-xi1-est}
N-\frac{1}{10^3} N^{-\frac{1}{2}}
\le 
\xi-\xi_1
\le
N+\frac{1}{10^3}(2N)^{-\frac{1}{2}}
\quad
\text{ and }
\quad
-\frac12
\frac{1}{10^3} (2N)^{-\frac{1}{2}}
\le 
\xi_1-\frac{1}{2}\xi
\le
\frac{1}{10^3} N^{-\frac{1}{2}},
\end{align}
which combined with triangle inequality gives us that 
\begin{align}
\label{u4-counter-xi-xi1-tau-tau_1-est}
|\tau-\tau_1-4(\xi-\xi_1)^3|
\le&
|\tau-\xi^3|
+
|\tau_1-4\xi_1^3|
+
|d_4(\xi,\xi_1)|
\nn
\\
=&
|\tau-\xi^3|
+
|\tau_1-4\xi_1^3|
+
|12\xi(\xi_1-\xi/2)^2|
\le
3.
\end{align}
Estimates \eqref{u4-counter-xi-xi1-est} and \eqref{u4-counter-xi-xi1-tau-tau_1-est} implies that $(\xi-\xi_1,\tau-\tau_1)\in A^+_4$.
Now, letting $(\xi,\tau)\in \widetilde{B}^+\widetilde{C}^+\widetilde{D}^+\widetilde{E}^+$ and $(\xi_1,\tau_1)\in B^+C^+D^+E^+$, using $c_f(\xi-\xi_1,\tau-\tau_1)=1$ and $c_g(\xi_1,\tau_1)=1$, from \eqref{u4-def-theta} we have 
\begin{align}
\label{u4-counter-theta-est}
&\| 
\Theta
\|_{L^2_{\xi,\tau}}
\ge
\Big\|\chi_{(\xi,\tau)\in \widetilde{B}^+\widetilde{C}^+\widetilde{D}^+\widetilde{E}^+}
\frac{\xi(1+|\xi|)^s}{(1+|\tau-\xi^3|)^{b}}
\\
\times&
\int_{(\xi_1,\tau_1)\in B^+C^+D^+E^+}\frac{(1+|\xi-\xi_1|)^{-s}(1+|\xi_1|)^{-s}}{(1+|\tau-\tau_1-4(\xi-\xi_1)^{3}|)^{b}(1+|\tau_1-4\xi_1^{3}|)^{b}}d\xi_1 d\tau_1
\Big\|_{L^2_{\xi,\tau}}.
\nn
\end{align}
Since the area of rectangles $B^+C^+D^+E^+$, $\widetilde{B}^+\widetilde{C}^+\widetilde{D}^+\widetilde{E}^+$ are $c N^{-1/2}$,
using
estimates \eqref{u4-counter-xi-tau-est}--\eqref{u4-counter-theta-est}, we get 
\begin{align*}
\| 
\Theta
\|_{L^2_{\xi,\tau}}
\ge&
\Big\|
\frac{N(1+N)^{s}}{1}
\int_{(\xi_1,\tau_1)\in B^+C^+D^+E^+}
\frac{(1+N)^{-2s}}{1}
d\xi_1 d\tau_1
\cdot
\chi_{(\xi,\tau)\in \widetilde{B}^+\widetilde{C}^+\widetilde{D}^+\widetilde{E}^+}
\Big\|_{L^2_{\xi,\tau}}
\gtrsim
N^{-s}N^{\frac14},
\end{align*}
which is the desired estimate \eqref{u4-counter-example-bounded-from-below}. 
\,\,
$\square$

%
%
%
%
%
%
%
%
%
%
%
\subsection{Counterexample for $0<\alpha<4$ and $\alpha\neq 1$.}
In this subsection, we prove the following result
\begin{lemma}
\label{u0-counter-example-b-small}
Let $0<\alpha<4$ and $\alpha\neq 1$.
If $s<0$, then the bilinear estimate \eqref{bi-est-X-1} (with $b'=b$) fails for any  $b\in\rr$.
\end{lemma}
\nin
{\bf Proof  of Lemma \ref{u0-counter-example-b-small}.}
Using the notation $c_{f,\alpha}$, $c_{g,\alpha}$ defined in  \eqref{eq:c_u}, we see that
the $L^2$ formulation \eqref{bilinear-est-L2-form}
of  bilinear estimate \eqref{bi-est-X-1} when
$\supp\, c_{f,\alpha}\subset \{(\xi_2,\tau_2)\in\rr^2 : |\xi_2|\ge 1\}$, $\supp\, c_{g,\alpha}\subset \{(\xi_1,\tau_1)\in\rr^2 : |\xi_1|\ge 1\}$ and dropping $\chi_{|\xi|<1}(1+|\tau|)^{\theta-1}$,
reads as follows
\begin{align}
\label{u0-counter-1st-estimate}
&
\Big\| \frac{\xi(1+|\xi|)^s}{(1+|\tau-\xi^3|)^{b}}
\int_{\rr^2}\frac{
(1+|\xi-\xi_1|)^{-s}(1+|\xi_1|)^{-s}
c_{f,\alpha}(\xi-\xi_1,\tau-\tau_1) c_{g,\alpha}(\xi_1,\tau_1)}{(1+|\tau-\tau_1-\alpha(\xi-\xi_1)^{3}|)^{b}(1+|\tau_1-\alpha\xi_1^{3}|)^{b}}d\xi_1 d\tau_1
 \Big\|_{L^2_{\xi,\tau}}
 \nn
\\
\lesssim&
\big\|c_{f,\alpha}\big\|_{L^2_{\xi,\tau}}
\big\|c_{g,\alpha}\big\|_{L^2_{\xi,\tau}}.
\end{align}
Like in the counterexample for $\alpha=4$, in this counterexample $\xi_1$ is near the zero of  the Bourgain quantity
$
d_{\alpha}(\xi,\xi_1)
=
\xi(\xi_1-r_1\xi)(\xi_1-r_2\xi),
$
where 
$
r_1
=
\frac12-\frac{\sqrt{-3+12\alpha^{-1}}}{6}
$
and
$
r_2
=
\frac12+\frac{\sqrt{-3+12\alpha^{-1}}}{6}
$
(see property \eqref{d-le-4-xi1}).
In fact, 
letting  $N\in \zz^+$, we define
$
c_{f,\alpha}(\xi_2,\tau_2)
\doteq
\chi_{A_\alpha}(\xi_2,\tau_2),
$
where  
$
A_\alpha
\doteq
\big\{
(\xi_2,\tau_2)
\in
\rr^2
:
(1-r_2)N-2N^{-2}
\le
\xi_2
\le
(1-r_2)N+2N^{-2},
\alpha [(1-r_2)N]^3
-
10^3
\le
\tau_2
\le
\alpha [(1-r_2)N]^3
+
10^3
\big\}.
$
Also, we define
$
c_{g,\alpha}(\xi_1,\tau_1)
\doteq
\chi_{\widetilde{A}_\alpha}(\xi_1,\tau_1),$%
where  
$
\widetilde{A}_\alpha
\doteq
\big\{
(\xi_1,\tau_1)
\in
\rr^2
:
r_2N-\frac{1}{(10+|r_1|+|r_2|)^3} N^{-2}
\le
\xi_1
\le
r_2N+\frac{1}{(10+|r_1|+|r_2|)^3} N^{-2},
\alpha(r_2N)^3
-
1
\le
\tau_2
\le
\alpha(r_2N)^3
+
1
\big\}.
$
Furthermore, we define the domain for $(\xi,\tau)$  as follows
$
\Delta
=
\big\{
(\xi,\tau)\in\rr^2
:
N-\frac{1}{(10+|r_1|+|r_2|)^3}N^{-2}<\xi<N+\frac{1}{(10+|r_1|+|r_2|)^3}  N^{-2},
 N^3-1 \le \tau\le  N^3+ 1
\big\}.
$
All domains  $A_\alpha$, $\widetilde A_\alpha$ and $\Delta$ are right rectangle with dimensions $\frac{\varepsilon}{N^2}\times c$, where $\varepsilon$ is a small number related to the zeros of Bourgain quantity $d_{\alpha}$ and $c$ is a constant. 
And, the centers of these rectangles are at $((1-r_2)N, \alpha[(1-r_2)N]^3)$, $(r_2N, \alpha (r_2N)^3)$, $(N, N^3)$, respectively.
In addition, we have
\begin{align}
\label{u0-L2-norm-cf}
\|c_f\|_{L^2_{\xi_2,\tau_2}}^2
\lesssim
N^{-2}
\quad
\text{and}
\quad
\|c_g\|_{L^2_{\xi_1,\tau_1}}^2
\lesssim
N^{-2}.
\end{align}
Next, denoting the quantity in the $L^2$-norm of estimate \eqref{u0-counter-1st-estimate} by
\begin{align}
\label{u0-counter-def-theta}
\Theta(\xi,\tau)
\doteq&
\frac{\xi(1+|\xi|)^s}{(1+|\tau-\xi^3|)^{b}}
\int_{\rr^2}\frac{
(1+|\xi-\xi_1|)^{-s}(1+|\xi_1|)^{-s}
c_{f,\alpha}(\xi-\xi_1,\tau-\tau_1) c_{g,\alpha}(\xi_1,\tau_1)}{(1+|\tau-\tau_1-\alpha(\xi-\xi_1)^{3}|)^{b}(1+|\tau_1-\alpha\xi_1^{3}|)^{b}}d\xi_1 d\tau_1,
\end{align}
we obtain the following key estimate for $\Theta$
\begin{align}
\label{u0-counter-example-bounded-from-below}
\| 
\Theta
\|_{L^2_{\xi,\tau}}
\gtrsim
N^{-s-2}.
\end{align}
We prove estimate \eqref{u0-counter-example-bounded-from-below} later. Now, combining estimate \eqref{u0-counter-example-bounded-from-below} with \eqref{u0-L2-norm-cf}, we see that if  the bilinear estimate \eqref{u0-counter-1st-estimate} holds, then  we must have
$
N^{-s-2}
\lesssim
N^{-2}.
$
Since $N\gg 1$, we must have
$
-s-2
\le
-2
$
or
$s\ge 0$.
This completes the proof of  Lemma \ref{u0-counter-example-b-small}, once we prove estimate \eqref{u0-counter-example-bounded-from-below}.
\,\,
$\square$

\vskip0.05in
\noindent
{\bf Proof of Estimate \eqref{u0-counter-example-bounded-from-below}.} First, we prove that if $(\xi,\tau)\in \Delta$ and $(\xi_1,\tau_1)\in \widetilde{A}_{\alpha}$, then we have $(\xi-\xi_1,\tau-\tau_1)\in A_{\alpha}$.
In fact, for $(\xi,\tau)\in \Delta$ and $(\xi_1,\tau_1)\in \widetilde{A}_{\alpha}$ we have
\begin{align}
\label{u0-counter-xi-xi1-est}
(1-r_2)N
-
\frac{2}{(10+|r_1|+|r_2|)^3 }N^{-2}
\le
\xi-\xi_1
\le
(1-r_2)N
+
\frac{2}{(10+|r_1|+|r_2|)^3} N^{-2}.
\end{align}
Also, we have $\tau-\tau_1-\alpha(\xi-\xi_1)^3
=
(\tau-\xi^3)-(\tau_1-\alpha\xi_1^3)-
d_{\alpha}(\xi,\xi_1)$. Using the triangle inequality and the fact   that $\xi_1$ is near the zero of $d_\alpha(\xi,\xi_1)$, we get 
$
|\tau-\tau_1-\alpha(\xi-\xi_1)^3|
\le
10^3,
$
which combined 
 with estimate \eqref{u0-counter-xi-xi1-est} implies that $(\xi-\xi_1,\tau-\tau_1)\in A_{\alpha}$. Therefore, for $(\xi,\tau)\in \Delta$ and $(\xi_1,\tau_1)\in \widetilde{A}_{\alpha}$ we have 
$c_{f,\alpha}(\xi-\xi_1,\tau-\tau_1)=1$ and $c_{g,\alpha}(\xi_1,\tau_1)=1$. Thus, from \eqref{u0-counter-def-theta} we get 
\begin{align*}
\| 
\Theta(\xi,\tau)
\|_{L^2_{\xi,\tau}}
\ge&
\Big\|
\frac{\xi
(1+|\xi|)^s}{(1+|\tau-\xi^3|)^{b}}
\int_{\widetilde{A}_\alpha}\frac{
(1+|\xi-\xi_1|)^{-s}(1+|\xi_1|)^{-s}
}{(1+|\tau-\tau_1-\alpha(\xi-\xi_1)^{3}|)^{b}(1+|\tau_1-\alpha\xi_1^{3}|)^{b}}d\xi_1 d\tau_1
\cdot
\chi_{ (\xi,\tau)\in \Delta}
\Big\|_{L^2_{\xi,\tau}}
\\
\gtrsim&
\Big\|
\frac{N(1+N)^s}{1}
\int_{\widetilde{A}_\alpha}
\frac{(1+N)^{-2s}}{1}
d\xi_1d\tau_1
\cdot
\chi_{ (\xi,\tau)\in \Delta}
\Big\|_{L^2_{\xi,\tau}}
\gtrsim
N^{-2-s},
\end{align*}
since the area of the rectangles $\Delta$ and  $\widetilde{A}_\alpha$ are $c N^{-2}$.    This is  the desired estimate \eqref{u0-counter-example-bounded-from-below}.
\,\,
$\Box$

%
%
%
%
%
%
%
%
%
%
\subsection{Counterexample for $\alpha>4$ or $\alpha=1$.}
We prove that bilinear estimate \eqref{bi-est-X-1} fails for $s<-\frac34$. More precisely, we have the following result.
\begin{lemma}
\label{counter-example-b-small}
Let $\alpha>0$.
If $s<-\frac34$ and the bilinear estimate \eqref{bi-est-X-1} 
 holds (with $b'=b$), then $b>\frac12$.
\end{lemma}
\nin
{\bf Proof  of Lemma \ref{counter-example-b-small}.}
Using the notation  $c_{f,\alpha}(\xi,\tau)$
given by \eqref{eq:c_u}
we see that
the $L^2$ formulation \eqref{bilinear-est-L2-form}
of  bilinear estimate \eqref{bi-est-X-1} when $f=g$,
$\supp\, c_{f,\alpha}\subset\{(\xi_1,\tau_1)\in\rr^2 : |\xi_1|\ge 1\}$ and dropping $\chi_{|\xi|<1}(1+|\tau|)^{\theta-1}$,
 reads as follows
%
\begin{align}
\label{counter-1st-estimate}
\hskip-0.17in
\Big\|
\frac{\xi(1+|\xi|)^s}{(1+|\tau-\xi^3|)^{b}}
\hskip-0.05in
\int_{\rr^2}
\hskip-0.09in
\frac{
(1+|\xi-\xi_1|)^{-s}
(1+|\xi_1|)^{-s}
c_{f,\alpha}(\xi-\xi_1,\tau-\tau_1) c_{f,\alpha}(\xi_1,\tau_1)
d\xi_1 d\tau_1
}{
(1+|\tau-\tau_1-\alpha(\xi-\xi_1)^{3}|)^{b}(1+|\tau_1-\alpha\xi_1^{3}|)^{b}}
\Big\|_{L^2_{\xi,\tau}}
\hskip-0.05in
\hskip-0.07in
\lesssim
\hskip-0.03in
\big\|c_{f,\alpha}\big\|_{L^2_{\xi_1,\tau_1}}^2.
\end{align}
We will show  that if $s<-\frac 34$ and bilinear estimate \eqref{counter-1st-estimate} holds, then we must have  $b>\frac12$.  
Let $N\in \zz^+$ and $c_{f,\alpha}$ is defined by
$
c_{f,\alpha}(\xi_1,\tau_1)
\doteq
\chi_{A^+_\alpha}(\xi_1,\tau_1)+\chi_{A^-_\alpha}(\xi_1,\tau_1) \in L^2_{\xi_1,\tau_1},
$
where  $\chi_{A^+_\alpha}(\cdot)$ is the characteristic function of the set $A^+_\alpha$ 
given by
$
A^+_\alpha
\doteq
\big\{
(\xi_1,\tau_1)
\in
\rr^2
:
N-N^{-\frac{1}{2}}
\le
\xi_1
\le
N+N^{-\frac{1}{2}},
\,\,
|\tau_1-\alpha\xi_1^3|
\le
10(1+\alpha)^3
\big\}
$
and
$
A^-_\alpha
=
-A^+_\alpha
=
\big\{
(\xi_1,\tau_1)
\in
\rr^2
:
(-\xi_1,-\tau_1)\in A^{+}_\alpha
\big\}.
$
By a straightforward  computation we get $c_{f,\alpha} \in L^2_{\xi_1,\tau_1}$.
In fact,  we have
\begin{align}
\label{L2-norm-cf}
\|c_{f,\alpha}\|_{L^2}^2
=
\int_{A^+_\alpha\cup A^-_\alpha} 1^2 d\xi_1d\tau_1
=
2\int_{A^+_\alpha} 1 d\xi_1d\tau_1
=
2\int_{N-N^{-\frac{1}{2}}}^{N+N^{-\frac{1}{2}}}\int_{\alpha\xi_1^3-10(1+\alpha)^3}^{\alpha\xi_1^3+10(1+\alpha)^3}
1
d\tau_1d\xi_1
\lesssim
N^{-\frac{1}{2}}.
\end{align}
Next,  following the arguments in \cite{kpv1996} 
we restrict  $(\xi_1, \tau_1)$ in an appropriate rectangle  contained 
in $A_\alpha^+$, and $(\xi, \tau)$ in a similar rectangle centered at the origin, so that $(\xi-\xi_1,\tau-\tau_1)\in A_\alpha^-$.
Then, we bound the left-hand side of  \eqref{counter-1st-estimate}
by a quantity which is bigger that   $N^{-2s-\frac{3}{2}b-5/4}$,
thus obtaining the estimate
\begin{align}
\label{counter-example-bounded-from-below}
\hskip-0.15in
\Big\|
\frac{\xi(1
\hskip-0.03in + \hskip-0.03in
|\xi|)^s}{(1
\hskip-0.03in + \hskip-0.03in
|\tau-\xi^3|)^{b}}
\hskip-0.05in
\int_{\rr^2}
\hskip-0.09in
\frac{
(1+|\xi
\hskip-0.02in - \hskip-0.02in
\xi_1|)^{-s}
(1+|\xi_1|)^{-s}
c_{f,\alpha}(\xi
\hskip-0.02in - \hskip-0.02in
\xi_1,\tau
\hskip-0.02in - \hskip-0.02in
\tau_1) c_{f,\alpha}(\xi_1,\tau_1)
d\xi_1 d\tau_1
}{
(1+|\tau-\tau_1-\alpha(\xi-\xi_1)^{3}|)^{b}(1+|\tau_1-\alpha\xi_1^{3}|)^{b}}
\Big\|_{L^2_{\xi,\tau}}
\hskip-0.12in
\gtrsim
\hskip-0.03in
N^{-2s-\frac{3}{2}b-5/4}.
\end{align}
Using estimates \eqref{L2-norm-cf} and \eqref{counter-example-bounded-from-below}  we see that if  the bilinear estimate \eqref{counter-1st-estimate} holds, then  we would have
$
N^{-2s-\frac{3}{2}b-5/4}
\lesssim
N^{-\frac{1}{2}}.
$
Since $N\gg 1$, we must have
$
-2s-\frac{3}{2}b-\frac{5}{4}
\le
-\frac12
$
or
$
b\ge -\frac12-\frac43s.
$
So, if $s<-\frac34$, then 
$
b
>
\frac12.
$
This completes the proof of  Lemma \ref{counter-example-b-small}.
\,\,
$\square$

%
%
%
%
%
%
%
%
%
%
%
\section{Derivation of Fokas UTM formula for the Dirichlet  Problem
}
\label{sec:utm-deriv}
\setcounter{equation}{0}
In this section, we provide an outline of the Fokas method
for solving 
ibvp \eqref{LKdV-v}, which we restate  as follows by using a more standard notation
\begin{subequations}
\label{L-alpha-KdV}
\begin{align}
\label{LKdV eqn}
&\p_tu+\alpha\p_x^3u
=
f,
\quad
0<x<\infty,
\,\,
0<t<T,
\\
\label{LKdV iv}
&u(x,0)
=
u_0(x),
\quad
0<x<\infty,
\\
\label{LKdV bc}
&u(0,t)
=
g_0(t),
\quad
0<t<T.
\end{align}
\end{subequations}
Utilizing the exponential solutions $\widetilde{u}=e^{-i\xi x-i\alpha\xi^3t}$ ($\xi\in \cc$) to the adjoint  equation
$
\partial_t \widetilde{u}
+
\alpha\partial^3_x \widetilde{u}
=
0
$
and doing some manipulations we obtain the divergence form
\begin{align}
\label{divergence form}
(e^{-i\xi x-i\alpha\xi^3t}u)_t+\alpha(e^{-i\xi x-i\alpha\xi^3t}[\p_x^2 u+i\xi\p_x u-\xi^2 u])_x
=
e^{-i\xi x-i\alpha\xi^3t}f.
\end{align}
Integrating the divergence form \eqref{divergence form} from $x=0$ to $\infty$
(assuming vanishing at  $\infty$)
 gives
\begin{align}
\label{tode}
(e^{-i\alpha\xi^3t}\hat u(\xi,t))_t
=
e^{-i\alpha\xi^3t}\hat f(\xi,t)
+
\alpha
e^{-i\alpha\xi^3t}g(\xi,t),
\end{align}
where the half-line Fourier transforms
$\hat u(\xi,t)$, $\hat f(\xi,t)$ are defined by \eqref{half-line-FT-ic}, 
and
$
g(\xi,t)
\doteq
\p_x^2u(0,t)+i\xi\p_xu(0,t)-\xi^2u(0,t).
$
Integrating \eqref{tode} from $0$ to $t$, $0< t< T$, we find the following global relation:
\begin{align}
\label{global relation}
e^{-i\alpha\xi^3t}\hat u(\xi,t)
=
\hat u_0(\xi)+
\alpha[\widetilde g_2(\xi^3,t)+i\xi\widetilde g_1(\xi^3,t)
-\xi^2\widetilde g_0(\xi^3,t)]
+
F(\xi,t),
\quad
{\rm Im}\xi\leq 0,
\end{align}
where $F$ is given by \eqref{time-trans-F2} and 
$
\widetilde g_j(\xi,t)
\doteq
\int_0^te^{-i\alpha\xi\tau}\p_x^j u(0,\tau)d\tau.
$
Inverting equation \eqref{global relation} we find the following integral representation for the solution to the forced ibvp \eqref{L-alpha-KdV}
\begin{align}
\label{sln-repr}
u(x,t)
=
\frac{1}{2\pi}\int_{-\infty}^\infty e^{i\xi x+i\alpha\xi^3t}
[\hat u_0(\xi)+F(\xi,t)]d\xi
+
\frac{\alpha}{2\pi}\int_{-\infty}^\infty e^{i\xi x+i\alpha\xi^3t}
\widetilde{g}(\xi^3,t)d\xi,
\end{align}
where 
$
\widetilde{g}
\doteq 
\widetilde g_2(\xi^3,t)+i\xi\widetilde g_1(\xi^3,t)
-\xi^2\widetilde g_0(\xi^3,t).
$

Next, we notice that the integral representation \eqref{sln-repr} includes the 
time transforms $\widetilde g_1$ and $\widetilde g_2$ of the unknown data  
$\p_xu(0,t)$ and $\p_x^2u(0,t)$, which we can not remove
 if the $\xi$-integral is over the real line. 
 So, following Fokas method (see \cite{f1997}), we use the analyticity of the half-line Fourier transform to deform the 
 $\xi$-integration  from the real line to $\p D^+$, where $D^+$ is the domain shown in Figure \ref{kdv-domain}.  Then, the  integral representation \eqref{sln-repr} becomes
\begin{align}
\label{sln-defrom}
u(x,t)
=
\frac{1}{2\pi}\int_{-\infty}^\infty e^{i\xi x+i\alpha\xi^3t}
[\hat u_0(\xi)+F(\xi,t)]d\xi
+
\frac{\alpha}{2\pi}\int_{\p D^+} e^{i\xi x+i\alpha\xi^3t}
\widetilde g(\xi^3,t)d\xi.
\end{align}
Now, we are able to remove the unknown data. In fact, using the  symmetry transformation $\xi \mapsto \sigma \xi$  and 
$\xi \mapsto \sigma^2 \xi$, where  
$\sigma=e^{\frac{2\pi}{3}i}$,   from the global relation 
\eqref{global relation} we get
\begin{align}
\label{glo-1}
e^{-i\alpha\xi^3t}\hat u(\sigma \xi,t)
=
\hat u_0(\sigma\xi)+F(\sigma\xi,t)
&+
\alpha
[\widetilde g_2(\xi^3,t)+i\sigma\xi\widetilde g_1(\xi^3,t)
-\sigma^2\xi^2\widetilde g_0(\xi^3,t)],
\quad
{\rm Im}(\sigma\xi)\leq 0,
\\
\label{glo-2}
e^{-i\alpha\xi^3t}\hat u(\sigma^2 \xi,t)
=
\hat u_0(\sigma^2 \xi)+F(\sigma^2 \xi,t)
&+
\alpha
[\widetilde g_2(\xi^3,t)+i\sigma^2 \xi\widetilde g_1(\xi^3,t)
-\sigma\xi^2\widetilde g_0(\xi^3,t)],
\quad
{\rm Im}(\sigma^2\xi)\leq 0.
\end{align}
Since $1+\sigma+\sigma^2=0$, 
multiplying equation \eqref{glo-1} by $\sigma$, 
equation \eqref{glo-2} by $\sigma^2$, 
and adding them to the equation $\alpha\widetilde g(\xi^3,t)
=
\alpha[\widetilde g_2(\xi^3,t)+i\xi\widetilde g_1(\xi^3,t)
-\xi^2\widetilde g_0(\xi^3,t)]$ we get the relation
\begin{align}
\label{global-solv}
\alpha
\widetilde g(\xi^3,t)
&=
\sigma[\hat u_0(\sigma\xi)+F(\sigma\xi,t)]+\sigma^2[\hat u_0(\sigma^2 \xi)+F(\sigma^2 \xi,t)]-3
\alpha
\xi^2\widetilde g_0(\xi^3,t)
\nonumber
\\
&-
\sigma e^{-i\alpha\xi^3t}\hat u(\sigma\xi,t)
-\sigma^2e^{-i\alpha\xi^3t}\hat u(\sigma^2 \xi,t),
\qquad
\xi\in \p D^+.
\end{align}
Substituting relation \eqref{global-solv} into the
 integral representation   \eqref{sln-defrom}
  and using the Cauchy Theorem, we remove the last two terms  and obtain  the first 
version of the UTM solution formula for ibvp \eqref{L-alpha-KdV}
\begin{align}
\label{kdv-utm-der}
&
u(x,t)
=
S_{\alpha}\big[u_0,g_0;f\big](x, t)
\doteq
\frac{1}{2\pi}\int_{-\infty}^\infty e^{i\xi x+
i\alpha\xi^3t
}
[\widehat {u}_0(\xi)+F(\xi,t)]d\xi
\\
+&
\frac{1}{2\pi}\int_{\p D^+} 
e^{i\xi x+
i\alpha\xi^3t
}
\big\{
\sigma[\widehat u_0(\sigma \xi)
+
F(\sigma \xi,t)]+\sigma^2[\widehat u_0(\sigma^2 \xi)+F(\sigma^2 \xi,t)]
-
3
\alpha
\xi^2 \tilde g_0(\xi^3,t)
\big\}d\xi.
\nonumber
\end{align}
Finally, using the Cauchy Theorem, again, we  replace  in 
formula \eqref{kdv-utm-der} the term
$
\int_{\p D^+} e^{i\xi x+i\alpha \xi^3t}\xi^2 \widetilde g_0(\xi^3,t)d\xi
$
by
$
\int_{\p D^+} e^{i\xi x+i\alpha\xi^3t}
$
$
\xi^2 \widetilde g_0(\xi^3,T)d\xi,
$
and get the desired UTM solution formula \eqref{kdv-utm-v},
after going back to the original ibvp notation.
\,\, $\square$

%
%
%
%
%
%
%
%
%
%
\vskip0.15in
\noindent
{\bf Acknowledgements.} The first author was partially supported by a grant from the Simons Foundation (\#524469 to Alex Himonas).

\vskip0.15in
\noindent
\noindent
{\bf Disclosure statement} 

\nin
The authors declare that they have no conflict of interest.

\vskip0.15in
\noindent
\noindent
{\bf Funding} 

\nin
Himonas was partially supported by
 Simons Foundation  grant  \#524469.

%

%
%
%
%
%
%
%
%
%
%
%

%

\vspace{3mm}
\noindent
A. Alexandrou Himonas  \hfill Fangchi Yan\\
Department of Mathematics  \hfill Department of Mathematics\\
University of Notre Dame  \hfill West Virginia University\\
Notre Dame, IN 46556  \hfill Morgantown, WV 26506 \\
E-mail: \textit{himonas.1$@$nd.edu}  \hfill E-mail: \textit{fyan1@alumni.nd.edu}

\end{document}